\documentclass{article}

\usepackage{graphicx} % Required for inserting images
\usepackage{subcaption}
\usepackage{verbatim}
\usepackage{csquotes}

\usepackage{tikz}
\usetikzlibrary{calc}
\usetikzlibrary{cd}
\usetikzlibrary{nfold}
\usetikzlibrary{arrows.meta, graphs, graphs.standard, quotes}
\usetikzlibrary{spath3}
\definecolor{Mycolor1}{HTML}{082A54}
\definecolor{Mycolor2}{HTML}{E02B35}
\definecolor{Mycolor3}{HTML}{F0C571}
\definecolor{Mycolor4}{HTML}{59A89C}
\definecolor{Mycolor5}{HTML}{A559AA}
\title{Polynomial Swap Coboundary Expansion with Applications to LTCs and PCPs}
\author{Yotam Dikstein \thanks{Tel Aviv University, Tel Aviv, Israel, yotam.dikstein@gmail.com} \orcidlink{0000-0002-6248-6574}, Izhar Oppenheim\thanks{Ben-Gurion University of the Negev, Be'er-Sheva, Israel.} \orcidlink{0000-0002-1849-1851}, Inga Valentiner-Branth \thanks{UC Berkeley, Berkeley, United States, ingavb@berkeley.edu} \orcidlink{0000-0001-7744-8733}}
\date{November 12, 2026}

\def\showauthornotes{0}

\def\showkeys{0}
\def\showdraftbox{0}
\def\showcolorlinks{1}
\def\usemicrotype{1}
\def\showfixme{0}

\usepackage[
backend=biber,
maxcitenames=50,
maxbibnames=99,
style=alphabetic,
sorting=ynt
]{biblatex}
\newcommand{\ordr}{34}

\providecommand{\NN}{\mathbb{N}}

\usepackage{nag}

\usepackage{xspace,enumerate}

\usepackage[dvipsnames]{xcolor}

\usepackage[T1]{fontenc}
\usepackage[full]{textcomp}

\usepackage{enumitem}

\makeatletter
\newcommand{\customlabel}[2]{#2\def\@currentlabel{#2}\label{#1}}
\makeatother
\makeatletter

\usepackage[american]{babel}

\usepackage{mathtools}

\usepackage{amsthm}
\usepackage{amsmath}

\newtheorem{theorem}{Theorem}[section]
\newtheorem*{theorem*}{Theorem}

\newtheorem{proposition}[theorem]{Proposition}
\newtheorem*{proposition*}{Proposition}
\newtheorem{lemma}[theorem]{Lemma}
\newtheorem*{lemma*}{Lemma}
\newtheorem{corollary}[theorem]{Corollary}
\newtheorem*{corollary*}{Corollary}
\newtheorem*{conjecture*}{Conjecture}
\newtheorem{fact}[theorem]{Fact}
\newtheorem*{fact*}{Fact}

\newtheorem*{hypothesis*}{Hypothesis}

\theoremstyle{definition}
\newtheorem{definition}[theorem]{Definition}
\newtheorem*{definition*}{Definition}

\newtheorem{example}[theorem]{Example}

\theoremstyle{remark}
\newtheorem{claim}[theorem]{Claim}
\newtheorem*{claim*}{Claim}
\newtheorem{remark}[theorem]{Remark}
\newtheorem*{remark*}{Remark}
\newtheorem{observation}[theorem]{Observation}
\newtheorem*{observation*}{Observation}

 \usepackage[letterpaper,
 top=1in,
 bottom=1in,
 left=1in,
 right=1in]{geometry}

\usepackage[varg]{pxfonts} % varg - uses nicer g,v,y,w

\usepackage[sc,osf]{mathpazo}
\usepackage{lmodern}

\ifnum\showkeys=1
\usepackage[color]{showkeys}
\fi

\ifnum\showcolorlinks=1
\usepackage[
colorlinks=true,
urlcolor=blue,
linkcolor=blue,
citecolor=OliveGreen,
]{hyperref}
\fi

\ifnum\showcolorlinks=0
\usepackage[
colorlinks=false,
pdfborder={0 0 0}
]{hyperref}
\fi

\usepackage{prettyref}

\newcommand{\savehyperref}[2]{\texorpdfstring{\hyperref[#1]{#2}}{#2}}

\newrefformat{eq}{\savehyperref{#1}{\textup{(\ref*{#1})}}}
\newrefformat{lem}{\savehyperref{#1}{Lemma~\ref*{#1}}}
\newrefformat{def}{\savehyperref{#1}{Definition~\ref*{#1}}}
\newrefformat{thm}{\savehyperref{#1}{Theorem~\ref*{#1}}}
\newrefformat{cor}{\savehyperref{#1}{Corollary~\ref*{#1}}}
\newrefformat{cha}{\savehyperref{#1}{Chapter~\ref*{#1}}}
\newrefformat{sec}{\savehyperref{#1}{Section~\ref*{#1}}}
\newrefformat{subsec}{\savehyperref{#1}{Subsection~\ref*{#1}}}
\newrefformat{app}{\savehyperref{#1}{Appendix~\ref*{#1}}}
\newrefformat{tab}{\savehyperref{#1}{Table~\ref*{#1}}}
\newrefformat{fig}{\savehyperref{#1}{Figure~\ref*{#1}}}
\newrefformat{hyp}{\savehyperref{#1}{Hypothesis~\ref*{#1}}}
\newrefformat{alg}{\savehyperref{#1}{Algorithm~\ref*{#1}}}
\newrefformat{rem}{\savehyperref{#1}{Remark~\ref*{#1}}}
\newrefformat{item}{\savehyperref{#1}{Item~\ref*{#1}}}
\newrefformat{obs}{\savehyperref{#1}{Observation~\ref*{#1}}}
\newrefformat{step}{\savehyperref{#1}{step~\ref*{#1}}}
\newrefformat{conj}{\savehyperref{#1}{Conjecture~\ref*{#1}}}
\newrefformat{fact}{\savehyperref{#1}{Fact~\ref*{#1}}}
\newrefformat{prop}{\savehyperref{#1}{Proposition~\ref*{#1}}}
\newrefformat{prob}{\savehyperref{#1}{Problem~\ref*{#1}}}
\newrefformat{claim}{\savehyperref{#1}{Claim~\ref*{#1}}}
\newrefformat{relax}{\savehyperref{#1}{Relaxation~\ref*{#1}}}
\newrefformat{red}{\savehyperref{#1}{Reduction~\ref*{#1}}}
\newrefformat{part}{\savehyperref{#1}{Part~\ref*{#1}}}
\newrefformat{ex}{\savehyperref{#1}{Example~\ref*{#1}}}
\newrefformat{para}{\savehyperref{#1}{Paragraph~\ref*{#1}}}
\newrefformat{ass}{\savehyperref{#1}{Assumption~\ref*{#1}}}
\newrefformat{question}{\savehyperref{#1}{Question~\ref*{#1}}}

\newcommand{\Sref}[1]{\hyperref[#1]{\S\ref*{#1}}}

\usepackage{nicefrac}
\ifnum\usemicrotype=1
\usepackage{microtype}
\fi

\ifnum\showauthornotes=1
\newcommand{\Authornote}[2]{{\sffamily\small\color{red}{[#1: #2]}}}
\newcommand{\Authornotecolored}[3]{{\sffamily\small\color{#1}{[#2: #3]}}}
\newcommand{\Authorcomment}[2]{{\sffamily\small\color{gray}{[#1: #2]}}}
\newcommand{\Authorstartcomment}[1]{\sffamily\small\color{gray}[#1: }

\newcommand{\Authorfnote}[2]{\footnote{\color{red}{#1: #2}}}
\newcommand{\Authorfixme}[1]{\Authornote{#1}{\textbf{??}}}
\newcommand{\Authormarginmark}[1]{\marginpar{\textcolor{red}{\fbox{\Large #1:!}}}}
\else
\newcommand{\Authornote}[2]{}
\newcommand{\Authornotecolored}[3]{}
\newcommand{\Authorcomment}[2]{}
\newcommand{\Authorstartcomment}[1]{}

\newcommand{\Authorfnote}[2]{}
\newcommand{\Authorfixme}[1]{}
\newcommand{\Authormarginmark}[1]{}
\fi

\newcommand{\Ydnote}{\Authornotecolored{Purple}{YD}}
\newcommand{\ydnote}{\Ydnote}
\newcommand{\Ionote}{\Authornotecolored{Red}{IO}}
\newcommand{\ionote}{\Ionote}
\newcommand{\Ivnote}{\Authornotecolored{Green}{IVB}}
\newcommand{\ivnote}{\Ivnote}

\ifnum\showfixme=0

\fi

\usepackage{boxedminipage}

\newcommand\sett[2]{\left\{ #1 \left| \; \vphantom{#1 #2} \right. #2  \right\}}
\newcommand{\set}[1]{\{#1\}}

\def\dim{\mathrm{ dim}}
\def\sp{\mathrm{ span}}

\newcommand{\F}{\mathbb{F}}

\newcommand{\Psymb}{\mathbb{P}}

\DeclareMathOperator*{\ProbOp}{\Psymb}

\renewcommand{\Pr}{\ProbOp}

\newcommand{\Prob}[2][]{\Pr_{{#1}}\left[#2\right]} % use by \Prob[x]{event}
 \usepackage{dsfont}
\usepackage{mathrsfs}

\newcommand{\textparen}[1]{\text{(#1)}}

\ifx\because\undefined
\newcommand{\because}[1]{\textparen{because #1}}
\else
\renewcommand{\because}[1]{\textparen{because #1}}
\fi

\newcommand\bdot\bullet

\DeclareMathOperator{\poly}{poly}

\DeclareMathOperator{\polylog}{polylog}
\DeclareMathOperator{\supp}{supp}
\DeclareMathOperator{\dist}{dist}

\DeclareMathOperator{\col}{col}
\DeclareMathOperator{\diam}{diam}

\newcommand{\N}{\mathbb N}
\newcommand{\R}{\mathbb R}

\newcommand{\cC}{\mathcal C}

\newcommand{\cI}{\mathcal I}

\renewcommand{\leq}{\leqslant}

\renewcommand{\geq}{\geqslant}

\ifnum\showdraftbox=1

\else

\fi

\let\epsilon=\varepsilon

\numberwithin{equation}{section}

\newcommand{\MYstore}[2]{%
  \global\expandafter \def \csname MYMEMORY #1 \endcsname{#2}%
}

\newcommand{\MYload}[1]{%
  \csname MYMEMORY #1 \endcsname%
}

\newcommand{\MYnewlabel}[1]{%
  \newcommand\MYcurrentlabel{#1}%
  \MYoldlabel{#1}%
}

\newcommand{\MYdummylabel}[1]{}

\newcommand{\torestate}[1]{%
  \let\MYoldlabel\label%
  \let\label\MYnewlabel%
  #1%
  \MYstore{\MYcurrentlabel}{#1}%
  \let\label\MYoldlabel%
}

\newcommand{\restatetheorem}[1]{%
  \let\MYoldlabel\label
  \let\label\MYdummylabel
  \begin{theorem*}[Restatement of \prettyref{#1}]
    \MYload{#1}
  \end{theorem*}
  \let\label\MYoldlabel
}

\newcommand{\restatelemma}[1]{%
  \let\MYoldlabel\label
  \let\label\MYdummylabel
  \begin{lemma*}[Restatement of \prettyref{#1}]
    \MYload{#1}
  \end{lemma*}
  \let\label\MYoldlabel
}

\newcommand{\restateprop}[1]{%
  \let\MYoldlabel\label
  \let\label\MYdummylabel
  \begin{proposition*}[Restatement of \prettyref{#1}]
    \MYload{#1}
  \end{proposition*}
  \let\label\MYoldlabel
}

\newcommand{\restateclaim}[1]{%
  \let\MYoldlabel\label
  \let\label\MYdummylabel
  \begin{claim*}[Restatement of \prettyref{#1}]
    \MYload{#1}
  \end{claim*}
  \let\label\MYoldlabel
}

\newcommand{\restatecorollary}[1]{%
  \let\MYoldlabel\label
  \let\label\MYdummylabel
  \begin{corollary*}[Restatement of \prettyref{#1}]
    \MYload{#1}
  \end{corollary*}
  \let\label\MYoldlabel
}

\newcommand{\restatefact}[1]{%
  \let\MYoldlabel\label
  \let\label\MYdummylabel
  \begin{fact*}[Restatement of \prettyref{#1}]
    \MYload{#1}
  \end{fact*}
  \let\label\MYoldlabel
}

\newcommand{\restatedefinition}[1]{% % overwrite label command with dummy
\let\MYoldlabel\label 
\let\label\MYdummylabel 
\begin{definition*}[Restatement of \prettyref{#1}] 
    \MYload{#1} 
\end{definition*} 
\let\label\MYoldlabel 
} 

\newcommand{\restate}[1]{%
  \let\MYoldlabel\label
  \let\label\MYdummylabel
  \MYload{#1}
  \let\label\MYoldlabel
}

\let\origparagraph\paragraph
\renewcommand{\paragraph}[1]{\origparagraph{#1.}}

\allowdisplaybreaks

\newcommand{\cclassmacro}[1]{\texorpdfstring{\textbf{#1}}{#1}\xspace}
\newcommand{\p}{\cclassmacro{P}}

\let\pref=\prettyref

\newcommand{\dir}[1]{\overset{\to}{#1}}
\newcommand{\bproof}[1]{\begin{proof}[Proof of \pref{#1}]}
\newcommand{\eproof}{\end{proof}}

\newcommand{\Img}{{\textrm {Im}}}

\newcommand{\Dehn}{\operatorname{Dehn}}
\newcommand{\Free}{\operatorname{Free}}
\newcommand{\Area}{\operatorname{Area}}
\newcommand{\Unip}{\operatorname{Unip}}
\newcommand{\elm}{\operatorname{elm}}
\newcommand{\Elm}{\operatorname{Elm}}

\def\S{X}

\def\T{T}

\newcommand{\CC}{\mathcal{CC}}

\newcommand{\KMS}[2]{\mathrm{KMS}(#1,#2)}
\DeclareMathOperator{\Sl}{SL} % I would prefer to use \SL but it seems to be used otherwise
\newcommand{\chev}{\operatorname{Chev}}
\newcommand{\flag}{\operatorname{Flag}}
\newcommand{\Ao}{\mathring \Phi}
\DeclareMathOperator{\stab}{Stab}
\newcommand{\ordX}{\overset{\rightarrow}{X}}
\DeclareMathOperator{\Aut}{Aut}
\newcommand{\TR}{\mathop{\sim}\limits^{(TR)}}
\newcommand{\BT}{\mathop{\sim}\limits^{(BT)}}
\DeclareMathOperator{\Rad}{Rad}
\newcommand{\Cone}{\operatorname{Cone}}

\usepackage{orcidlink}

\begin{document}

\maketitle

\begin{abstract}
    %Swap coboundary expansion, introduced independently by Bafna--Minzer and Dikstein--Dinur \cite{BafnaM2023,DiksteinD2023agr}, is a natural local-to-global property that underlies recent constructions of two-query PCPs with arbitrarily small soundness and near-linear proof length. Obtaining strong quantitative bounds for this notion has proved challenging: previous works established only subexponential lower bounds of \(\exp(-O(\sqrt{r}))\) for the swap coboundary expansion of Chapman--Lubotzky complexes \cite{DiksteinDL2024,BafnaLM2024}.

    Swap coboundary expansion, introduced independently by Bafna--Minzer and Dikstein--Dinur \cite{BafnaM2023,DiksteinD2023agr}, is a local-to-global property that underlies agreement tests in the low-acceptance regime and recent constructions of succinct PCPs. Previous works obtained bounds of \(\exp(-O(\sqrt{r}))\) for the \(r\)-swap coboundary expansion of Chapman--Lubotzky complexes \cite{DiksteinDL2024,BafnaLM2024}. We establish the first inverse-polynomial bounds for bounded-degree complexes: suitable skeleta of Kac--Moody--Steinberg (KMS) complexes have \(r\)-swap coboundary expansion at least \(r^{-O(1)}\).

    Our theorem supplies a key ingredient in forthcoming work of Ben Yaacov--Dikstein--Liu \cite{BenYaacovDL2026}, which obtains local tests for sparse algebraic codes with guarantees even at arbitrarily small acceptance probability. These tests use \(\exp(\poly(1/\varepsilon))\) queries and guarantee that any word accepted with probability at least \(\varepsilon\) agrees with a codeword on an \(\varepsilon^{O(1)}\) fraction of coordinates. Separately, substituting our bound into the existing analysis of Bafna--Minzer--Vyas--Yun \cite{BafnaMVY2025} improves the soundness dependence of their two-query PCP construction, giving perfect completeness, soundness \(\varepsilon\), and proof length \(n(\log n)^{\exp(\poly(1/\varepsilon))}\).

    Our proof exploits the product structure of the groups arising in KMS complexes. The main new ingredient is a decomposition lemma for Dehn functions, which separates the analysis of group relations into independent problems on much smaller factors, with only polynomial overhead.

    %Alternative to last paragraph:
    %Our proof exploits the algebraic structure of KMS complexes to turn a complicated word problem into a collection of elementary ones: by sorting group elements into commuting matrix blocks, we reduce the analysis to smaller groups with uniformly linear Dehn functions.
\end{abstract}
\section{Introduction}
High-dimensional expanders (HDXs) are generalizations of expander graphs into hypergraphs and other high dimensional structures. These objects support global structural guarantees from local constraints, and therefore are extremely useful in constructing locally testable and decodable codes, PCPs and other types of property tests \cite{DinurLZ2025,DiksteinHIP2026,BafnaMVY2025}. 

In this work, we study swap coboundary expansion, a high dimensional analogue of expansion that captures topological properties. We prove that families of Kac--Moody--Steinberg (KMS) complexes with arbitrarily many vertices and uniformly bounded degree are swap coboundary expanders with parameter \(\beta=r^{-\ordr}\), where \(r\) is the face dimension (an \(r\)-dimensional face has \(r+1\) vertices). This improves the previously known subexponential bounds for similar complexes to an inverse-polynomial bound \cite{DiksteinDL2024,BafnaLM2024,ODonnellS2025}.

%In this work, we study swap coboundary expansion, a high dimensional analogue of expansion that captures topological properties. We prove that families of Kac–Moody–Steinberg (KMS) complexes are swap coboundary expanders with parameter \(\beta=r^{-\ordr}\), improving the previously known subexponential bounds for similar complexes to an inverse-polynomial bound \cite{DiksteinDL2024,BafnaLM2024,ODonnellS2025}.

Very informally, coboundary expansion is a statement about when a labeling of directed edges of a simplicial complex can be explained by a labeling defined on the vertices (aka coboundaries). Given a labeling of directed edges by group elements, we test whether it comes from a vertex assignment by evaluating it on triangles of the complex: on each triangle, we multiply the values on its three directed edges. If the labeling is induced by a vertex assignment, then this product is always equal to the identity. Coboundary expansion is a robust converse to this statement — if most triangle evaluations yield the identity, then the labeling must be close to a coboundary.

Swap coboundary expansion is the analogous statement in a higher-order setting, where the objects being labeled are higher dimensional sets rather than vertices and edges.

This local-to-global mechanism is the core tool underlying recent constructions of agreement tests, locally testable codes, and probabilistically checkable proofs \cite{DinurLZ2025,BafnaMVY2025}. In these settings, one is given a collection of local functions defined on overlapping parts of a larger domain, and the goal is to “stitch” them into a single global function that simultaneously agrees with many of the local functions. Swap coboundary expansion provides the structural guarantee that makes this stitching possible.

\subsection{Definitions and our results}
We now formally define coboundary expansion. Let \(X\) be a simplicial complex with vertices \(X(0)\), edges \(X(1)\) and triangles \(X(2)\). Let \(\dir{X}(1)\) be the set of directed edges (i.e.\ all ordered pairs \(uv\) such that \(\set{u,v} \in X(1)\)). Let \(\mathrm{Sym}_m\) be the symmetric group on \(m\) elements, and \(\beta \in (0,1]\). For a function \(g:X(0) \to \mathrm{Sym}_m\) we denote by \(\delta g:\dir{X}(1) \to \mathrm{Sym}_m\) the function \(\delta g(uv)=g(v)g(u)^{-1}\). This function is called \emph{the coboundary of \(g\)}. A function \(f:\dir{X}(1) \to \mathrm{Sym}_m\) is called antisymmetric if \(f(uv)=f(vu)^{-1}\) for every \(uv \in \dir{X}(1)\). %Notice that \(\delta g\) is always antisymmetric. 
We say that \(X\) is a \((\beta,m)\)-coboundary expander if for every antisymmetric function \(f:\dir{X}(1) \to \mathrm{Sym}_m\) there exists \(g:X(0)\to \mathrm{Sym}_m\) such that
\begin{equation}\label{eq:intro-cbdry-exp}
    \beta \Prob[uv \in \dir{X}(1)]{f(uv) \ne g(v)g(u)^{-1}} \leq \Prob[uvw \in X(2)]{f(uv)f(wu)f(vw) \ne Id}.\footnote{The probabilities in this inequality are with respect to the natural probability distribution on \(X\). For the purpose of this introduction we encourage the readers to think of the uniform probability over directed edges or triangles respectively. See \pref{sec:prelims} for details.}
\end{equation}
In words, if an antisymmetric function fails the triangle identity on a \(\delta\)-fraction of triangles, then it can be changed on at most a \(\delta/\beta\)-fraction of edges to obtain a coboundary. This gives a local-to-global test: sample a random triangle \(uvw\) and check whether \(f(uv)f(wu)f(vw)=Id\). For swap coboundary expansion, our theorem makes the correction loss \(1/\beta\) polynomial in the dimension of the faces; this is the quantitative gain used in the applications below.
%
%In words, if an antisymmetric function multiplies to the identity on most triangles \(uvw\), then it agrees with a coboundary on almost all edges. This definition has a local-to-global flavor: to test closeness of \(f\) to the set of coboundaries, one samples a random triangle and checks how many triangles \(uvw\) have the property that \(f(uv)f(wu)f(vw) = Id\). %If \(f\) is given via some truth table, this test requires probing \(f\) at \(3\) positions only, and of course it does not require us to actually find any coboundary \(g\) that is close to \(f\). 
%This is considered a local-to-global definition as the local property of satisfying most triangles given the global conclusion of being close to some unknown coboundary.

As for swap coboundary expansion, let \(X\) be a \(d\)-dimensional simplicial complex. Let \(X(r)\) denote the faces of \(X\) of size \(r+1\). The \(r\)-faces complex, denoted \(F^r X\), is a complex whose vertices are \(F^r X(0) = X(r)\). Its \(k\)-faces are all \(t=\set{s_0,s_1,\dots,s_k}\) such that every pair of \(s_i,s_j \in t\) are disjoint and such that the union \(\bigcup_{i=0}^k s_i \in X\).

\begin{definition}[Swap coboundary expander]
    A simplicial complex \(X\) is an \((r,\beta,m)\)-swap coboundary expander if \(F^rX\) is a \((\beta,m)\)-coboundary expander.
\end{definition}

Our main theorem is about the KMS complexes, a variant of the complexes appearing in \cite{kaufmanO2018,GraveDePeraltaV2025,OppenheimV2025}.
\begin{comment} 
\begin{theorem}[Main] \label{thm:main}\ydnote{TODO: change after details are done.}
    For every \(r\) there exists \(q_0,d_0 = \poly (r)\) such that the following holds. Let \(q \geq q_0\) be a prime power and \(d\geq d_0\). Let \(X\) be the KMS complex of type \(A\) associated with dimension \(d\) and prime \(q\). Then \(X\) is a \((r,\Omega(r^{-10}),m)\)-swap coboundary expander for any \(m \leq \max \set{\exp(\log^{0.5} q), \mathrm{Char}(\F_q)}\).\footnote{The requirement that \(m \leq \exp(\log^{0.5} q)\) is required in the theorem so that the cohomology vanishes. When this does not hold, these complexes are swap \emph{cocycle} expanders, another important notion that has agreement testing applications.}

    In particular, there exists a sequence of \(d\)-dimensional complexes \(X\) such that the edge-degree of every vertex is \(\exp(O(d^2 \log q))\) that are simultaneously \(\frac{1}{q}\)-local spectral expanders and \((r,\Omega(r^{-10}),m)\)-swap coboundary expanders.
\end{theorem}
\end{comment}

\begin{theorem}[Main] \label{thm:main}
For every \(r\), there exist \(q_0,d_0=\poly(r)\) such that the following holds. Let \(q\geq q_0\) be a prime power and \(d\geq d_0\). Let \(X\) be a congruence KMS complex of type \(A_d\) over \(\F_q\). Then \(X\) is an \((r,\Omega(r^{-\ordr}),m)\)-swap coboundary expander for every integer \(1\leq m<\mathrm{Char}(\F_q)\).

In particular, if moreover \(q>d^2\), there exists a sequence of \(d\)-dimensional complexes whose number of vertices tends to infinity and whose vertex degrees are at most \(\exp(O(d^2\log q))\), all of which are simultaneously \(\frac{1}{\sqrt{q}-(d-1)}\)-local spectral expanders and \((r,\Omega(r^{-\ordr}),m)\)-swap coboundary expanders for every integer \(1\leq m<\mathrm{Char}(\F_q)\).
\end{theorem}
%As one might expect, we only use some abstract properties of these KMS complexes. We give a more abstract version of this theorem in \ydnote{Decide whether we actually do this?}.

\subsection{Applications of our theorem}
\ydnote{I switched application order. See comment below.}
Both our applications rely on the same mechanism: using swap coboundary expansion to convert local consistency into global structure.

\paragraph{Derandomized low-degree testing} The Reed-Muller error correcting code \(RM(m,d,q)\) consists of all (truth tables of) degree \(d\) polynomials over \(\F_q^m\). This code plays an important role in constructions of PCPs and of locally testable codes \cite{KiwiMS2002}. In particular, locally testing whether a function is (close to a) low-degree polynomial by only querying few values from its truth table is a problem that has been extensively studied \cite{RubinfeldS1996,RazS1997,AroraS2003,MoshkovitzR2008,HarshaKSS2024}. Such a test is called a \emph{low degree test}.

One of the most common tests is the point vs.\ subspace test. The tester in this test receives as input a function \(f:\F_q^m \to \F_q\) which we want to test and determine whether it is a low degree function. In addition, it receives an additional `witness' \(F\): for every \(k\)-dimensional subspace (for some constant parameter \(k\)) \(W \subseteq \F_q^m\), \(F\) has a degree \(d\) polynomial \(F_W:W \to \F_q\). The tester samples a random \(W \subseteq \F_q^m\) and a point \(x \in W\) and accepts if \(F_W(x)=f(x)\). 

Recently, \cite{DinurLZ2025} analyzed a derandomized version of this test, where both \(f\) and the witness \(F\) are only defined on a small subset of the points \(S\), and a small subset of all possible subspaces. As a result they get a low degree test whose witness \(F\) has linear size in \(f\) and has the following guarantee. For every function \(f:S \to \F_q\), if 
\[\exists F \; \; \Prob[x,W]{(f,F) \text{ passes the test}} \geq 1-\varepsilon,\]
then there exists a codeword \(g:S\to\F_q\), whose restriction to every prescribed local domain has low degree, such that
\[\dist(f,g) \leq O(\varepsilon).\]
Here the distance is the normalized Hamming distance over \(S\).

A very natural question is what happens in the high error regime, i.e.\ what happens if 
\[\exists F \; \; \Prob[x,W]{(f,F) \text{ passes the test}} \geq \varepsilon\]
where \(\varepsilon > 0\) is just a small constant (or even subconstant). The analysis of \cite{DinurLZ2025} does not generalize to the this setup. In particular, if the probability of passing the test is less than \(\frac{4}{5}\), the analysis in \cite{DinurLZ2025} does not predict any conclusion). However, in upcoming work, \cite{BenYaacovDL2026} show that a modified version of the code in \cite{DinurLZ2025} still has a non-trivial guarantee in this case: passing the test with probability \(\varepsilon\) guarantees that there exists a codeword \(g:S\to\F_q\) such that
\[\dist(f,g)\leq 1-\poly(\varepsilon).\]

The proof of \cite{BenYaacovDL2026} uses swap coboundary expansion of the KMS complexes underlying the code to stitch local corrections into a single codeword. Our theorem establishes this property for these complexes, for which even an inverse-exponential swap coboundary expansion bound was not previously known.

The polynomial strength of our bound also improves the query complexity of the resulting local tester. This tester reads all the values of \(f\) on a sampled local domain and checks the prescribed low-degree constraints there. Our bound yields \(\exp(\poly(1/\varepsilon))\) queries. If one instead had an expansion bound of the previously known subexponential strength, the same analysis would yield \(\exp(\exp(\poly(1/\varepsilon)))\) queries. Thus our bound gives a testing guarantee at smaller acceptance thresholds for a given query budget.

\paragraph{Succinct PCPs}
A probabilistically checkable proof (PCP) allows a randomized polynomial-time verifier to check a proof by reading only a few of its symbols. The number of symbols read is the query complexity, while the total number of symbols in the proof is its length. For satisfiability, perfect completeness means that every satisfiable instance has a proof that is always accepted; soundness error \(\varepsilon\) means that, for every unsatisfiable instance, every purported proof is accepted with probability at most \(\varepsilon\).

A central goal in PCP constructions is to simultaneously achieve:
\begin{enumerate}
    \item A constant number of queries (ideally \(2\) queries),
    \item Arbitrarily small soundness error \(\varepsilon>0\) (even subconstant error), and
    \item Near-linear proof length.
\end{enumerate}
Achieving all three simultaneously was a longstanding challenge. A breakthrough of \cite{BafnaMVY2025} achieves this combination for every fixed \(\varepsilon>0\): they construct a two-query PCP with soundness error \(\varepsilon\) and proof length
\begin{equation}\label{eq:intro-proof-length}
    |w|=O(|x|\log^{C(\varepsilon)}|x|),
\end{equation}
where \(|x|\) is the input size and \(|w|\) is the proof length.

More recently, \cite{DiksteinDHS2026} obtained two-query PCPs with proof length \(n(\log n)^{O(1)}\), where the exponent is independent of \(\varepsilon\), for soundness down to \(1/\polylog(n)\). Their construction uses HDXs of constant dimension, independent of the soundness parameter, and does not require bounding swap coboundary expansion.

Although this newer construction achieves better parameters, our theorem also improves the parameters of the original construction of \cite{BafnaMVY2025}. This improvement follows simply by substituting our expansion bound into their existing analysis: it requires no new PCP analysis and does not use the new arguments of \cite{DiksteinDHS2026}.

To explain this dependence, the core technical ingredient in \cite{BafnaMVY2025} is a family of agreement tests analyzed by \cite{DiksteinD2023agr,BafnaM2023,DiksteinD2023swap,DiksteinDL2024,BafnaLM2024}. These tests check whether a collection of local functions on overlapping domains is consistent with a single global function. Their quantitative guarantees depend on swap coboundary expansion, and this dependence propagates to the proof-length exponent \(C(\varepsilon)\). The previously available expansion bounds give \(C(\varepsilon)=\exp(\exp(\poly(1/\varepsilon)))\). Our polynomial bound reduces this to \(C(\varepsilon)=\exp(\poly(1/\varepsilon))\), removing one exponential from the dependence of the polylogarithmic overhead on the soundness parameter.

\begin{theorem}[Improved parameters for the PCP of \cite{BafnaMVY2025}]\torestate{\label{thm:pcp}
    There exists an absolute constant \(a>0\) such that, for every fixed \(\varepsilon\in(0,1/2)\) and all sufficiently large \(n\), the PCP construction of \cite{BafnaMVY2025}, instantiated using \pref{thm:main}, gives a polynomial-time reduction from a \(3\)-SAT instance \(\varphi\) of size \(n\) to a two-query PCP with perfect completeness, soundness at most \(\varepsilon\), and proof length at most \[n(\log n)^{\exp(O(\varepsilon^{-a}))}.\] The proof alphabet has size at most \(\exp\left((1/\varepsilon)^{O(\log(1/\varepsilon))}\right)\).}
\end{theorem}

\subsection{Proof overview}
To give an overview of the proof of \pref{thm:main}, we present an analysis of a family of complexes that are subcomplexes of links of the KMS complex. We dub these `matrix block complexes'. We will discuss how to bound their coboundary expansion. In the end of the overview, we will explain how a series of reductions (all of which appear in previous works) lift coboundary expansion lower bounds on these complexes to coboundary expansion lower bounds on the KMS complexes.

The key idea is that the groups defining the matrix block complexes possess a natural product structure. Our main new ingredient is a decomposition lemma that allows us to analyze the factors separately, with only a polynomial overhead in \(r\). This reduces the problem of proving coboundary expansion for an \(r\)-faces complex to independent problems on smaller groups, whose defining subgroups have only \(O(\log r)\) diagonal blocks.

What does this reduction buy us? The block boundaries are randomly positioned. For example, if all boundaries of a subgroup lie in the first fifth of the matrix, then the remaining four fifths form a single large block. The probability of this event is exponentially small in the number of boundaries. With only \(O(\log r)\) boundaries, it is therefore inverse polynomial in \(r\), so polynomially many sampled subgroups suffice to find such a subgroup. Applying the same idea to several prescribed regions gives the arrangement of large blocks needed for our direct analysis. These large blocks let us handle whole regions of the matrix together, using an argument in the spirit of \cite{KaufmanO2021} whose bounds are independent of the matrix size. Combining these bounds using the decomposition lemma yields the desired inverse-polynomial expansion bound.

%The key idea is that the groups defining the matrix block complexes possess a natural product structure. We show that this algebraic structure can be exploited to decompose the problem of proving coboundary expansion for an \(r\)-faces complex into many simpler coboundary-expansion problems on face complexes indexed by \(r'=O(\log r)\)-faces rather than \(r\)-faces. We then analyze these smaller instances directly in the spirit of \cite{KaufmanO2021}. This analysis requires an ambient complex whose dimension is exponential in \(r'\), and yields a coboundary-expansion bound that is inverse polynomial in this dimension, or equivalently exponential in \(r'\). Since \(r'=O(\log r)\), this translates into the desired inverse polynomial bound in r.

%The resulting complexes admit a particularly simple geometric structure, allowing for a direct analysis à la \cite{KaufmanO2021}. While this direct analysis yields a mere exponential lower bound in the structural parameter, the exponent is in \(O(\log r)\) due to the reduction, and therefore polynomial in \(r\). \ionote{I don't think this is what we ended up doing. The Dehn bound for each upper triangular block is linear with a universal constant. The $r$ actually comes into effect in two ways - the sorting in the direct sums and the size of $T$ such that we get niceness} \ydnote{How is this?}

We begin by introducing the matrix block complexes and recalling the connection between coboundary expansion and Dehn functions for coset complexes. We then explain how the product structure leads to the decomposition described above, and how this reduction exposes the simple geometric structure needed for the analysis.

Finally, we give a disclaimer. To explain our argument clearly, we skip the definition of the KMS complexes in this overview. Instead, we will explain how to bound the swap coboundary expansion for \emph{a vertex link} of the KMS complex. The same decomposition and Dehn-function arguments apply to the full KMS complexes; the additional vanishing-cohomology input is discussed at the end of the overview.

\subsubsection{Setup: matrix block complexes and the Dehn function}
The matrix block complexes are coset complexes.
\begin{definition}[Coset complexes]
Let \(\Gamma\) be a group and \(H_0,H_1,\dots,H_d\leq\Gamma\) be distinct subgroups. The coset complex \(X=\CC(\Gamma;H_0,H_1,\dots,H_d)\) is the complex whose vertices are cosets
\[X(0)=\sett{AH_i}{0\leq i\leq d,\; A\in\Gamma}.\]
Its \(k\)-faces are
\[X(k)=\sett{\set{AH_{j_0},AH_{j_1},\dots,AH_{j_k}}}{0\leq j_0<j_1<\dots<j_k\leq d,\; A\in\Gamma}.\]
\end{definition}
Coset complexes have been studied in the context of high dimensional expanders in many works such as \cite{kaufmanO2018,KaufmanO2021,ODonnellP2022,GraveDePeraltaV2025,KaufmanOW2025,OppenheimV2025, ODonnellS2025}. %For now we just state a simple claim that we will need to describe the links of these complexes.
% \begin{claim} \label{claim:links-of-coset-complex}
%     Let \(X\) be a coset complex and let \(A H_1 \in X(0)\) be a vertex. Then the link \(X_{A H_1}\) is isomorphic to the coset complex \(\CC(H_1; H_1 \cap H_2, H_1 \cap H_3, \dots ,H_1 \cap H_d)\).
% \end{claim}
% This gives us an easy way to describe links of a coset complex. This is important because the local-to-global machinery described later allows us to convert coboundary expansion of links (which are often simpler in some sense) to coboundary expansion of \(X\) itself.

Fix a dimension \(d>0\) and a field \(\F_q\). Let \(\mathrm{Unip}_{d+1}\subseteq SL_{d+1}(\F_q)\) be the subgroup consisting of upper triangular unipotent matrices. A vertex link of the \(d\)-dimensional KMS complex is isomorphic to the coset complex \(\CC(\mathrm{Unip}_{d+1};H_1,H_2,\dots,H_d)\), where, for \(1\leq i\leq d\),
\[H_i=\sett{\begin{pmatrix}
A_1&0\\
0&A_2
\end{pmatrix}}{A_1\in\mathrm{Unip}_i,\; A_2\in\mathrm{Unip}_{d+1-i}}.\]

The swap complex of this complex is not a coset complex, but nevertheless it has group-inherited geometry. The vertices of the \(r\)-swap complex of this complex can be identified with cosets of the subgroups \(H_J= \cap_{i \in J}H_i\), for all subsets \(J \in \binom{\{1,\dots,d\}}{r+1}\). The faces in the \(r\)-swap complex are all \(\set{A H_{J_0}, A H_{J_1}, \dots, A H_{J_\ell}}\) for \(A \in \mathrm{Unip}_{d+1}\) and \emph{pairwise disjoint} index sets \(J_0,J_1,\dots,J_\ell \in \binom{\{1,\dots,d\}}{r+1}\). The pairwise disjointness is the reason this is not a coset complex per se. Albeit, if \(J_0,J_1,\dots,J_T\) are pairwise disjoint, the complex induced by the cosets of \(H_{J_0},H_{J_1},\dots,H_{J_T}\) is indeed a coset complex. In a nutshell, the matrix block complexes are these subcomplexes. A reduction in \cite{DiksteinD2023swap} shows that, assuming sufficiently strong local spectral expansion, if most collections of pairwise disjoint color sets \(J_1,\dots,J_T\) induce coboundary expanders, then the whole faces complex is a coboundary expander.

Therefore, the main object in this proof overview are actually these subcomplexes, which we call matrix block complexes and denote by \(Y^{\cI}\). Every subcomplex is indexed by collections \(\cI = \set{J_1,J_2,\dots,J_T} \subset \binom{\{1,\dots,d\}}{r+1}\). That is, \(\cI\) is a collection of \(T=\poly(r)\) subsets of \(\{1,\dots,d\}\), each of size \(|J_i|=r+1\) which are pairwise disjoint. 

Let us describe these subgroups explicitly. For \(J=\set{i_1<\dots<i_{r+1}}\), set \(i_0=0\) and \(i_{r+2}=d+1\). Then
\[H_J=\sett{\begin{pmatrix}
A_1&0&\dots&0\\
0&A_2&\dots&0\\
0&0&\ddots&0\\
0&0&\dots&A_{r+2}
\end{pmatrix}}{A_\ell\in\mathrm{Unip}_{i_\ell-i_{\ell-1}}\text{ for every }1\leq\ell\leq r+2}.\]
In other words, we can decompose
\[H_J=\mathrm{Unip}_{i_1}\oplus\mathrm{Unip}_{i_2-i_1}\oplus\dots\oplus\mathrm{Unip}_{d+1-i_{r+1}}\]
as diagonal blocks.

The matrix block complexes are \(Y^{\cI} \coloneqq \CC(\mathrm{Unip}_{d+1}; \set{H_J}_{J \in \cI})\). The main technical result of the paper is that most \(Y^{\cI}\) are coboundary expanders.
\begin{lemma}[Informal, combining \pref{lem:most-colors-are-nice} and \pref{prop: nice colored faces are CBE}] \label{lem:main-technical}
For suitable polynomial choices of \(T\) and \(d\) in \(r\), and under the field-size and coefficient assumptions of \pref{thm:main}, the following holds. With probability at least \(90\%\) over a uniformly random collection \(\cI\) of \(T\) mutually disjoint \((r+1)\)-element subsets of \(\set{1,\dots,d}\), the complex \(Y^{\cI}\) is an \((r^{-O(1)},m)\)-coboundary expander.
\end{lemma}

To prove \pref{lem:main-technical} we use the group theoretic properties of this complex. More specifically, we use the connection between coboundary expansion of a coset complex and its associated Dehn function (or `word problem'), as observed by \cite{KaufmanO2021}. 
Let \(X=\CC (\Gamma;H_0,H_1,\dots,H_d)\) be as above. 
Let 
    \[g_1 \cdot g_2 \cdot \dots \cdot g_\ell = 1_{\Gamma}\]
where every \(g_j \in \bigcup_{i=0}^d H_i\).\footnote{Formally, the object we are looking at is a finite sequence of group elements \((g_1,g_2,\dots,g_\ell)\) whose product, taken in the order above, is equal to the identity.} 
The area of the product is the minimal number of `reduction steps' necessary to go from this product to the trivial product \(1_{\Gamma}\), where a reduction step is one of the following
\begin{enumerate}
    \item Add or remove \(g g^{-1}\) for some \(g \in \bigcup_{i=0}^d H_i\) somewhere in the sequence.
    \item Replace the pair \(g_j g_{j+1}\) with \(g'\) such that \(g_j g_{j+1} = g'\) and all three \(g_j,g_{j+1},g'\) are in the same subgroup \(H_i\) for some \(i=0,1,\dots,d\).
\end{enumerate}
The Dehn function associated with the complex \(X\), 
\[\mathrm{Dehn}_X(t) = \max \sett{area(g_1\cdot  g_2 \cdot \dots \cdot g_\ell)}{\ell\leq t, \; \;g_1,g_2,\dots,g_\ell \in \bigcup_{i=0}^d H_i, \; \; g_1 \cdot g_2 \cdot \dots \cdot g_\ell = 1_\Gamma}.\]
In words, this is the maximal number of reduction steps necessary in order to reduce a word of length at most \(t\) to the trivial product (and define \(\mathrm{Dehn}_X(t) = \infty\) if no such reduction exists). For a formal definition see \pref{subsec:dehn}.

A theorem of \cite{KaufmanO2021} states that if \(X\) has a small diameter and if every short word can be reduced to the identity in a few steps, then the coset complex is a coboundary expander.
\begin{theorem}[{\cite{KaufmanO2021}}] \label{thm:dehn-informal}
Let \(X\) be a \(d\)-dimensional coset complex, with diameter \(k\) and \(\mathrm{Dehn}_X(2k+1)<\infty\). Then \(X\) is a \((\poly(\frac{1}{d},\frac{1}{k},\frac{1}{\mathrm{Dehn}_X(1+2k)}),m)\)-coboundary expander for every \(m > 1\).
\end{theorem}
Let us denote by \(\Delta_X=\mathrm{Dehn}_X(1+2\mathrm{diam}(X))\). That is, the smallest integer such that every product of \(1+2\mathrm{diam}(X)\) elements as above that are equal to the identity, can be reduced to the identity in \(\Delta_X\) reduction steps. Next we see how the Dehn function interplays with the product structure of the links of \(Y^{\cI}\). %We record for later a trivial bound that states that the Dehn function of \(X\) is upper bounded by a linear function whose slope is \((1+2\mathrm{diam}(X)))\Delta_X\):
% \begin{equation}\label{eq:dehn-bound-with-diameter}
%     \mathrm{Dehn}_X(t) \leq (1+2\mathrm{diam}(X)))\Delta_X \cdot t.
% \end{equation}

\subsubsection{The product decomposition}
Once sufficiently strong local spectral expansion and vanishing first cohomology have been established, the local-to-global machinery of \cite{KaufmanKL2016,EvraK2016,DiksteinD2023cbdry} reduces the remaining task to proving coboundary expansion of every vertex link of \(Y^{\cI}\). We will therefore analyze the Dehn functions of these links. Specifically, we will apply \pref{thm:dehn-informal} with both the diameter and \(\Delta_{Y^{\cI}_v}\) bounded by \(\poly(r)\). In this overview we focus on the Dehn bound; the diameter bound follows from a simpler argument.

Fix \(J_1\in\cI\). The link of a vertex \(AH_{J_1}\) is isomorphic to the coset complex with ambient group \(H_{J_1}\) and defining subgroups \(H_{J_1}\cap H_{J'}\), where \(J'\) ranges over \(\cI\setminus\set{J_1}\). The key structural observation is that the ambient group and all these subgroups decompose along the same diagonal blocks. We use this common product structure to bound the Dehn function of the link.

%The key structural observation is that for any color sets \(J_1,J'\), the groups \(H_{J_1}\) and \(H_{J_1} \cap H_{J'}\) have a group product structure. We can use this product structure to get control on the Dehn function. %\ivnote{What are $J_1,J'$ here? Not sure if used before, but in any case I think we should recall it here.}

\begin{figure}[ht]
\centering
\begin{tikzpicture}[
    blk/.style  = {draw, fill=gray!25},
    dblk/.style = {draw, fill=gray!45},
    grp/.style  = {draw, dashed},
    lbl/.style  = {inner sep=1pt}
  ]

  %%% left matrix: H_{J_1}
  \draw[blk] (0,0)       rectangle node{$G_1$}     (1.30,-1.30);
  \draw[blk] (1.30,-1.30) rectangle node{$G_2$}     (2.57,-2.57);
  \draw[blk] (3.35,-3.35) rectangle node{$G_{r+2}$} (4.97,-4.80);
  \foreach \t in {2.77,2.96,3.15} \fill (\t,-\t) circle (1.3pt);
  \node at (-0.45,-2.44) {$\left(\rule[-2.65cm]{0pt}{5.3cm}\right.$};
  \node at ( 5.42,-2.44) {$\left.\rule[-2.65cm]{0pt}{5.3cm}\right)$};
  \node at (2.49,-5.55) {$H_{J_1}$};

  %%% right matrix: H_{J_1} \cap H_{J'}
  \begin{scope}[xshift=7.2cm]
    % group boxes first, so the solid block outlines are drawn over them;
    % they stop at the G_1/G_2 tangency point (1.15) instead of crossing there
    \draw[grp] (-0.08,0.08) rectangle (1.15,-1.15);
    \draw[grp] (1.15,-1.15) rectangle (2.65,-2.65);
    \draw[grp] (3.27,-3.27) rectangle (4.88,-4.88);
    \node[lbl,anchor=south west] at (1.15,-1.15) {\small$G_1$};
    \node[lbl,anchor=south west] at (2.65,-2.65) {\small$G_2$};
    \node[lbl,anchor=south east] at (4.72,-3.27) {\small$G_{r+2}$};

    \draw[dblk] (0,0)        rectangle node{\small$G_{1,1}$}       (0.70,-0.70);
    \draw[blk]  (0.70,-0.70) rectangle node{\tiny$G_{1,2}$}        (1.15,-1.15);
    \draw[blk]  (1.15,-1.15) rectangle node{\tiny$G_{2,1}$}        (1.57,-1.57);
    \draw[dblk] (1.57,-1.57) rectangle node{\scriptsize$G_{2,2}$}  (2.15,-2.15);
    \draw[blk]  (2.15,-2.15) rectangle node{\tiny$G_{2,3}$}        (2.57,-2.57);
    \draw[dblk] (3.35,-3.35) rectangle node{\tiny$G_{r+2,1}$}      (3.97,-3.97);
    \draw[blk]  (3.97,-3.97) rectangle node{\scriptsize$G_{r+2,2}$}(4.80,-4.80);

    \foreach \t in {2.81,2.96,3.11} \fill (\t,-\t) circle (1.3pt);
    \node at (-0.55,-2.44) {$\left(\rule[-2.65cm]{0pt}{5.3cm}\right.$};
    \node at ( 5.35,-2.44) {$\left.\rule[-2.65cm]{0pt}{5.3cm}\right)$};
    \node at (2.40,-5.55) {$H_{J_1}\cap H_{J'}$};
  \end{scope}
\end{tikzpicture}
\caption{The blocks}
\label{fig:blocks}
\end{figure}
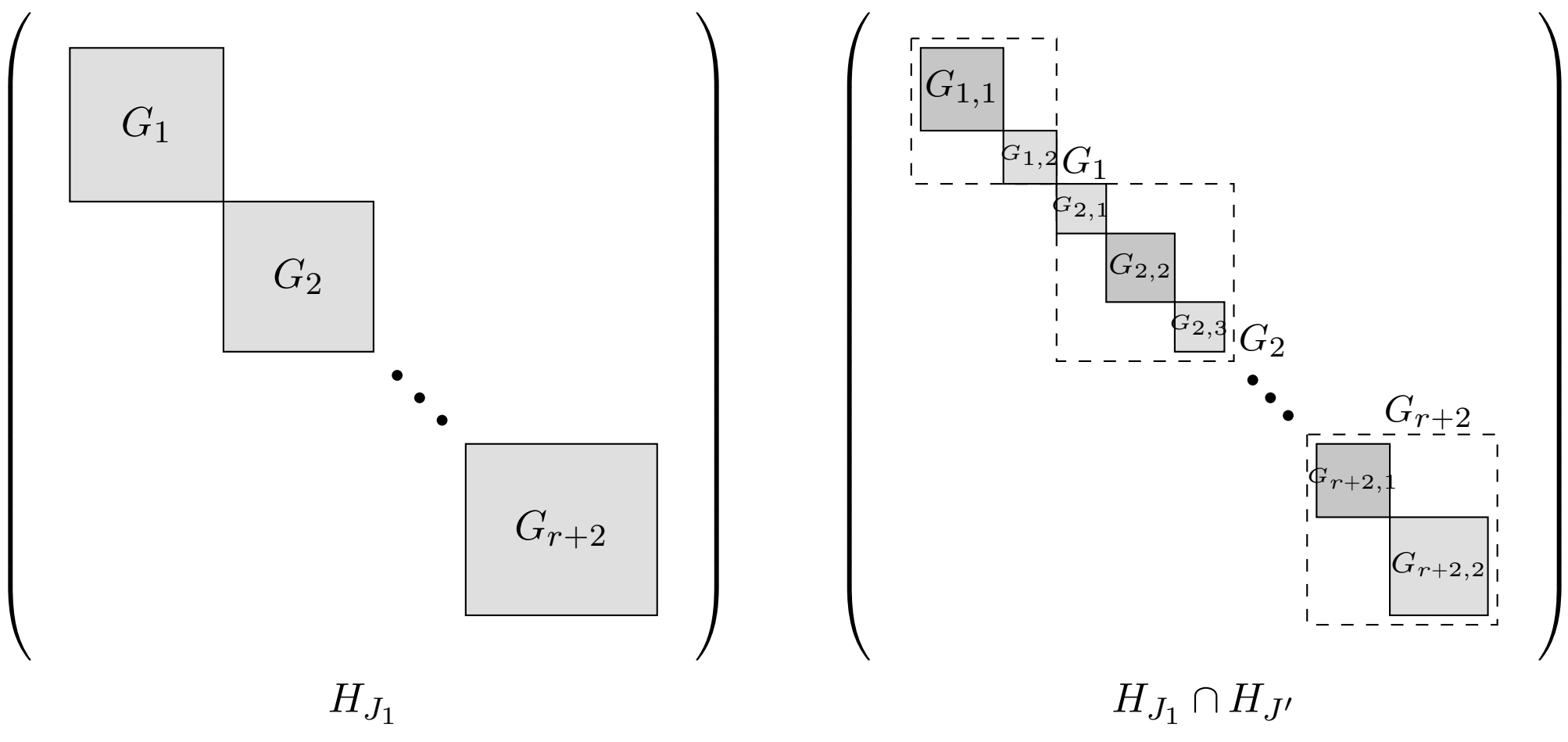

Indeed, recall that \[H_{J_1} = \prod_{\ell=1}^{r+2} G_\ell\] 
where every \(G_\ell \cong \mathrm{Unip}_{j_\ell}\) for some integer \(j_\ell\). These correspond to the different blocks in the matrices in \(H_{J_1}\). The same holds for \(H_{J_1} \cap H_{J'}\) which is equal to \( H_{J_1 \cup J'}\). One observes that every block in \(H_{J_1} \cap H_{J'}\) sits entirely inside one of the blocks in \(H_{J_1}\). For a pictorial example, see \pref{fig:blocks}.

% Requires: \usepackage{tikz}
% Blocks sit on the diagonal: a block spanning [a,b] is the square (a,-a)--(b,-b),
% so consecutive blocks are tangent exactly when one's b equals the next one's a.

Therefore, we can decompose \(H_{J_1}\cap H_{J'}=\prod_{\ell=1}^{r+2}H_\ell(J')\), where each \(H_\ell(J')\leq G_\ell\) is itself a direct product of upper triangular unipotent groups.

Let \(Y^0 = \CC(H_{J_1};\set{H_{J_1 \cup J'}}_{J' \in \cI \setminus \set{J_1}})\) be the link of any coset \(A H_{J_1}\) and let \(Y_\ell = \CC(G_\ell;\set{H_\ell(J')}_{J' \in \cI \setminus \set{J_1}})\) (these depend on \(\cI\) even though it doesn't appear in the notation).

The heart of the proof is the following bound on \(\Delta_{Y^0}\) using the Dehn functions of the individual factors \(\mathrm{Dehn}_{Y_\ell}\) on words whose length is \(1+2\mathrm{diam}(Y^0)\). More explicitly we show that with high probability over the choice of \(\cI\),
\begin{equation} \label{eq:sum-of-dehn}
    \Delta_{Y^{0}} \leq \sum_{\ell=1}^{r+2}\mathrm{Dehn}_{Y_{\ell}}(1+2\mathrm{diam}(Y^0)) + \poly(r).
\end{equation}
Thus, instead of analyzing the Dehn function of the original link directly, it suffices to analyze the Dehn function for much simpler complexes \(Y_\ell\).

To prove \eqref{eq:sum-of-dehn}, we sort a word of length at most \(1+2\diam(Y^0)\) into separate words, one for each factor \(G_\ell\). Although elements from distinct factors commute, exchanging them must be justified using the allowed reduction steps, which operate within a single defining subgroup. For example, suppose \(ghk=Id\) and write \(g=g_1g_2g_3\), \(h=h_1h_2h_3\), and \(k=k_1k_2k_3\), where \(g_i,h_i,k_i\in G_i\). Each letter can be split into its components using reductions within the subgroup containing that letter. We then transform \(g_1g_2g_3h_1h_2h_3k_1k_2k_3\) into the `sorted product' \((g_1h_1k_1)(g_2h_2k_2)(g_3h_3k_3)\). Since \(ghk=Id\), each component word equals the identity and can be reduced separately using the Dehn bound for \(Y_i\). In general there are \(r+2\) components, but the same argument applies.

We perform the sorting by a bubble-sort algorithm, repeatedly exchanging adjacent elements from different factors. With high probability over \(\cI\), each exchange can be implemented by splitting the two elements into a bounded number of pieces. Each pair of pieces that we need to exchange lies in a common subgroup \(H_{J_1}\cap H_{J'}\), for a suitable \(J'\in\cI\setminus\set{J_1}\), and still belongs to distinct factors. We can therefore exchange the pieces using reductions within that subgroup, and then recombine them. Since the original word length and the number of components are polynomial in \(r\), the total sorting cost is \(\poly(r)\). Each resulting component word has length at most the original word length, so reducing these words separately gives \eqref{eq:sum-of-dehn}.

\subsubsection{The direct analysis}
What do we gain from \eqref{eq:sum-of-dehn}? The gain is that the complexes \(Y_\ell\) are much simpler than \(Y^0\). They have a similar structure to \(Y^0\), being coset complexes defined by block-unipotent matrix subgroups. However, although the number of subgroups remains \(T-1=\poly(r)\), a random choice of \(\cI\) ensures with high probability that each subgroup \(H_\ell(J')\) has only \(O(\log r)\) blocks. More specifically, if \(G_\ell\cong\mathrm{Unip}_{d'}\), then with high probability \(d'=O(d\log r/r)\), and each \(J'\) places only \(O(\log r)\) indexes inside this block. A probabilistic argument then shows that, for every \(\ell\), there are three subgroups \(H_\ell(J'),H_\ell(J''),H_\ell(J''')\) whose indexes are confined to the first, last, and middle fifth of the block, respectively. Up to rounding, the first two subgroups contain a bottom block and a top block of size at least \(4d'/5\), respectively, while the third contains both a top block and a bottom block of size at least \(2d'/5\). This is an advantage because we can now use any relation that lives in these huge blocks to reduce our words.

Our goal is now to reduce any word of length \(t\) whose product is the identity using \(O(t)\) allowed reduction steps, independently of the matrix size. The huge blocks are useful because one allowed step can perform a calculation on many matrix entries at once: whenever two neighboring letters belong to the same defining subgroup, we may replace them by their product. For example, consider two matrices that differ from the identity only within the same rectangular block above the diagonal, containing arrays \(A\) and \(B\). Their product simply has \(A+B\) in that block. If the rectangle lies inside a huge diagonal block of one of our subgroups, this multiplication is one allowed reduction step, regardless of how many entries \(A\) and \(B\) contain.

We exploit this by viewing each matrix as a \(5\times5\) array of blocks of size approximately \(d'/5\), and treating entire blocks as entries. The example above illustrates how addition of these entries can be performed in one step. More generally, we show that the arrangement and overlap of the three large-block subgroups provide the relations needed for the block-matrix calculations in our reduction, each at a bounded cost. Thus, enlarging the matrices only enlarges the arrays involved in these calculations; it does not increase the number of reduction steps needed to perform them. Using these relations, we obtain
\begin{equation}\label{eq:constant-slope-in-block-dehn}
    \Dehn_{Y_\ell}(t)\leq Ct
\end{equation}
for a universal constant \(C>0\). This argument is in the spirit of \cite{KaufmanO2021}.

 We stress that in the more common case where there are many small blocks, even proving a bound on the Dehn function is more challenging (partially from the fact that even the diameter in this case is super constant), but we completely sidestep this issue by using the big-block triple of subgroups.

Therefore, together with a bound on the diameter of \(Y^0\) we get a bound of
\[\Delta_{Y^0} \overset{\eqref{eq:sum-of-dehn}, \eqref{eq:constant-slope-in-block-dehn}}{\leq} C(r+2)(1+2\mathrm{diam}(Y^0)) + \poly(r) = \poly(r).\]

\subsubsection{Putting everything together}
The argument above bounds the coboundary expansion of every vertex link of most subcomplexes \(Y^{\cI}\). Together with sufficiently strong local spectral expansion and vanishing first cohomology, the local-to-global machinery of \cite{KaufmanKL2016,EvraK2016,DiksteinD2023cbdry} implies that these \(Y^{\cI}\) are themselves coboundary expanders. The reduction in \cite{DiksteinD2023swap} then transfers these bounds to the faces complex of \(\CC(\mathrm{Unip}_{d+1};H_1,\dots,H_d)\). In the body of the paper, we apply the same decomposition and Dehn-function arguments to the KMS complexes, together with a vanishing-cohomology argument based on \cite{OppenheimV2025}. This last input accounts for the restriction on the coefficient group in \pref{thm:main}.

% Existing local-to-global theorems reduce swap coboundary expansion of KMS complexes to coboundary expansion of links in their faces complexes \cite{KaufmanKL2016,EvraK2016,OppenheimV2025}. Although these links are not coset complexes, they decompose into matrix block complexes of the form \(Y^{\cI}\).\ydnote{Inga - please suggest a different phrasing} A theorem of \cite{DiksteinD2023swap} then transfers the expansion bounds from these auxiliary complexes to the links, completing the proof.

\subsection{Related works}
Coboundary expansion was defined by \cite{LinialM2006} and independently by \cite{Gromov2010} to study topological properties of simplicial complexes. The definition we use in this paper is due to \cite{DinurM2019} (see also \cite{ChapmanL2023stability}). Swap coboundary expansion was later defined by \cite{BafnaM2023} and \cite{DiksteinD2023agr} independently for purposes of agreement testing. Both works observed that an \(\exp(-O(r))\) lower bound on swap coboundary expansion is insufficient for the intended agreement-testing applications. Subsequent works obtained the first subexponential lower bounds on swap coboundary expansion for bounded-degree simplicial complexes, together with the first applications to agreement testing \cite{DiksteinD2023swap,DiksteinDL2024,BafnaLM2024} and PCP constructions \cite{BafnaMVY2025}. O'Donnell and Singer \cite{ODonnellS2025} obtained a subexponential bound for a variant of the \cite{kaufmanO2018}-coset complexes constructed by \cite{KaufmanOW2025}, following the proof techniques of \cite{BafnaLM2024}. Hopkins and Ray \cite{HopkinsR2026} constructed swap coboundary expanders from projections of spherical buildings, with degree quasipolynomial in the number of vertices and expansion \(\exp(-O(\sqrt{r}\log^2 r))\). Our work is the first to achieve polynomial bounds on bounded-degree complexes. As mentioned above, the work \cite{DiksteinDHS2026} used a different route for constructing PCPs that does not use swap coboundary expansion. We also mention the work of \cite{BenYaacovDL2026} which uses \pref{thm:main} to analyze locally testable codes.

Bounding coboundary expansion of various complexes has proved difficult. Gromov suggested `the cone method' \cite{Gromov2010}, which is the main method we rely on in this work. This method bounds coboundary expansion for symmetric complexes by relating expansion to generalized diameter bounds. This method was later applied to various complexes, improved and generalized (see e.g.\ \cite{LubotzkyMM2016,KaufmanO2021,KozlovM2019,DiksteinD2023cbdry}). A second method is the local-to-global method which originated in \cite{KaufmanKL2016} for dimension \(1\) (which is what is relevant to this paper) and \cite{EvraK2016} (for all dimensions). This method also underwent several quantitative improvements and was generalized to cover all groups as in our current definition of coboundary expansion \cite{KaufmanM2018,KaufmanM2021,KaufmanM2022,DiksteinD2023cbdry}. Our proof combines cone-type arguments with structural properties of KMS complexes.

Finally, the local-to-global method discussed above requires a separate proof for vanishing \(1\)-cohomology. The first constructions of swap coboundary expanders appealed to the congruence subgroup property of the symplectic group \cite{DiksteinDL2024,BafnaLM2024}. More elementary vanishing-cohomology arguments were obtained for the Kaufman--Oppenheim construction by \cite{KaufmanOW2025} and for congruence KMS complexes by \cite{OppenheimV2025}; the latter applies in odd characteristic \(p\) when the coefficient group has no nontrivial elements of order \(p\).

\subsection{Open questions}
It is not clear what is the optimal coboundary expansion constant \(\beta\) one can get on these KMS complexes. Our current cone construction gives an inverse-polynomial bound, with losses coming from the diameter and the number of word reductions used to contract loops. At present, it is conceivable that the true expansion is bounded below by a constant independent of \(r\), or even tends to \(1\) as the field size and dimension grow.

Another interesting question is whether one can get a blackbox reduction from coboundary expansion of \(X\) to swap coboundary expansion. A black-box reduction achieving \(\beta = \exp(-O(r))\) was established \cite{DiksteinD2023swap}. However, all current applications require substantially stronger quantitative bounds.

\subsection{Acknowledgements}
We thank Irit Dinur for insightful discussions on this result.
This work was done in part while the second and third named authors were visiting the Simons Institute for the Theory of Computing. The second named author was partially supported by ISF grant no. 242/24. The third named author was supported by the FWO and the F.R.S.--FNRS under the Excellence of Science (EOS) program (project ID~40007542) for the main part of the project.

\subsection{Organization of this paper}
\pref{sec:prelims} contains the necessary preliminaries. In \pref{sec:overvew}, we state the main technical ingredients and use them to prove the main theorem. These ingredients are proved in \pref{sec: color sampling lemma}, \pref{sec: Unip sec}, \pref{sec:dirprod} and \pref{sec: bounding the Dehn function}. \pref{sec:ai} contains the AI methodology statement. In \pref{app:pcp}, we prove the PCP theorem by adapting the argument of \cite{BafnaMVY2025} to our complex. Finally, in \pref{app:dehn-equivalence}, we prove that the cone method is equivalent to bounding the Dehn function in coset complexes. One direction of this equivalence already essentially appeared in \cite{KaufmanO2021}; we include both directions for completeness and because the equivalence may be of independent interest.

\section{Preliminaries} \label{sec:prelims}
Throughout we use the notation \([m]=\set{0,1,\dots,m}\) where \(m > 0\) is an integer.

\subsection{Simplicial complexes} \label{subsec: simplicial complex}

We start by giving the definition of simplicial complexes and recall some basic properties. 
\begin{definition}
    A \emph{simplicial complex} $X$ consists of a set of vertices $V$ and a collection of subsets of vertices $X \subseteq \mathcal{P}(V)$ such that 
    \begin{enumerate}
        \item $\{v\} \in X$ for all $v \in V$, 
        \item for all $A \in X, B \subset A$ we have $B \in X$ (i.e. $X$ is closed under taking subsets). 
    \end{enumerate}
\end{definition}
 The elements of $X$ are called simplices or faces. Given a simplex $\sigma \in X$ we define its dimension $\dim(\sigma) = \vert \sigma \vert -1$. We will denote $X(k) = \{\sigma \in X \mid \dim(\sigma) = k\}$ for all $k \in \N_0$. We set $X(-1)=\{\emptyset\}$.

As is customary, we will identify the singletons $\{v\}$ with the vertices $v \in V$. In particular, we identify $X(0)$ with $V$ and also write $X(0)$ for the set of vertices. 

We call a simplicial complex finite if its vertex set is finite. We call it $d$-dimensional for some $d \in \N$ if $\max_{\sigma \in X} \dim(\sigma) = d$. Further, $X$ is pure $d$-dimensional if every face is contained in a $d$-dimensional one. 

The \emph{degree} of a simplex $\sigma$ in a pure $d$-dimensional complex $X$ is the number of maximal faces containing $\sigma$, i.e. $\deg(\sigma) = \vert \{\tau \in X(d) \mid \sigma \subseteq \tau \} \vert$. 

The $k$-\emph{skeleton} of $X$, for some $1\leq k \leq d$, is $X^{\leq k}=\bigcup_{i=0}^k X(i)$. This is a $k$-dimensional simplicial complex obtained by forgetting the higher-dimensional structure of $X$. In particular, the 1-skeleton of $X$ is its underlying graph. 

We call a simplicial complex \emph{connected} if its 1-skeleton is connected as a graph. 

\begin{definition} \label{def: link of SC}
    Given a simplicial complex $X$ and a simplex $\sigma \in X$, the \emph{link} of $\sigma$ in $X$ is defined by
    \[X_{\sigma} = \{\tau \in X \mid \tau \cap \sigma = \emptyset, \tau \cup \sigma \in X\}.\]
\end{definition}

\begin{definition} \label{def: partite and color}
    We call a simplicial complex $X$ to be $d+1$-\emph{partite} if there exists a surjective map $\col: X(0) \to [d]=\{0,\dots,d\}$ such that the induced map $\col: X \to \mathcal{P}([d]), \{v_0,\dots,v_k\} \mapsto \{\col(v_0),\dots,\col(v_k)\}$ satisfies $\vert \sigma \vert = \vert \col (\sigma) \vert$ for all $\sigma \in X$, i.e. each simplex consists of vertices with different colors. The map $\col$ is called the \emph{type} or \emph{color} function.
\end{definition}

We use coset complexes to construct nice simplicial complexes from a group and a family of subgroups. 

\begin{definition}\label{def: coset complex}
Let $G$ be a group and $\mathcal{H}=(H_0,\dots , H_d)$ be a family of subgroups of $G$. Then the coset complex $\CC(G;\mathcal{H})$ is defined as the simplicial complex  
\begin{itemize}
\item with vertex set $\bigsqcup_{i=0}^d G/H_i$

\item such that a set of vertices $\left\{ g_1H_{i_1},\dots,g_kH_{i_k} \right\}$ forms a $(k-1)$-simplex in $\CC(G; \mathcal{H})$ if and only if $\bigcap_{j=1}^k g_j H_{i_j} \neq \emptyset$ and all $i_j$'s are distinct.
\end{itemize}
\end{definition}

\begin{remark}\label{rmk: action gives coco}
Given a simplex $\left\{ g_1H_{i_1},\dots,g_kH_{i_k} \right\} \in \CC(G;\mathcal{H})$ we define its color to be the set of indices $\{i_1,\dots, i_k\}$. 
The coset complex $\CC(G;\mathcal{H})$ is a pure $d$-dimensional, $(d+1)$-partite, simplicial complex. The group $G$ acts by $g. \left\{ g_1H_{i_1},\dots,g_kH_{i_k} \right\} = \left\{ gg_1H_{i_1},\dots,gg_kH_{i_k} \right\}$. This action is color preserving and transitive on the maximal faces. The vertex stabilizers of the vertices $1H_i$ are precisely the subgroups $H_i$.

Indeed, every $d$-dimensional, $d+1$-partite simplicial complex with color preserving group action which is transitive on the maximal faces can be described as a coset complex. 
\end{remark}
A key property of coset complexes is that we have a very good description of their links, as the following well-known lemma shows. 
\begin{lemma}[{see e.g. \cite[Lemma 3.4]{HarshaS2019}}] \label{lem: link in coco}
    Let $G$ be a group with subgroups $\mathcal{H}=(H_0,\dots , H_d)$. Set $X=\CC(G;\mathcal{H})$ and let $\sigma \in X$ be a simplex of type $T \subset \{0,\dots,d\}$. Then the link of $\sigma$ in $X$ can again be described as a coset complex: 
    \[X_{\sigma} \cong \CC(H_T;(H_{T \cup \{i\}})_{i \in \{0,\dots,d\}\setminus T})\]
    where $H_T = \cap_{j \in T}H_j$.
\end{lemma}

We equip a finite, pure \(d\)-dimensional complex $X$ with a weight function called Garland weights, which we later use to define a norm on the cochains. Let \(w_d:X(d) \to (0,1]\) be a probability distribution over \(d\) faces. We can extend \(w_d\) to a probability distribution over chains of faces \(\sett{(\sigma_0,\sigma_1,\dots,\sigma_d)}{ \forall i \; \sigma_i \in X(i), \; \forall i<j, \; \sigma_i \subset \sigma_j}\) by the following probabilistic experiment:
\begin{enumerate}
    \item Sample \(\sigma_d \sim w_d\).
    \item Order the vertices of \(\sigma_d\) uniformly at random \((v_0,v_1,\dots,v_d)\).
    \item Output \(\sigma_i = \set{v_0,v_1,\dots,v_i}\) for all \(i=0,1,\dots,d\).
\end{enumerate}
The Garland weights \(w_i: X(i) \to (0,1]\) are the marginal distributions of this joint distribution, explicitly
%For every $\sigma \in X(i)$ we set
\[w_i(\sigma) = \frac{\sum_{\tau \in X(d), \tau \supset \sigma} w_d(\tau)}{\binom{d+1}{i+1}}.\]
When the distribution is not otherwise specified, we use \(w_d(\tau) = \frac{1}{|X(d)|}\) the uniform distribution as top-level distribution. In the following two cases we will choose a different top-level distribution by convention.
Given a complex $X$ with distributions $w_i$ as above, let $1\leq n\leq d$. When we consider the $n$-skeleton of $X$, $X^{\leq n}$ we equip it with the top-level distribution \(w^{\leq n}_n:X(n) \to (0,1]\) which is the \(n\)-th marginal of the top level distribution of \(X\), i.e. $w^{\leq n}_n = w_n$. In particular, we get $w^{\leq n}_i = w_i$ for all $i\leq n$.
Note that the resulting distribution is different then if we would consider $X^{\leq n}$ as an individual complex, ignoring its connection to $X$. In that case, we would choose the top level distribution $w_n^{\leq n}$ to be the uniform distribution. The resulting distribution is however comparable to the one we defined respecting the connection to $X$. A detailed comparison can be found in \cite[Remark 2.3.15]{ingathesis}.

When \(X_{\sigma_i}\) is a link of \(\sigma_i \in X(i)\) then we equip \(X_{\sigma_i}\) with the conditional distribution 
\[w_{d-i-1}(\tau) \propto \begin{cases}
   w_d(\tau \cup \sigma_i) & \tau \cup \sigma_i \in X(d) \\ 0 &\text{otherwise}
\end{cases}\]
where \(w_d\) is the top level distribution over \(X\).

\subsection{High-dimensional expansion}
In this section, we recall the definition of local spectral expansion, coboundary and cosystolic expansion with (potentially) non-abelian coefficients and swap coboundary expansion.
\subsubsection{Local spectral expansion}
 Let $X$ be a finite, pure $d$-dimensional simplicial complex together with the weight function $w:X \to \R_{\geq 0}$ from above. For each $\tau \in X$ of dimension $k$, $w$ induces a weight function on $X_{\tau}$, which we denote by $w_\tau$. Let $X^{\leq 1}_\tau$ denote the 1-skeleton of $X_{\tau}$. We consider the random walk on $X^{\leq 1}_\tau$ with respect to $w_\tau$ and denote the second largest eigenvalue of the random walk matrix by  $\lambda_2(X^{\leq 1}_\tau)$. 
\begin{definition}
    Let $X$ be a finite, pure $d$-dimensional simplicial complex with weight function $w$ and let $\varepsilon \in [0,1)$. We call $X$ an $\varepsilon$-(one sided) \emph{local spectral expander} if for every $\tau \in X$ with $-1\leq  \dim \tau \leq d-2$ we have 
    \[\lambda_2( X^{\leq 1}_\tau) \leq \varepsilon.\]
\end{definition}
Note that in particular the 1-skeleton of $X = X_{\emptyset} $ should be an $\varepsilon$-expander.

\subsubsection{Coboundary and cosystolic expansion with non-abelian coefficients}
For this section, we fix a finite, pure $d$-dimensional simplicial complex $X$ and an arbitrary group $\Lambda$. We denote by $\ordX$ the set of ordered simplices of $X$. I.e. for a simplex $\{v_0,\dots, v_k \} \in X$ we have for every permutation $\pi$ that $v_{\pi(0)} \dots v_{\pi(k)} \in \ordX$. We define the spaces of $-1,0$ and 1-cochains:

\begin{align*}
C^{-1} (X,  \Lambda) &= \lbrace \phi : \lbrace \emptyset \rbrace  \rightarrow \Lambda \rbrace,\\
C^{0} (X,  \Lambda) &= \lbrace \phi : X(0) \rightarrow \Lambda \rbrace,\\
C^1 (X,  \Lambda) &= \lbrace \phi :   \ordX(1) \rightarrow \Lambda \mid \forall uv \in \ordX (1),  \phi (uv) = ( \phi (vu))^{-1} \rbrace.
\end{align*}

Furthermore, we define the coboundary maps $d_k : C^{k} (X,  \Lambda) \rightarrow C^{k+1} (X,  \Lambda)$ as follows: 
\begin{itemize}
\item $k=-1$: For every $\phi \in C^{-1} (X, \Lambda)$ and every $v \in X(0)$,  we define $ d_{-1} \phi (v) = \phi ( \emptyset)$. 
\item $k=0$: For every $\phi \in C^{0} (X,  \Lambda)$ and every $v_0 v_1 \in \ordX (1)$,  we define $d_{0} \phi (v_0 v_1) = \phi (v_1) \phi (v_0)^{-1}$. 
\item $k=1$: We also define the map $d_1 : C^1 (X,  \Lambda) \rightarrow \lbrace \phi : \ordX (2) \rightarrow \Lambda \rbrace$ as follows: for every $\phi \in C^{1} (X,  \Lambda)$ and every $v_0 v_1v_2 \in \ordX (2)$,  we define
$d_{1} \phi (v_0v_1v_2) = \phi (v_0 v_1)  \phi (v_2v_0) \phi (v_1 v_2)$.
\end{itemize}

Note that $(d_0 \circ d_{-1}) = 1$ and $d_1 \circ d_0 = 1$. Thus we can define cocycles and coboundaries for $k=0,1$ by \(Z^k (X, \Lambda) = \lbrace \phi \in C^{k} (X,  \Lambda) : d_k \phi \equiv 1 \rbrace, B^k (X, \Lambda)  = \Img(d_{k-1}) \) and get $B^k (X, \Lambda)  \subseteq Z^k (X, \Lambda)$ (as sets). 

We further define an action of $C^{k-1} (X, \Lambda)$ on $Z^{k} (X, \Lambda)$:
\begin{itemize}
\item An element $\psi \in C^{-1} (X, \Lambda)$ acts on $Z^{0} (X, \Lambda)$ by
\[\psi . \phi (v) = \psi (\emptyset)^{-1} \phi (v),  \forall \phi \in Z^{0} (X, \Lambda),  \forall v \in X (0).\]
\item An element $\psi \in C^{0} (X, \Lambda)$ acts on $Z^{1} (X, \Lambda)$ by
\[\psi . \phi (v_0 v_1) = \psi (v_1) \phi (v_0 v_1) \psi (v_0)^{-1},  \forall \phi \in Z^{1} (X, \Lambda),  \forall v_0 v_1 \in \ordX (1).\]
\end{itemize}

We denote by $[ \phi ]$ the orbit of $\phi \in Z^{k} (X, \Lambda)$ under the action of $C^{k-1} (X, \Lambda)$ and define the reduced $k$-th cohomology to be
\[H^k (X, \Lambda) = \lbrace [ \phi ] : \phi \in Z^{k} (X, \Lambda) \rbrace.\]
We will say that the $k$-th cohomology is trivial if $H^k (X, \Lambda) = \lbrace [ \phi \equiv 1]  \rbrace$, denoted $H^k(X,\Lambda) = 1$. Note that the cohomology is trivial if and only if $B^k (X, \Lambda) = Z^k (X, \Lambda)$.

We define the support of cochains as follows. For $\phi \in C^0 (X, \Lambda)$, we define
\[\supp (\phi) = \lbrace v \in X(0) : \phi (v) \neq 1 \rbrace.\]
Also, for $\phi \in C^{1} (X, \Lambda)$, we define
\[\supp (\phi) = \lbrace \lbrace u,v \rbrace \in X(1) : \phi (uv) \neq 1 \rbrace,\]
and
\[\supp (d_1 \phi) =  \lbrace \lbrace v_1,v_2,v_3 \rbrace \in X(2): (d_1\phi) (v_1v_2v_3) \neq 1 \rbrace.\]

With this notation, we define the following norms: for $\phi \in  C^{k} (X, \Lambda)$, $k=0,1$,  we set
\[\Vert \phi \Vert = \sum_{\tau \in \supp (\phi)} w (\tau).\]
Also,  for $\phi \in  C^{1} (X, \Lambda)$, we define
\[\Vert d_1 \phi \Vert = \sum_{\tau \in \supp (d_1 \phi)} w (\tau).\]

Lastly, for two maps $\phi, \psi \in  C^{k} (X, \Lambda),  k=0,1$,  we define $\dist (\phi, \psi)$ as follows: For $k=0$,  we set
\[\dist (\phi, \psi) = w \left( \lbrace v \in X(0) \mid \phi (v) (\psi (v))^{-1} \neq 1 \rbrace \right).\]
For $k=1$,  we define 
\[\dist (\phi, \psi) = w \left( \lbrace \lbrace u,v \rbrace \in X(1) \mid \phi (uv) (\psi (uv))^{-1} \neq 1 \rbrace \right).\]

\begin{definition}
Using the notation above,  we will define for $k=0,1$ the $k$-\emph{cocycle expansion constant} of $X$ with respect to $\Lambda$ as 
\[h^{k}(X, \Lambda) = \min_{\phi \in C^{k} (X, \Lambda) \setminus Z^{k} (X, \Lambda)} \dfrac{\Vert d_k \phi \Vert}{\min_{\psi \in Z^k (X,\Lambda)} \dist (\phi,  \psi)}.\]
\end{definition}

\begin{definition}
    For a constant $\beta >0$,  we say that
\begin{itemize}
\item The complex $X$ is a \textit{$(\Lambda,  \beta)$ $1$-coboundary expander},  if $h^{0} (X, \Lambda) \geq \beta$, $h^{1} (X, \Lambda) \geq \beta$ and $H^1(X,\Lambda)$ is trivial.
\item The complex $X$ is a \textit{$(\Lambda,  \beta)$ $1$-cocycle expander},  if $h^{0} (X, \Lambda) \geq \beta$,  $h^{1} (X, \Lambda) \geq \beta$. If furthermore for every $\phi \in Z^1 (X,  \Lambda) \setminus B^1 (X,  \Lambda)$ it holds that $\Vert \phi \Vert \geq \beta$ then we call $X$ a \textit{$(\Lambda,  \beta)$ $1$-cosystolic expander}.
\end{itemize}

\end{definition}

\ivnote{I called the statement below an observation and added a proof because the AI was complaining about the statement. Feel free to remove the proof if you think it's unnecessary. }
\begin{observation}
A lower bound on $h^{0}$ is usually easy to verify: If the $1$-skeleton of $X$ is a connected finite graph, then for any non-trivial group $\Lambda$, $h^{0} (X, \Lambda)$ is exactly the Cheeger constant of the weighted $1$-skeleton (where each edge $\lbrace u,v \rbrace$ is weighted by $w(\lbrace u,v \rbrace)$).
\end{observation}
\begin{proof}
    Note that the weighted Cheeger constant is 
    \[h(X^{\leq 1}) = \min_{\emptyset \neq A \subsetneq X(0)} \frac{w(E(A,A^c))}{\min(w(A),w(A^c))} \]
    where $A^c=X(0) \setminus A$ denotes the complement of $A$ and $E(A,A^c)=\{\{v,w\} \in X(1) \mid v \in A, w \notin A\}$. 

    To see that $h^0(X,\Lambda) \leq h(X^{\leq 1})$, fix $\emptyset \neq A \subsetneq X(0) $ minimizing for $h(X^{\leq 1})$ and define a cochain $f:X(0) \to \Lambda$ by fixing an arbitrary $1 \neq g\in \Lambda$ and setting $f(v)=\begin{cases}
        g & \text{if } v \in A \\ 1 &\text{else}
    \end{cases}$.
    Then $\supp(d_0f) = E(A,A^c)$ and hence $\Vert d_0f \Vert = w(E(A,A^c))$. Furthermore, note that since $X^{\leq 1}$ is connected $Z^0(X,\Lambda)=B^0(X,\Lambda)$ are all constant functions. Thus 
    \[\min_{b\in Z^0(X,\Lambda)} \dist(f,b) = \min(\dist(f,[v \mapsto 1]), \dist(f,[v \mapsto g]), \dist(f,[v \mapsto h] : h \neq g,1) = \min (w(A),w(A^c)).\]
    Together we get
    \[h(X^{\leq 1}) = \frac{w(E(A,A^c))}{\min(w(A),w(A^c))} = \frac{\Vert d_0f \Vert}{\min_{b \in Z^0(X,\Lambda)} \dist(f,b) } \geq h^0(X,\Lambda).\]

    To see the other inequality, $h^0(X,\Lambda) \leq h(X^{\leq 1})$, fix an arbitrary cochain $f \in C^0(X,\Lambda) \setminus Z^0(X,\Lambda).$ Then $f$ partitions the set of vertices
    \[X(0) = \sqcup_{g \in \Lambda} V_g \text{ with } V_g = f^{-1}(\{g\}).\]
    Firstly, observe that for a constant map $b=[v \mapsto g]$ we have $\dist(f,b) = w(\{v\in X(0) \mid f(v) \neq g\}) = w(V_g^c)$. Next, we use the definition of $d_0f$ and that $h(X^{\leq 1}) \min(w(A),w(A^c)) \leq w(E(A,A^c))$ for all $\emptyset \neq A \subset X(0)$ to obtain
    \[\Vert d_0f \Vert = \frac12 \sum_{g \in \Lambda} w(E(V_g, V_g^c))\geq \frac12 \sum_g h(X^{\leq 1})\min(w(V_g), w(V_g^c)).\]
To compare this to the distance of $f$ to the cocycles/coboundaries, we consider two cases. First, assume that $\min(w(V_g), w(V_g^c)) = w(V_g)$ for all $g$. Thus, in particular $w(V_g) \leq \frac12$ and hence $\min_{b \in Z^0(X,\Lambda)} \dist (f,b) = \min_g w(V_g^c) \geq \frac12$. And we have $\sum_{h\in \Lambda} \min(w(V_h),w(V_h^c))= \sum_{h\in \Lambda}w(V_h) = w(X(0)) =1.$

On the other hand, assume that there exists a $g \in \Lambda$ such that $w(V_g)>\frac12$ (note the there can exist at most one such $g$). Then 
\[\sum_{h\in \Lambda} \min(w(V_h),w(V_h^c)) = w(V_g^c) + \sum_{h\neq g} w(V_h = 2(w(V_g^c)) \geq 2\min_b \dist(f,b).\] 
Together, we get that
\begin{align*}
   \Vert d_0f \Vert \geq \frac12 \sum_g h(X^{\leq 1})\min(w(V_g), w(V_g^c)) \geq \frac12 h(X^{\leq 1}) 2\min_b \dist(f,b)
\end{align*}
and hence $h^0(X,\Lambda) \geq h(X^{\leq 1}).$
    
\end{proof}

\begin{comment}
In \cite{DiksteinD2023cbdry}, the following result was proven, improving the parameters over the result of \cite{KaufmanM2021}.
\begin{theorem}[{\cite[Theorem 8]{DiksteinD2023cbdry}}]
\label{thm: local to global nonabelian}
Let $0 \leq \lambda <1,  \beta >0$ be constants and $\Lambda$ be a group.  For a finite pure $n$-dimensional simplicial complex $X$ with $n \geq 3$,  if
\begin{itemize}
    \item for every vertex $v$,  $h^1 (X_{v}, \Lambda) \geq \beta$ and $H^1(X_{v},\Gamma)$ is trivial,
    \item and $X$ is a $\lambda$-local spectral expander,  
\end{itemize}
    then
\[h^1 (X, \Lambda) \geq \frac{(1-\lambda) \beta}{24} - e \lambda,\]
where $e \approx 2.71$ is the Euler constant.
\end{theorem}
\end{comment}
\subsubsection{Swap coboundary expansion}

Swap coboundary expansion of a complex $X$ in essence means that the $r$-faces complex of $X$ is a coboundary expander. This notion was defined independently by \cite{DiksteinD2023agr} and \cite{BafnaM2023} in order to analyze property tests used in PCPs.

We start by defining the faces complex and related notions to then state the precise definition of swap coboundary expansion below in \pref{def: swap CBE}.

\begin{definition}[{\cite[Definition 2.11]{DiksteinD2023swap}}]
Let $X$ be a $d$-dimensional simplicial complex and $1 \leq r \leq d$. We denote by $F^rX$ the simplicial complex with vertex set $X(r)$ and where a set of vertices $\{\sigma_0,\dots, \sigma_k\}$ forms a simplex if and only if the $\sigma_i$ are pairwise disjoint and $\bigcup_{i=0}^k \sigma_i \in X$. We call $F^rX$ the $r$-\emph{faces complex} of $X$.    
\end{definition}
The dimension of $F^rX$ is \(t= \lfloor \frac{d+1}{r+1} \rfloor -1\). Let \(t'=(t+1)(r+1)-1\).  We equip $F^rX$ with the following measure. On top level faces of $F^rX$, the measure is given by sampling a \(t'\)-dimensional simplex in \(X\) (according to the measure of \(X\) and then partitioning it into \(t+1\) faces of dimension \(r\) uniformly at random. For lower dimensional faces, we consider the induced measure as usual. Explicitly, if $w_i$ denote the weight functions on $X$ we get
\[
w_{F^rX,\ell}(\{\sigma_0,\dots,\sigma_{\ell}\}) = \frac{(\ell +1)! ((r+1)!)^{\ell+1}}{((l+1)(r+1))!}w_{(r+1)(\ell+1)-1}(\cup_i \sigma_i). 
\]
\begin{definition}\label{def: swap CBE}
    Let $X$ be finite, pure $d$-dimensional simplicial complex, $\Lambda$ a group and $\beta >0, 1 \leq r\leq d$ be constants. We call $X$ a \emph{$(\beta,r,\Lambda)$-swap coboundary expander} if $F^rX$ is a $(\beta,\Lambda)$ 1-coboundary expander. 
\end{definition}

\begin{observation} \label{obs: LSE of full faces complex}
For a $d$-dimensional $(d+1)$-partite complex $Y$, mapping each face to its set of colors induces a surjective graph homomorphism from the one-skeleton of $F^rY$ to the Kneser graph $K(d+1,r+1)$, and similarly for links with the corresponding parameters. Here $K(d+1,r+1)$ is the graph whose vertices are the $(r+1)$-element subsets of $\set{0,1,\dots,d}$, with two vertices adjacent if their sets are disjoint.

Thus, whenever the corresponding Kneser graph is disconnected, so is the link. In particular, if $r+1$ divides $d+1$, the codimension-$2$ links map surjectively onto $K(2(r+1),r+1)$, which is a disjoint union of edges, and hence these links are disconnected. Even when $d+1$ is not a multiple of $r+1$, a similar spectral argument shows that the codimension-$2$ links cannot be $\lambda$-one-sided spectral expanders for any $\lambda<\frac14$.
\end{observation}

The faces complex of a coset complex is not again a coset complex. Thus, it is much easier to analyze the colored faces complex of a coset complex, which we will define now, and which will turn out to be again a coset complex.

\begin{definition}
    Let $X$ be a finite, pure $d$-dimensional, $d+1$-partite simplicial complex with color function $\col: X \to \mathcal{P}(\{0,\dots,d\})$. Let $1 \leq r,T \leq d$. 
    Fix a set of colors $\cI = \{I_1,\dots, I_T\}$ such that $I_j \subseteq \{0,\dots,d\}$, \(\vert I_j \vert = r+1\), and such that the $I_j$ are pairwise disjoint. 
    We denote by $X[I_j] = \{\sigma \in X \mid \col(\sigma) = I_j\}$ the faces with color $I_j$.
    
    Then we define the \emph{colored faces complex} $F^{\cI}X$ of $X$ with respect to $\cI$ as simplicial complex  
    \begin{itemize}
        \item with vertex set $\bigcup_{j=0}^T X[I_j]$,
        \item such that a set of vertices $\{\sigma_0,\dots, \sigma_\ell\}$ forms a simplex if there exist distinct $j_0,\dots,j_\ell$ such that $\sigma_i \in X[I_{j_i}]$ and $\bigcup_i \sigma_i \in X$.
    \end{itemize}
\end{definition}
Note that $F^{\cI}X$ is $\vert \T \vert$-partite and $(T-1)$-dimensional. We equip $F^{\cI}X$ with the measure obtained by the measure on \(F^r X\) conditioned on sampling a face whose vertices have colors \(\cI\). Denote the weight on $X$ by $w_i$, we get for $\sigma=\{\sigma_0,\dots,\sigma_{\ell}\} \in F^{\cI}X$ the weight
\[
w_{F^{\cI}X,\ell}(\sigma) = \frac{w_{(\ell+1)(r+1)-1}(\cup_i \sigma_i)}{\sum_{\{\tau_0,\dots,\tau_{\ell}\} \in F^{\cI}X(\ell)}w_{(\ell+1)(r+1)-1}(\cup_i \tau_i)}.
\]

\begin{fact}\label{lem: colored faces coset complex}
    Let $X = \CC(G;(H_i)_{i \in [d]})$ be a coset complex and $\cI = \{I_1,\dots, I_T\}$ a set of disjoint colors as above. Set $H_I = \cap_{i \in I}H_i$ for each $I \in \cI$. Then 
    \begin{align*}
        F^{\cI}X \cong \CC(G;(H_I)_{I\in \cI}).
    \end{align*}
\end{fact}

\begin{definition}
    Let $1\leq r \leq d \in \N$ and $T\in \N$. We denote by $\binom{[d]}{r+1}$ the set of subsets of size $r+1$ of $[d] =\{0,\dots,d\}$. We also introduce the following notation for the set of all sets of $T$ disjoint colors of size $r+1$ of $[d]$:
    \[\binom{[d]}{r+1}^{[T]} = \left\{\cI \subseteq \binom{[d]}{r+1} \mid \vert \cI \vert = T \text{ and } \forall I,J \in \cI, I\neq J: I \cap J = \emptyset \right\}.\]
\end{definition}

\begin{comment}
We need the following result about the expansion of the $r$-face complex given expansion of the colored face complex. \ivnote{Check termnilogy}

\begin{theorem}[{\cite[Lemma 8.2]{DiksteinD2023swap}}] \label{thm: swap CSE from CSE for many colored faces}
Let $X$ be a finite, pure $d$-dimensional, $(d+1)$-partite simplicial complex. Let $1 \geq r \geq ??? $ and $\cC = \binom{[d]}{r+1}$ be the set of subsets of $[d]=\{0,\dots,d\}$ of size $r+1$, i.e. all possible sets of colors of $r$-dimensional faces of $X$. Let $m = \sqrt{r+1}$. We sample a set $\cI$ of $m$ pairwise disjoint colors from $\cC$ uniformly. If $\displaystyle \Psymb_{\cI}[h^1(F^IX,\Lambda) \geq \beta] = p$ then $h^1(F^rX,\Lambda) = \Omega(\beta p^2)$.    
\end{theorem}
\end{comment}

\subsection{KMS complexes} \label{subsec:KMS-complexes}
We now describe the main objects of our study, the congruence KMS complexes. 
We restrict to the congruence KMS complexes of type $A_d$, hence we work with the Chevalley group of type $A_d$ which is $\Sl_{d+1}$.
We fix a prime power $q$ and a size $d \geq 2$. For a polynomial $f\in \F_q[t]$ and $1 \leq i,j \leq d+1, i\neq j$ we denote by $e_{i,j}(f)$ the matrix in $\Sl_{d+1}(\F_q[t])$ with ones on the diagonal and $f$ at position $(i,j)$. We then define for $i=1,\dots, d$
\[U_i = \langle e_{i,i+1}(\lambda) \mid \lambda \in \F_q \rangle \leq \Sl_{d+1}(\F_q[t])\]
and 
\[U_0 = \langle e_{d+1,1}(\lambda t) \mid \lambda \in \F_q \rangle\leq \Sl_{d+1}(\F_q[t]).\]
For each subset $J \subseteq  [d], i\in [d]$ we define
\[U_J = \langle U_i \mid i\in J \rangle \leq \Sl_{d+1}(\F_q[t]), \qquad H_i = U_{[d]\setminus \{i\}}\]
and we call $\KMS{d}{q} = U_{[d]}$ the \emph{concrete KMS group in $\Sl_{d+1}(\F_q[t])$}.

\begin{remark}
    For $d\geq 3$ the concrete KMS group $\KMS{d}{q}$ is isomorphic to the KMS group of type $\widetilde{A}_d$ over $\F_q$ as defined in \cite{GraveDePeraltaV2025,oppenheimVB2025new,OppenheimV2025}.    
\end{remark}
\begin{example}
For the case $d=2$ we have 
    \[U_{\{1,2\}}= \begin{pmatrix} 1 & \F_q & \F_q \\ 0 & 1 & \F_q \\ 0 & 0 & 1\end{pmatrix}, \qquad U_{\{0,1\}}= \begin{pmatrix} 1 & \F_q & 0 \\ 0 & 1 & 0 \\ \F_qt & \F_qt & 1\end{pmatrix}, \qquad U_{\{0,2\}}= \begin{pmatrix} 1 & 0 & 0 \\ \F_qt & 1 & \F_q \\ \F_qt & 0 & 1\end{pmatrix}.\]

    Here we use a shortened notation, by which we for example mean 
\[\begin{pmatrix} 1 & 0 & 0 \\ \F_qt & 1 & \F_q \\ \F_qt & 0 & 1\end{pmatrix} = \left\{ \begin{pmatrix} 1 & 0 & 0 \\ \lambda_1 t & 1 & \lambda_2 \\ \lambda_3 t & 0 & 1\end{pmatrix}\mid \lambda_i \in \F_q\right\}.\]
Then $\KMS{2}{q}$ is the group generated by $U_{\{1,2\}},U_{\{0,1\}},U_{\{0,2\}}$ in $\Sl_3(\F_q[t])$. It consists of all matrices that are unipotent upper triangular when considered modulo $t$.
\end{example}

Fix an irreducible polynomial $f\in \F_q[t]$ of degree at least two. The quotient map $\F_q[t] \to \F_q[t]/(f)$ induces a map 
\[\pi_f: \Sl_{d+1}(\F_q[t]) \to \Sl_{d+1}(\F_q[t]/(f)).\]
We set $H^f_i = \pi(H_i) = \pi_f(U_{[d]\setminus \{i\}})$ for all $i\in [d]$. The \emph{congruence KMS complex of type $A_d$ over $\F_q[t]/(f)$} is 
\[X(A_d, \F_q,f) = \CC(\Sl_{d+1}(\F_q[t]/(f));(H^f_i)_{i \in [d]}).\]

For $J \subset [d]$ we set $H_J = \cap_{i \in J} H_i$. Note that $H_J = U_{[d]\setminus J}$ and $\pi_f(H_J) = \cap_{i\in J} H_i^f$.

Let $v \in X(A_d,\F_q,f)(0)$ be a vertex. Then the link of $v$ in $X= X(A_d,\F_q,f)$ can be described as follows. 
Let $U^+$ be the subgroups of unipotent upper triangular matrices in $\Sl_{d+1}(\F_q)$ and let $K_i = \langle e_{j,j+1}(\lambda)\mid \lambda \in \F_q, j \in \{1,\dots,d\} \setminus \{i \} \rangle$.
Then 
\[X_{v} \cong \CC(U^+;(K_i)_{i \in \{1,\dots,d\}}).\]

This complex is isomorphic to a subcomplex of the spherical building of type $A_d$ called the opposite complex. However, we will not make use of this fact in this work.

\subsection{Previous results on KMS complexes}
In \cite{GraveDePeraltaV2025} the authors show that KMS complexes are local spectral expanders. More precisely and adapted to our case they show the following. 
\begin{theorem}
    Let $q \geq 4, d \geq 2$ and $f\in \F_q[t]$ be an irreducible polynomial of degree at least 2. Further, assume that ${\sqrt{q}} \geq d$. Then $X(A_d,\F_q,f)$ is a $\frac{1}{\sqrt{q}-(d-1)}$-local spectral expander.    
\end{theorem}

\ydnote{I changed to a better estimate on on spectral expansion:}
We use a slightly better local-to-global theorem by \cite{LeakeO2025} as a black box to obtain local spectral expansion that is independent of dimension.

\begin{theorem} \label{thm: KMS complexes are local spectral}
Let $d\geq 2$, let $q>16$ be a prime power, and let $f\in\F_q[t]$ be irreducible with $\deg f\geq 2$. Set $X=X(A_d,\F_q,f)$. Then for every face $\tau\in X$, including the empty face, with $|\tau|\leq d-1$, we have
    \[\lambda_2(X_\tau^{\leq 1})\leq \frac{2}{(d-|\tau|)(\sqrt{q}-2)}.\]
In particular \(X\) is a \(\frac{2}{\sqrt{q}-2}\)-one-sided local spectral expander.
Moreover, for any $2\leq n\leq d$, then $X^{\leq n}$, with the weights inherited from $X$, is a $\frac{1}{d-n+1}$-two-sided local spectral expander.
\end{theorem}

\begin{proof}
The links of $X$ are connected. A codimension-two link whose two remaining colors are adjacent in the cyclic diagram of type $\widetilde A_d$ has second eigenvalue at most $q^{-1/2}$; for nonadjacent colors, the link is complete bipartite and its second eigenvalue is zero \cite{GraveDePeraltaV2025}. Thus the pairwise spectral influence matrix is bounded entrywise by $q^{-1/2}$ times the adjacency matrix of a cycle, and its largest eigenvalue is at most $2/\sqrt q$. Applying \cite[Theorem 1.3]{LeakeO2025} gives the stated bound.

For the two-sided conclusion, recall that the smallest eigenvalue of the weighted graph of a pure $D$-dimensional complex is at least $-1/D$, since its edge distribution is a mixture of the edge distributions of $D$-simplices.
Here $\dim X_\tau=d-|\tau|$. When $q\geq 16$, we have $\frac{2}{\sqrt{q}-1}\leq 1$, so every nontrivial eigenvalue of $X_\tau^{\leq 1}$ has absolute value at most $1/(d-|\tau|)$. Every face whose link must be checked in $X^{\leq n}$ has size at most $n-1$, and its weighted link graph is unchanged by taking the skeleton. This proves the claim.
%Here $\dim X_\tau=d-|\tau|$. When $q\geq 16$, we have $\eta_q\leq 1$, so every nontrivial eigenvalue of $X_\tau^{\leq 1}$ has absolute value at most $1/(d-|\tau|)$. Every face whose link must be checked in $X^{\leq n}$ has size at most $n-1$, and its weighted link graph is unchanged by taking the skeleton. This proves the claim.
\end{proof}

\cite{OppenheimV2025} also proved these complexes are 1-coboundary expanders. We will only use part of the proof method for the vanishing of cohomology later, but state the result for completenes.
\begin{theorem} \label{thm: CBE for KMS}
    Let $d \geq 3$. There exist $q_0,\varepsilon>0$ such that 
    \begin{itemize}
        \item for every prime power $q=p^m \geq q_0$,
        \item for every group $\Lambda$ without non-trivial elements of order $p$ (e.g. if $\Lambda$ is finite and $\vert \Lambda \vert <p$),
        \item for every irreducible polynomial $f\in \F_q[t]$ of degree at least 2
    \end{itemize}
    we have
    \[h^1(X(A_d,\F_q,f),\Lambda) \geq \varepsilon \text{ and } H^1(X(A_d,\F_q,f),\Lambda) = 1.\]
    In particular, $X(A_d,\F_q,f)$ is a $(\Lambda,\varepsilon)$ 1-coboundary expander.
\end{theorem}

\begin{comment}
\subsection{Density of opposite complex inside building}
Let $\Delta$ be the building of type $A_n$ over $\F_q$ with opposite complex $\Delta^{op}$. Set
$$F_n(q) = \frac{g^{n+1}-1}{q-1} = q^n+q^{n-1} + \dots + q+1.$$
Then the number of chambers, i.e of $n$-dimensional simplices in $\Delta$ and $\Delta^{op}$ is 
$$\vert \Delta(n) \vert = F_1(q)F_2(q) \cdots F_n(q), \ \vert \Delta^{op}(n) \vert = q^{\frac{n(n+1)}{2}}.$$

We can deduce that $\frac{\vert \Delta^{op}(n) \vert}{\vert \Delta(n)\vert}$ behaves like $1-\frac{n}{q}$ as $q \to \infty$. \ivnote{This holds for any type of building I think.}\ydnote{Do you know whether this is proven somewhere?}

For the proof of good cone function, need: Diameter of $P^I$ bounded, $\Psymb_{g\in G}[g.w \notin P] \leq \frac{1}{q-1}$ for a vertex $w$ and $G$ all $\Sl_n(\F_q)$.
$\Psymb_{g\in G}[g.w \notin P] = \frac{\vert \{g.w \mid g \in G\} \setminus P \vert}{\vert \{g.w \mid g\in G\}\vert}$. Let $\dim w = i$, then $\{g.w \mid g\in G\} = S[i]$ by transitivity of the action. So this is $1-\frac{\vert P[i] \vert}{\vert S[i] \vert}$. 

Maybe we can also compute this for arbitrary types using the sizes of the groups i.e. $\vert S[i] \vert = \vert G \vert / \vert P_i \vert$ and similarly with $U^+$.  This should also be known. 

\ivnote{Don't use $P$ for the opposite vertices, there are way too many other things named $P$.}
\end{comment}

\subsection{Dehn functions}

The notion of Dehn functions will be useful to use below as a tool to bound coboundary expansion in coset complexes. In the section, we cover the basic definitions and some results regarding Dehn function that we will need in the sequel.

Let $G = \langle S \vert R \rangle$ be a finitely presented group, where $S$ is a symmetric generating set of $G$ and $R$ is a set of relations. Denote $\Free (S)$ to be the free group with a generating set $S$. 

\begin{definition}[Word length]
Any element in $\Free (S)$ can be written as a word in the alphabet $S$ and for a word $w$ in the alphabet $S$, we denote $\vert w \vert_S$ to be the length of the word (i.e., the number of letters).
\end{definition}

Without loss of generality, we can assume that $R \subseteq \Free (S)$ are cyclically reduced words.  A \textit{relation} in $G$ is a word $w \in \Free (S)$ such that $w=1$ in $G$. We note that every relation is in the normal closure of $R$ in $\Free (S)$. In the following, we will deal with different generating sets for the same group $G$, in that situation (when $G$ is obvious from the context), we will say that $w$ is a \textit{relation in the alphabet $S$} to specify to what generating set we are referring to.

Given a relation $w \in \Free (S)$, the \textit{area with respect to $R$} of $w$, denoted $\Area_{R} (w)$, is the minimal number $k \in \mathbb{N} \cup \lbrace 0 \rbrace$ such there are $\rho_1,...,\rho_k \in R \cup R^{-1}$ and $u_1,...,u_k \in \Free (S)$ such that
$$w = u_1^{-1} \rho_1 u_1 ... u_k^{-1} \rho_k u_k \text{ in } \Free (S).$$

\begin{definition}[The Dehn function]
\label{Dehn func def}
The \textit{Dehn function of $G$ with respect to $S,R$} is the function $\Dehn_{S,R} : \mathbb{N} \rightarrow \mathbb{N}$ defined as
$$\Dehn_{S,R} (m) = \max \lbrace \Area_R (w) :w = 1 \text{ in } G, \vert w \vert_S \leq m \rbrace.$$
\end{definition}

The idea behind the Dehn function is that it counts how many reduction moves are needed to reduce $w$ to the trivial word. More precisely, following \cite[Chapter 8]{clay2017office}, we define the following moves on a word $w$:
\begin{enumerate}
\item Free reduction: remove a sub-word of the form $s s^{-1}$ within the word $w$, where $s \in S$.
\item Free expansion: Replace $w=uv$ (where $u,v$ are words) with $w' = u s s^{-1} v$ or with $w' = u e v$.
\item Applying a relation: replace a subword $w'$ in $w$ with a new subword $w''$ where $w' (w'')^{-1}$ or $w'' (w')^{-1}$ is a cyclic permutation of a word in $R \cup R^{-1}$.
\end{enumerate}

Every $w=1$ in $G$ can be reduced to $1$ by applying the above moves, and a different way to define $\Dehn_{S,R} (m)$ is as follows: $\Dehn_{S,R} (m)$ is the minimal number $N$ such that for every $w$ with $\vert w \vert \leq m$,  and $w=1$, the word $w$ can be reduced to $1$ by applying most $N$ relations. It does not take into account the number of free reductions or free expansions that are applied, and we will refer to these relations as \textit{trivial relations}. 

Given a presentation as above and two words $w, w'$, we will say that \textit{the word $w'$ can be deduced from $w$}, if one can apply the above moves on $w$ and get $w'$. 

In the sequel, it will be useful to compare Dehn function of the same group with two different presentations. Let $G$ be a group with two different presentations in terms of generators and relations: $G = \langle S \vert R \rangle$ and $G = \langle S' \vert R' \rangle$. It is well-known that $\Dehn_{S,R}$ and $\Dehn_{S',R'}$ are the same up to a linear distortion, but we will need to keep track of the constants of this distortion. We define homomorphism $\phi : \Free (S) \rightarrow \Free (S')$ and $\phi ' : \Free (S') \rightarrow \Free (S)$ as follows: Every element $s \in S$ is in particular an element in $G$ and as such can be written as a word in $S'$. For every element $s \in S$ we choose $w_s \in \Free (S')$ such that $s=w_s$ in $G$ and $\vert w_s \vert_{S'}$ is minimal. We define $\phi (s) = w_s$ and extend $\phi$ to a homomorphism $\phi : \Free (S) \rightarrow \Free (S')$. Similarly, for every $s' \in S'$, we choose $w_{s'} \in \Free (S)$ such that $s' = w_{s'}$ in $G$ and $\vert w_{s'} \vert_S$ is minimal and define a homomorphism $\phi ' : \Free (S') \rightarrow \Free (S)$ as above.

\begin{proposition}
\label{prop: comparing dehn for dif pres prop}
 Let $G$ be a group with two different presentations in terms of generators and relations: $G = \langle S \vert R \rangle$ and $G = \langle S' \vert R' \rangle$. For choices of $\phi, \phi '$ as above, we define the following constants:
\begin{itemize}
    \item $C_1 = \max_{s \in S} \vert \phi (s) \vert_{S'}$.
    \item We note that for every $\rho ' \in R'$, $\phi ' (\rho ')$ is a relation in the alphabet $S$ and define $C_2 = \max_{\rho ' \in R'} \Area_R (\phi ' (\rho '))$.
    \item We note that for every $s \in S$,  $s (\phi ' (\phi (s)))^{-1}$ is a relation in the alphabet $S$ and define $C_3 = \max_{s \in S} \Area_{R} (s (\phi ' (\phi (s)))^{-1})$.
\end{itemize}
Then it holds for every $m \in \mathbb{N}$ that 
$$\Dehn_{S,R} (m) \leq C_2 \Dehn_{S',R'} (C_1 m) + C_3 m.$$
\end{proposition}
\ivnote{I feel like we should either motivate why this is well known/classical or put it in Chapter 5 where we need it.}
\begin{proof}
Fix $m \in \mathbb{N}$ and let $w$ be a relation in the alphabet $S$ with $\vert w \vert_S \leq m$. We denote $w = s_1 ... s_j$ where $s_1,...,s_j \in S$ (not necessarily distinct) where $j \leq m$. We note that by applying a free expansion $w = s_1 (\phi ' (\phi (s_1)))^{-1} (\phi ' (\phi (s_1))) s_2 ... s_j$. Thus, by the definition of $C_3$, we can deduce $(\phi ' (\phi (s_1))) s_2 ... s_j$ from $w$ by applying at most $C_3$ relations in $R$. Continuing in the same manner, we can deduce $(\phi ' (\phi (s_1))) ... (\phi ' (\phi (s_j))) = \phi ' (\phi  (s_1 ... s_j))$ by applying at most $C_3 j \leq C_3 m$ relations in $R$ (and $m$ free expansions that we do not count). Thus, 
$$\Area_{R} (s_1 ... s_j) \leq C_3 m + \Area_{R} (\phi ' (\phi  (s_1 ... s_j))).$$

We will finish the proof by showing that $\Area_{R} (\phi ' (\phi  (s_1 ... s_j))) \leq C_2 \Dehn_{S',R'} (C_1 m)$. We note that $w = s_1 ... s_j$ is a relation in the alphabet $S$ and thus $\phi (s_1 ... s_j)$ is a relation in the alphabet $S '$ and 
$$\vert \phi (s_1 ... s_j) \vert_{S'} \leq \sum_{i=1}^j \vert \phi (s_i ) \vert_{S'} \leq C_1 j \leq C_1 m.$$
By the definition of $\Dehn_{S', R'}$ there is $k \leq \Dehn (C_1 m)$, words $u_1,...,u_k \in \Free (S ')$ and relations $\rho_1 ',...,\rho_k ' \in R'$ such that 
$$\phi (s_1 ... s_j) = u_1^{-1} \rho_1 ' u_1 u_2^{-1} \rho_2 ' u_2 ... u_k^{-1} \rho_k ' u_k.$$
Therefore 
\begin{align*}
    & \Area_{R} (\phi ' (\phi  (s_1 ... s_j))) \leq \\
    & \Area_{R} (\phi ' (u_1^{-1}) \phi ' (\rho_1 ') \phi ' (u_1) ...  \phi ' (u_k^{-1}) \phi ' (\rho_k ') \phi ' (u_k) \leq \\
    & \sum_{i=1^k} \Area_R (\phi ' (\rho_i ')) \leq C_2 k \leq C_2  \Dehn (C_1 m),
\end{align*}
as needed.
\end{proof}

\subsection{Dehn functions, cone functions and coboundary expansion} \label{subsec:dehn}
Non-abelian cone functions were introduced in \cite{DiksteinD2023swap}, our definition is very similar. We follow the exposition in \cite{OppenheimV2025}.

Throughout, $X$ is an $n$-dimensional simplicial complex. We again denote by $\ordX(k)$ the set of ordered $k$-dimensional simplices of $X$, and for a simplex $\{v_0,\dots,v_k\} \in X$ we denote by $v_0\dots v_k \in \ordX$ the ordered simplex.

Given two vertices $u,v$ of $X$, a \emph{path from $u$ to $v$} is a $(k+1)$-tuple of vertices $(u_0 ;... ; u_k)$ where $u = u_0$ and $v =u_k$ such that for every $0 \leq i \leq k-1$,  $u_iu_{i+1} \in \ordX (1)$. The constant $k$ is called the \textit{length of the path}. Given two paths $P_0  =(u_0; ... ;u_k),  P_1 = (u_k ; ... ;u_m)$,  we define $P_0 \circ P_1 = (u_0; ... ;u_m)$. Also, given a path $P = (u_0; ... ;u_k)$, we denote $P^{-1} = (u_k ;...; u_0)$.

A \textit{loop around $v$} is a path from $v$ to itself.  We note that $(v)$ is also considered as a loop of length 0 and we will call it the \textit{trivial loop around $v$}.  

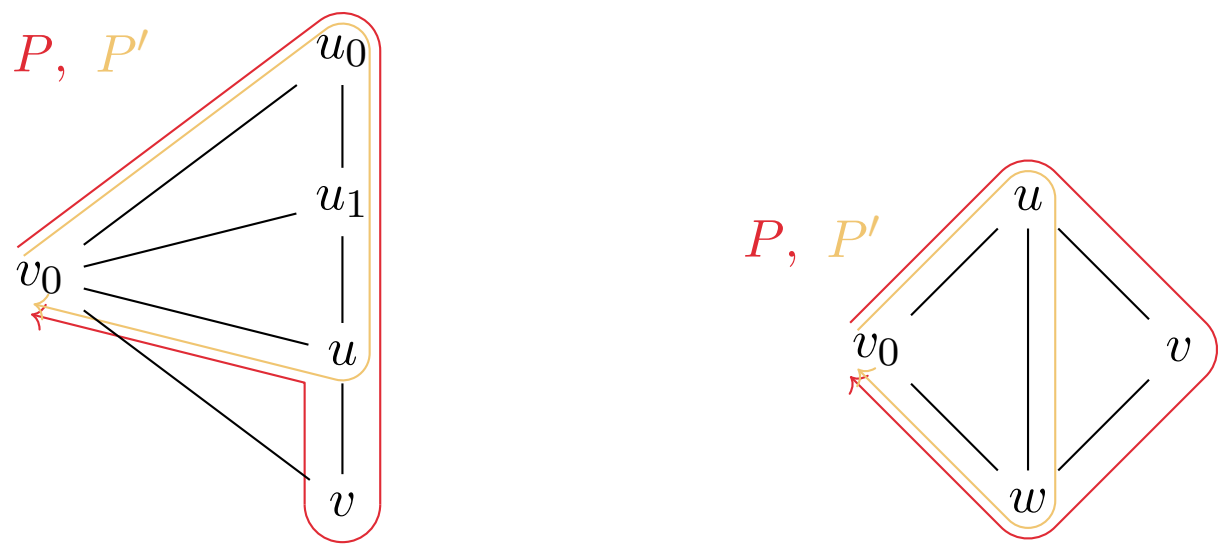
\begin{figure}[h] 
		\centering
			\tikzset{
  offset round/.code=
    \tikz@addmode{%
      \pgfsetroundjoin         \pgfgetpath   \tikz@temp
      \pgfsetpath\pgfutil@empty\pgfoffsetpath\tikz@temp{#1}}}
		\begin{tikzpicture}
  \draw[Mycolor2, -] plot[samples at={(-1,-1.5),(1,0),(1,-1),(1,-2),(1,-3)}] (\x) [offset round=+2.5mm];
  \draw[Mycolor2, ->] plot[samples at={(1,-3),(1,-2),(-1,-1.5)}] (\x) [offset round=+2.5mm];
  \draw[Mycolor2] (0.75,-3) arc (180:360:0.25);
  \draw[Mycolor3, ->] plot[samples at={(-1,-1.5),(1,0),(1,-1),(1,-2),(-1,-1.5)}] (\x) [offset round=+1.8mm];

	\node (a) at (-1,-1.5) {$v_0$};
  \node (b) at (1,0) {$u_0$};
  \node (c) at (1,-1) {$u_1$};
%	\node (d) at (1,-2) {$u_{n-1}$};
  \node (e) at (1,-2) {$u$};
    \node (f) at (1,-3) {$v$};
	%\node[Mycolor1] (n) at (-1.55,0) {$P_0,$};
	\node[Mycolor2] (o) at (-1,-0.05) {$P,$};
	\node[Mycolor3] (p) at (-0.45,0) {$P'$};
  \graph { (a) -- {(b), (c),  (e), (f)}; (b) -- (c) -- (e) -- (f) };
\end{tikzpicture} \hspace*{2cm}
		\begin{tikzpicture}
  \draw[Mycolor2, ->] plot[samples at={(0,-1),(1,0),(2,-1),(1,-2),(0,-1)}] (\x) [offset round=+2.5mm];
  \draw[Mycolor3, ->] plot[samples at={(0,-1),(1,0),(1,-2),(0,-1)}] (\x) [offset round=+1.8mm];

	\node (a) at (0,-1) {$v_0$};
  \node (b) at (1,0) {$u$};
  \node (c) at (1,-2) {$w$};
  \node (d) at (2,-1) {$v$};
	%\node[Mycolor1] (n) at (-1.55,0) {$P_0,$};
	\node[Mycolor2] (o) at (-0.7,-0.3) {$P,$};
	\node[Mycolor3] (p) at (-0.15,-0.25) {$P'$};
  \graph { (a) -- {(b), (c)}; (d) -- {(b), (c)};(b)--(c)};
\end{tikzpicture}
\caption{Visualization of $\BT$ on the left and $\TR$ on the right. \label{fig: BT and TR}}
\end{figure}

\begin{definition}
Fix $v_0$ to be a vertex in $X$.  We define the following symmetric relations between two loops $P, P'$ around $v_0$, see also their visualizations in \pref{fig: BT and TR}.
\begin{itemize}
\item We denote $P \BT P '$ if there exist vertices $u,v$ and paths $Q_0, Q_1$ from $v_0$ to $u$ such that $P = Q_0 \circ (u ; v; u) \circ Q_1^{-1}$ and $P ' = Q_0 \circ (u) \circ Q_1^{-1}$ ($(BT)$ stands for ``backtracking'').
\item We denote $P \TR P '$ if there is a triangle $uvw$ in $X$ and paths $Q_0, Q_1$ from $v_0$ to $u$ and $w$ correspondingly such that $P = Q_0 \circ (u ; v; w) \circ Q_1^{-1}$ and $P '= Q_0 \circ (u ; w) \circ Q_1^{-1}$ ($(TR)$ stands for ``triangle''). We denote by $t_{P,P'} = uvw$ the triangle used for the relation.
\end{itemize}
\end{definition}
	\makeatletter

\begin{remark}
Note that the relations above are not defined to be transitive. Indeed, in the sequel, it will be important to keep track of how many times the relation $\TR$ is used to pass from one loop to another.
\end{remark}

\begin{definition} \label{def: non-abelian cone}
A \emph{$0$-cone} on $X$ is a pair $\Cone_X = \left(v_0,  ( P_u )_{u \in X(0)}\right)$ such that
\begin{enumerate}
\item $v_0 \in X(0)$.
\item For every $u \in X(0), u \neq v_0$,  $P_u$ is a path from $v_0$ to $u$. For $v_0$, $P_{v_0}$ is the trivial loop around $v_0$.
\end{enumerate}

A \textit{$1$-cone} on $X$ is a triple $\Cone_X = \left(v_0,  ( P_u )_{u \in X(0)},  ( T_{uv} )_{uv \in \ordX (1)} \right)$ such that
\begin{enumerate}
\item $\left(v_0,( P_u )_{u \in X(0)}\right)$ is a $0$-cone on $X$.
\item For every $uv \in \ordX (1)$,  $T_{uv}$ is a sequence of loops around $v_0$,  $T_{uv} = (P_0,P_1,\dots,P_m)$ such that the following holds:
\begin{enumerate}
\item $P_0 = P_u \circ (u ; v) \circ P_v^{-1}$.
\item For every $0 \leq i \leq m-1$,  either $P_i \TR P_{i+1}$ or $P_i \BT P_{i+1}$.
\item $P_m$ is the trivial loop around $v_0$.
\end{enumerate}
\end{enumerate}

We refer to the fixed vertex $v_0$ as the apex of the cone. Furthermore, we call $T_{uv} = (P_0,P_1,...,P_m)$ as above a \emph{contraction} and denote the number of triangle relations in the contraction by
$$\vert T_{uv} \vert = \vert \lbrace 0 \leq i \leq m-1 : P_i \TR  P_{i+1} \rbrace \vert . $$
We further define the \textit{radii of the cone} as follows:
\begin{enumerate}
\item $\displaystyle \Rad_0 (\Cone_X) = \sup_{u \in X(0)} \operatorname{length}(P_u)$. 
\item $\displaystyle \Rad_1(\Cone_X) = \sup_{uv \in \ordX (1)} \vert T_{uv} \vert$.
\end{enumerate}
\end{definition}

\begin{definition}
    A topological space $T$ is \emph{$1$-connected} if it is path-connected and each continuous map from the $1$-sphere to $T$, $f:S^1 \to T$ is nullhomotopic, i.e. homotopic to a constant map that maps all points of $S^1$ to a single point $t_0 \in T$. We call a simplicial complex 1-connected if its geometric realization is 1-connected.
\end{definition}
\begin{remark} \label{rem: relations and homotopy}
    Let $P,P'$ be two loops around a vertex $v$ in $X$. If $P \BT P'$ or $P \TR P'$ then the corresponding paths in the geometric realization of $X$ are homotopic. Conversely, any two loops $P,P'$ in $\vert X \vert$ along the edges of $X$ around $v$ which are homotopic can be seen as loops in the abstract simplicial complex and we can find a sequence of loops $P=P_0, \dots, P_m = P'$ such that $P_i \BT P_{i+1}$ or $P_i \TR P_{i+1}$ for all $0 \leq i \leq m-1$.
\end{remark}
\begin{remark} \label{rem: connectednes implies cone}
    Any connected non-empty simplicial complex admits a $0$-cone: pick an arbitrary vertex $v_0$, and by connectivity, we can find a path from any other vertex $u$ to $v_0$. The converse is also true.

    Similarly, any $1$-connected (connected and simply connected) simplicial complex $X$ admits a $1$-cone function. By connectedness we get a $0$-cone. Then each loop $P_u \circ (u;v) P_v^{-1}$ is contractible in the geometric realization of $X$. By Remark \pref{rem: relations and homotopy}, we can realize this homotopy equivalence by a sequence of paths connected by the relations $\BT, \TR$.
\end{remark}

\begin{lemma}
    Let $X$ be a finite $n$-dimensional simplicial complex with 1-cone $\Cone_X = \left(v_0,  ( P_u )_{u \in X(0)},  \left( T_{uv} \right)_{uv \in \ordX (1)} \right)$. Then $\vert X \vert$ is 1-connected and hence $H^0(X,\Lambda)=1$ and $H^1(X,\Lambda)=1$ for any group $\Lambda$.
\end{lemma}
\begin{proof}
    Note that $X$ is connected and hence $\vert X \vert$ is path connected. Hence we can without loss of generality consider a loop around $v_0$ in $\vert X \vert$. This loop is homotopic to a loop along the edges of $X$. Let $P = (v_0;u_1;\dots;u_m;v_0)$ denote this loop. Then 
    \begin{align*} P=& (v_0) \circ (v_0;u_1,\dots, u_m,v_0) \sim_{T_{(u_1,v_0)}} P_{u_1} \circ (u_1;v_0) \circ (v_0;u_1,\dots, u_m,v_0) \\
        &\BT P_{u_1} \circ (u_1,\dots, u_m,v_0) \BT \dots \BT P_{u_1} (u_1;u_2) \circ P_{u_2}^{-1} \circ P_{u_2} \circ  (u_2,\dots, u_m,v_0) \\
        &\sim_{T_{(u_1,u_2)}} P_{u_2} \circ (u_2,\dots, u_m,v_0). \end{align*}
        Here we denoted by $\sim_{T_{(u_i,u_{i+1})}}$ the series of relations that are described by $T_{(u_i,u_{i+1})}$.

    Further iterating this process yields 
    $$P \sim P_{u_m} \circ (u_m;v_0) \sim_{T_{(u_m,v_0)}} (v_0).$$
    Thus, we have shown that $P$ is contractible and hence $\vert X \vert$ is simply connected.
    
    The vanishing of the cohomology groups follows from the fact that cohomology is homotopy invariant.
\end{proof}

A bound on the cone radius gives rise to a bound on the coboundary expansion, see \pref{prop: cone implies CBE} for the detailed statement. Thus, we can intuitively think of coboundary expansion as a quantitative way of (higher) connectedness.

\begin{comment}This result was first shown in \cite[Lemma 1.6,  Lemma 4.3]{DiksteinD2023swap} adapting similar results for coboundary expansion with $\F_2$ coefficients from \cite{KaufmanO2021}. 
\begin{proposition}[{\cite[Lemma 1.6,  Lemma 4.3]{DiksteinD2023swap}}] \label{prop: cone implies CBE}
Let $X$ be a pure $n$-dimensional finite simplicial complex, $n\geq 2$.  Assume that there is a 1-cone  on $X$ such that $\Rad_i (\Cone_X) < \infty, i=0,1$ and that the automorphisms group $\Aut(X)$ acts transitively on $X(n)$. Then for every group $\Lambda$,
$$h^1(X,  \Lambda) \geq \frac{1}{{n+1 \choose 3} \Rad_1 (\Cone_X)}.$$
\end{proposition}
\end{comment}
In particular, the existence of a cone function with bounded cone radius implies coboundary expansion.

Next, we recall a method for bounding the cone radius using the Dehn function of a specific presentation. We start by recalling some facts that describe on the group theoretic level when a coset complex is 1-connected.

\begin{lemma}[{see e.g. \cite[Lemma 3.3]{HarshaS2019}}]\label{prop: CC connected if subgroups generate G}
    Let $G$ be a group with subgroups $\mathcal{H}=(H_0,\dots , H_d)$. Then $\CC(G;\mathcal{H})$ is connected (i.e. the underlying graph is connected) if and only if $G = \langle H_0, \dots, H_d \rangle$.
\end{lemma}

\begin{definition}
    Let $G$ be a group and $H_i = \langle S_i \mid R_i\rangle, i \in I$ be a finite family of subgroups of $G$ given by their presentation. Let $\phi_{ij}: H_i \cap H_j \hookrightarrow H_i$ denote the embedding of $H_i \cap H_j$ into $H_i$ for all $i,j \in I$. The free product of the subgroups $H_i$ amalgamated along their intersections is the group with presentation
    \[\langle \bigcup_i S_i \mid \bigcup_i R_i \cup \{\phi_{ij}(h) = \phi_{ji}(h) \mid \text{ for all }i,j\in I, h \in H_i \cap H_j\} \rangle. \]
\end{definition}
Intuitively, in the above definition, we consider the free product of the $H_i$ but identify the elements in the intersection of the subgroups with their images in both subgroups.

 \begin{proposition}[{\cite[Theorem 2.4]{AbelsHolz}}]\label{prop: CC simply connected if G is amalgam of subgroups}
    Let $G$ be a group with subgroups $\mathcal{H}=(H_0,\dots , H_d)$. Then $\CC(G;\mathcal{H})$ is 1-connected if and only if $G$ is the free product of the subgroups $H_i$ amalgamated along their intersections. 
 \end{proposition}

For most of our work, we consider a specific presentation for a group associated to a coset complex. 

Let $G$ be a group with subgroups $H_i, i \in {I}$, where $I$ is a finite set with $|I| \geq 3$. Assume that $G$ is the free product of the subgroups $H_i$ amalgamated along their intersections. 
For each $i\in I$ let $R_i$ denote the set of all the nontrivial relations in the multiplication table of $H_i$ of length at most $3$.
Then we obtain a presentation of $G$ by considering all the elements of the $H_i$ as generators and $\bigcup_i R_i$ as the set of relations. (We do not need to consider the relations given by the identifications in the definition of the amalgamated product, as we already consider the $H_i$ as subgroups of $G$.) 

\begin{definition} \label{def: dehn function of coset complex}
    Given a 1-connected coset complex $X=\CC(G;(H_i)_{i \in I})$, we denote by $\Dehn_X$ the Dehn function of the presentation of $G$ obtained by considering all the elements of the $H_i$ as generators and $\bigcup_i R_i$ as the set of relations as above. 
\end{definition}

We will see below in \pref{thm: dehn implies cone}, following the work of \cite{KaufmanO2021}, that a bound on the Dehn function of a coset complex implies a bound on its cone radius.

\section{Proof of the Main Theorem} \label{sec:overvew}
In this section we prove \pref{thm: main thm ch 3}, our main theorem on swap coboundary expansion of congruence KMS complexes of type \(A_d\). The statement below specifies the parameter choices and the coefficient groups for which the result holds.
\ydnote{I removed \(T\) since it does not appear in the theorem.}
\begin{theorem} \label{thm: main thm ch 3}
Let $r \in \N$, $r \geq 10^4$. We fix the following parameters:
\begin{itemize}
    %\item $T \geq r^{10}$,
    \item $d \in \N, d \geq r^{55}$,
    \item a desired spectral expansion for the final complex $\frac{r+1}{d-r^{11}-r^{10}-2} \leq \gamma_1 <\frac{1}{r^{11}}$,
    \item $n \in \N,  r^{11}+r^{10}-1 \leq n \leq d- \frac{r+1}{\gamma_1}-1$,
    \item $q_0 \in \N$ large enough with respect to the other parameters, polynomial in $r$. \ydnote{Can we just set it?}
\end{itemize}
Then there exists $\beta>0$ which is of order $\frac{1}{r^{34}}$ such that 
    \begin{itemize}
        \item for all prime powers $q=p^m \geq q_0$, 
        \item for all groups $\Lambda$ without non-trivial elements of order $p$, 
        \item for all irreducible polynomials $f \in \F_q[t]$ with $\deg f \geq 2$, 
    \end{itemize}
    the $n$-skeleton of the KMS complex $X(A_d,\F_q, f)$ is a $(\beta,r,\Lambda)$-swap coboundary expander and $\gamma_1$-local spectral expander. 
\end{theorem}

\begin{corollary}
    In the above setting, the KMS complex $X(A_d,\F_q, f)$ is a $(\beta,r,\Lambda)$-swap coboundary expander, proving \pref{thm:main}. 
\end{corollary}
Note that it the full $r$-faces complex is in general not a local spectral expander, see \pref{obs: LSE of full faces complex}.

\begin{proof}
    We need to show that $F^rX$ is a $(\beta,\Lambda)$ 1-coboundary expander. Now, 1-coboundary expansion is fully determined by the 2-skeleton of the corresponding complex (using our convention about the weights of a skeleton). By \pref{lem: skeleton of faces complex} below, we have
    \[(F^rX)^{\leq 2} = F^r(X^{\leq 2r+r+2}) = F^r((X^{\leq n})^{\leq 2r+r+2})=(F^r(X^{\leq n}))^{\leq 2}\]
    given $2r+r+2 \leq  n$ which holds with our choice of constants.
\end{proof}

\paragraph{Roadmap} The proof proceeds in four steps. The first three are standard, the fourth is the main new step of this work. We explain the reductions in this section, stating the main technical ingredients whose proofs appear in later sections.

\begin{enumerate}
    \item \textbf{Reduction to colored faces complexes (\pref{subsec:overview-colored-reduction}).} We first reduce swap coboundary expansion to coboundary expansion of colored faces complexes for a sufficiently large fraction of choices of colors. We identify a collection of ``nice'' choices for which we will prove expansion. The color sampling lemma (\pref{lem:most-colors-are-nice}), proved in \pref{sec: color sampling lemma}, shows that a random choice is nice with high probability.

    \item \textbf{Reduction to vertex links (\pref{subsec:overview-links}).} For each nice choice of colors, we use the local-to-global theorem of \cite{DiksteinD2023swap} to reduce the expansion bound to coboundary expansion of every vertex link. These links have a simpler algebraic description, which allows us to analyze them in the remaining steps.

    \item \textbf{Reduction to diameter and Dehn bounds (\pref{subsec:overview-cones-dehn}).} We describe each vertex link as a coset complex and use the cone method to reduce its coboundary expansion to constructing a cone of bounded radius. In turn, this radius can be bounded using the diameter and the Dehn function of the coset complex, which measures how many applications of the defining relations are needed to reduce a word representing the identity. The precise connection between cones and Dehn functions is proved in \pref{app:dehn-equivalence}, extending the result of \cite{KaufmanO2021}.

    \item \textbf{Decomposition into bins (\pref{subsec:overview-bins}).} The group underlying each link decomposes as a direct product of smaller upper-triangular matrix groups, corresponding to intervals of colors called bins. Our main new ingredient is the decomposition lemma, proved in \pref{sec:dirprod}, which reduces the Dehn bound to separate bounds for these factors, with an additional cost polynomial in \(r\) and the word length. For each factor, we obtain a linear Dehn bound whose constant is independent of the matrix size. This uses the matrix-group theorem proved in \pref{sec: Unip sec}; its application to the factors, together with the required diameter bound, is proved in \pref{sec: bounding the Dehn function}.
\end{enumerate}

Finally, in \pref{subsec:overview-main-proof}, we combine these ingredients and track the parameters to complete the proof of \pref{thm: main thm ch 3}.

\subsection{Step 1: Reduction to colored faces complex} \label{subsec:overview-colored-reduction}
The following result by Dikstein-Dinur reduces the question to showing coboundary expansion for a large set of colored faces complexes, instead of considering the full faces complex $F^rX$. This has the advantage that the colored faces complex of a coset complex can again be described as a coset complex. 
\begin{theorem}[{\cite[Lemma 8.2]{DiksteinD2023swap}}] \label{thm: DD-reduction to colored faces}
    Let \(T, d,r \in \NN\), such that \(3 \leq rT \leq (d+1)^{1/5}\). Let \(\beta, p, \lambda \in (0,1]\). Let \(\Lambda\) be any group. Let \(X\) be a \((d+1)\)-partite simplicial complex such that 
    \begin{enumerate}
        \item \(F^r X\) is a \(\gamma\)-local spectral expander for \(\gamma < \frac{1}{rT}\), and
        \item \[ \Prob[\cI \in \binom{[d]}{r+1}^{[T]}]{h^1(F^{\cI} X, \Lambda) \geq \beta \land H^1(F^{\cI}X ,\Lambda) = 1} \geq p.\]
    \end{enumerate}
    Then \(H^1(F^r X, \Lambda) = 1\) and \(h^1(F^r X, \Lambda) \geq \Omega(p^2 \beta)\).
\end{theorem}
This statement is proven in \cite[Section 8.3]{DiksteinD2023swap} for the case where \(T = \sqrt{r+1}\) and for \(X\) which is a quotient of the affine building of type \(\tilde{A}\), but no properties of the building quotient were used in the proof other than the ones that appear in the theorem statement.

\paragraph{Proof intuition for \pref{thm: DD-reduction to colored faces}}
The proof of \pref{thm: DD-reduction to colored faces} in \cite[Section 8.3]{DiksteinD2023swap} is based on extending a coboundary from a colored faces complex to the entire faces complex. Given a \(1\)-cochain \(\phi\) on \(Z=F^rX\), the goal is to show that if \(\phi\) violates the cocycle condition on few triangles, then it is close to a coboundary.

The first step in the proof is to find a collection \(\cI\) of disjoint colors such that \(Y=F^{\cI}X\) has good coboundary expansion, and such that \(\phi\) violates the cocycle condition on relatively few triangles involving vertices of \(Y\). This is obtained by averaging over choices of colors, using the assumption that a \(p\)-fraction of the colored faces complexes are coboundary expanders. The argument also considers larger collections of colors, whose expansion is deduced from that of smaller restrictions.

Having found such a subcomplex \(Y\), its coboundary expansion implies that there is a \(0\)-cochain \(g\) on its vertices such that \(\phi=d_0g\) on almost all edges of \(Y\). The next step is to extend \(g\) to the remaining vertices. Every neighbor \(u\in Y\) of a vertex \(v\notin Y\) suggests a value for \(g(v)\), determined by \(g(u)\) and \(\phi(uv)\). The key observation is that two such suggestions agree whenever they come from vertices \(u,w\) for which the edge \(uw\) and triangle \(uvw\) satisfy the relevant coboundary and cocycle conditions. The spectral expansion of the links ensures that most suggestions are consistent, allowing \(g(v)\) to be defined by the most common suggestion.

Finally, the values of \(g\) are shown to agree with \(\phi\) on most edges whose endpoints both lie outside \(Y\). Indeed, for such an edge \(vw\), a triangle \(uvw\) with \(u\in Y\) allows one to recover \(\phi(vw)\) from the other two edges, provided the cocycle condition holds. The remaining technical part of the proof consists of showing that all these comparisons involve sufficiently many triangles, and bounding the resulting errors using local spectral expansion. This yields the claimed coboundary expansion bound.

To be able to apply \pref{thm: DD-reduction to colored faces}, we need to show local spectral expansion for the faces complex.

% \begin{lemma}
% Let \(r,n,d \in \N\), with \(r\geq 1\) and \(n\geq 2r+1\), and let \(\gamma_1 \in (0,1)\). Suppose that
% \[
% d \geq n+\frac{r+1}{\gamma_1}+1.
% \]
% There exists \(q_0 \in \N\), polynomially depending on \(r,n,d,\gamma_1\), such that for every prime power \(q=p^m \geq q_0\) and every irreducible polynomial \(f \in \F_q[t]\) with \(\deg f \geq 2\), the \(r\)-faces complex of the \(n\)-skeleton of \(X(A_d,\F_q,f)\), namely
% \[
% F^r\bigl(X(A_d,\F_q,f)^{\leq n}\bigr),
% \]
% is a \(\gamma_1\)-local spectral expander.
% \end{lemma}

\begin{lemma}  \label{lem: LSE for faces complex of skeleton}
Let $r,n,d\in\N$ satisfy $2r+1\leq n\leq d$, and let $\gamma_1\in(0,1)$. Suppose that \[d\geq n+\frac{r+1}{\gamma_1}+1.\]
For every prime power $q> 16$ and every irreducible $f\in\F_q[t]$ of degree at least two, the complex $F^r(X^{\leq n})$, where $X=X(A_d,\F_q,f)$, is a $\gamma_1$-local spectral expander.
\end{lemma}
\ydnote{Old proof commented out due to stronger local spectral expansion.}
\begin{proof}
By \pref{thm: KMS complexes are local spectral}, the complex $Y=X^{\leq n}$ is a $\frac{1}{d-n+1}$-two-sided local spectral expander. The graph underlying each link of $F^rY$ is the swap walk between disjoint $(r+1)$-vertex faces in a link of $Y$. The swap-walk bound of \cite[Lemma 3.3]{GurLL2022} therefore bounds its nontrivial eigenvalues in absolute value by \[\frac{r+1}{d-n+1}\leq\gamma_1.\]
\end{proof}

Thus, given that we choose the parameters $(r,n,d,\gamma, q_0)$ as in \pref{thm: main thm ch 3} and setting $X = X(A_d,\F_q,f)$ and $Y$ the $n$-skeleton of $X$, we get that $F^rY$ is a $\gamma_1$-local spectral expander. Given our choice of $r,n,d$ the colored faces complexes of $X$ and $Y$ are the same, as the following lemma shows.

\begin{lemma} \label{lem: skeleton of faces complex}
    Let $X$ be a $(d+1)$-partite complex, and fix integers $1\leq r,n< d, 1 \leq T \leq \frac{d+1}{r+1}$. Let $\cI \in \binom{[d]}{r+1}^{[T]}$ be a set of $T$ disjoint size $r+1$ color sets. Recall that $X^{\leq n}$ denotes the $n$-skeleton of $X$. Then 
    \begin{enumerate}
        \item \(\left(F^rX\right)^{\leq n} = F^r(X^{\leq nr+n+r})\),
        \item \(\left(F^{\cI}X\right)^{\leq n} = F^{\cI}(X^{\leq nr+n+r})\),
        \item if $n\geq T(r+1)-1$ then $F^{\cI}X = F^{\cI}(X^{\leq n}).$
    \end{enumerate}
    In all three cases, the equility also respects the weights on the complexes, using the convention introduced in \pref{subsec: simplicial complex} for the weights on skeleta. 
\end{lemma}

\begin{proof}
    Following the definition of $F^rX$ we have for a set $\{\sigma_0,\dots,\sigma_i\}$ with $\sigma_j \in X(r)$ pairwise disjoint
    \[\{\sigma_0,\dots,\sigma_i\} \in F^r X(i) \iff \cup_j \sigma_j \in X(ir+r+i)\]
    where we use that, since the $\sigma_j$ are pairwise disjoint, $\vert \cup_j \sigma_j \vert = \sum_j \vert \sigma_j \vert = (i+1)(r+1)$, and hence if $\cup_j \sigma_j \in X$ then it will be of dimension $(ir+i+r)$.

    Now when considering the $n$-skeleton of $F^rX$ we have
    \[\{\sigma_0,\dots,\sigma_i\} \in (F^r X)^{\leq n} \iff \cup_j \sigma_j \in X(ir+r+i) \land i \leq n \iff \cup_j \sigma_j \in X^{\leq nr+r+n}  \]
    which shows point 1.. The proof for point 2. is analogues. 

    Point 3. follows directly from point 2. as $\dim F^{\cI}X = T-1$, hence $F^{\cI}X=(F^{\cI}X)^{\leq T-1} = F^{\cI}(X^{\leq (T-1)r+(T-1)+r})=F^{\cI}(X^{\leq T(r+1)-1})$.

    Regarding the weights, note that for each $0\leq i \leq n$ by our convention, the weight on the $i$-simplices of $F^rX$ and $(F^rX)^{\leq n}$ coincide. 
    On the other hand, we have that $w_X,i = w_{X^{\leq nr+n+r}},i$ for all $0 \leq i \leq nr+n+r$ and thus they induce the same weight on the $r$-faces complex. Together we get
    \[w_{(F^rx)^{\leq n},i} = w_{F^rX,i} = w_{F^r(X^{\leq nr+n+r}),i}\]
    for all $0\leq i \leq n$ as required. The argument for the other two points is analogous. 
\end{proof}

Note that a statement similar to part 3. for $F^rX$ is not very useful, since $\dim F^rX = \lfloor \frac{d+1}{r+1} \rfloor -1$ and hence we would require $n \geq (\lfloor \frac{d+1}{r+1} \rfloor -1 +1)(r+1)-1$, which in the case where $d+1$ is divisible by $r+1$ results in requiring $n \geq d$ for which the statement is trivial.

 Using the above lemma, we will from now on consider $F^{\cI}X$ instead of the colored faces complex of the skeleton.
For part 2. of \pref{thm: DD-reduction to colored faces} we show coboundary expansion for all colored faces complexes given by ``nice'' sets of colors. Furthermore, we show that sets of colors are nice with high probability.

\ydnote{added somoe prose to explain the cyclic stuff:}
The conditions defining a nice collection will describe the matrix blocks that appear in vertex links. Before defining them, let us explain why these blocks correspond to intervals in a linear order, even though the root indices of the KMS complex are cyclic. Let \(X=X(A_d,\F_q,f)\), let \(s\in X(r)\) represent a vertex of \(F^{\cI}X\), and write \(J=\col(s)\). Choose a vertex \(u\in s\) and write \(s=s'\cup\{u\}\). A collection of vertices belongs to \((F^{\cI}X)_s\) exactly when its members have colors in \(\cI\setminus\{J\}\), are pairwise disjoint, and their union together with \(s\) is a face of \(X\). Equivalently, their union is a face of \(X_s\). Since \(X_s=(X_u)_{s'}\), we obtain \[(F^{\cI}X)_s=F^{\cI\setminus\{J\}}(X_s)=F^{\cI\setminus\{J\}}\bigl((X_u)_{s'}\bigr).\] Thus we can study the same link by first passing to \(X_u\) and then taking the link of \(s'\) there. The vertices of the resulting faces complex still represent faces with \(r+1\) vertices.

As described in \pref{subsec:KMS-complexes}, the vertex link \(X_u\) is a coset complex of upper triangular unipotent matrices. If \(u\) has color \(j_*\), its simple root indices follow the cyclic order \(j_*+1,\dots,d,0,\dots,j_*-1\), which becomes a path after deleting \(j_*\). We relabel this path as \(1,\dots,d\). The remaining \(r\) colors of \(s'\) divide it into \(r+1\) intervals, and the root groups in these intervals act on disjoint diagonal blocks. These are the bins used in the definition below.

\ydnote{added:}For the rest of the paper we will need a parameter polynomial in \(r\) to denote the number of colors we want to restrict to. We will keep this parameter abstract in statements but can set it to \(T=r^{10}\) for concreteness.
\begin{definition}\label{def:nice-colors}
\begin{enumerate}
    \item Let \(J\subseteq[d]\) have size \(r+1\). Choose \(j_0\in J\) and cyclically relabel all colors so that \(j_0\) becomes \(0\) and the colors following it become \(1,\dots,d\). Write the relabeled set as \(J=\{0=j_0<j_1<\cdots<j_r\}\), and set \(j_{r+1}=d+1\). We define the \emph{bins of \(J\)}, or \emph{\(J\)-bins}, by \[B_\ell=\{j_\ell+1,\dots,j_{\ell+1}-1\},\qquad 0\leq\ell\leq r.\] We write \(B=[a,b]\) as a shorthand for \(B=\{a,a+1,\dots,b\}\). When considering the bins of a fixed \(J\), we apply the same relabeling to every color set in \(\cI\). Choosing a different \(j_0\in J\) only cyclically permutes the bins and translates their coordinates, preserving the order within each bin. Thus the conditions below do not depend on this choice. Empty bins contribute trivial factors to the group decomposition and are omitted from the niceness requirements.
    \item Let \(\mathcal{I}\in\binom{[d]}{r+1}^{[T]}\) be a set of \(T\) disjoint color sets. We say \(\mathcal{I}\) is \emph{nice} if the following holds for every \(J\in\mathcal{I}\), in the path order chosen for \(J\).
    \begin{enumerate}
        \item Every bin \(B=[a,b]\) has size \(|b-a|>100\).
        \item For every bin \(B=[a,b]\) of \(J\) there exists a color \(I' \in \mathcal{I} \setminus \set{J}\) such that \(I' \cap B \subseteq [a, a+ \lfloor \frac{b-a}{5}\rfloor ]\).
        \item For every bin \(B=[a,b]\) of \(J\) there exists a color \(I' \in \mathcal{I} \setminus \set{J}\) such that \(I' \cap B \subseteq [b- \lfloor \frac{b-a}{5} \rfloor,b]\).
        \item For every bin \(B=[a,b]\) of \(J\) there exists a color \(I' \in \mathcal{I} \setminus \set{J}\) such that \(I' \cap B \subseteq [a+2\lceil\frac{b-a}{5}\rceil ,a+3\lfloor \frac{b-a}{5} \rfloor]\).
        \item For every pair of \(J\)-bins \(B=[a,b],B'=[a',b']\) and every \(\ell,\ell' \in \set{0,1,2,3}\) there exists some \(I' \in \mathcal{I} \setminus \set{J}\) such that \[I' \cap [a + \ell  \lfloor \frac{b-a}{5} \rfloor ,\min \lbrace a + (\ell+2) \lceil \frac{b-a}{5} \rceil +4,b \rbrace] = \emptyset, \text{ and}\]
        \[I' \cap [a' + \ell '  \lfloor \frac{b'-a'}{5} \rfloor ,\min \lbrace a' +  (\ell '+2) \lceil \frac{b'-a'}{5} \rceil +4,b' \rbrace] = \emptyset.\]
\ydnote{Changed to accomodate the cyclic bins +no small bins. I also changed last item to have two bins from the same J.}        %\item For every pair of bins \(B,B'\) and every \(\ell,\ell' \in \set{0,1,\dots,8}\) there exists some \(I' \in \mathcal{I} \setminus \set{I}\) such that \(I' \cap [a + \ell \frac{b-a}{10},a+(\ell+2) \frac{b-a}{10}] = \emptyset\) and \(I' \cap [a + \ell' \frac{b-a}{10},a+(\ell'+2) \frac{b-a}{10}] = \emptyset\) 
        %\ionote{I can make do with the above, but to simplify things I actually want something stronger (I think this stronger version should follow from the condition for the previous items). If we ignore rounding issue, what I want is the following:
        %For every pair of bins \(B,B'\) and every \(\ell,\ell' \in \set{0,1,2,3}\) there exists some \(I' \in \mathcal{I} \setminus \set{I}\) such that \(I' \cap [a + \ell \frac{b-a}{5},a+(\ell+2) \frac{b-a}{5}] = \emptyset\) and \(I' \cap [a + \ell' \frac{b-a}{5},a+(\ell'+2) \frac{b-a}{5}] = \emptyset\). With rounding issues (if we want to be very precise):
        %For every pair of bins \(B,B'\) and every \(\ell,\ell' \in \set{0,1,2,3}\) there exists some \(I' \in \mathcal{I} \setminus \set{I}\) such that \(I' \cap [a + \ell  \lfloor \frac{b-a}{5} \rfloor ,\min \lbrace a + (\ell+2) \lceil \frac{b-a}{5} \rceil +4,b \rbrace] = \emptyset\) and \(I' \cap [a + \ell '  \lfloor \frac{b-a}{5} \rfloor ,\min \lbrace a +  (\ell '+2) \lceil \frac{b-a}{5} \rceil +4,b \rbrace] = \emptyset\)}
    \end{enumerate}
    \end{enumerate}
\end{definition}

Indeed, sufficiently many sets of colors are nice: 
\begin{lemma} \label{lem:most-colors-are-nice}
    Let \(r\geq 10^4, T \geq r^{10}\) and \(d \geq 100(r+1)^2T^2\) then with probability \(\geq 0.9\) a randomly selected set of $T$ disjoint colors \(\mathcal{I} \in \binom{[d]}{r+1}^{[T]}\) is nice. 
\end{lemma}

We prove this lemma in \pref{sec: color sampling lemma}.

The main next step is now to show that the nice colored faces complexes are coboundary expanders.
\begin{proposition} \label{prop: nice colored faces are CBE}
    Let $r,T,d\in \N$ such that $T\geq 4$ and $(r+1)T \leq d+1$. There exist $\beta_1 >0, q_0 \in \N$ such that for all prime powers $q =p^m \geq q_0$, all groups $\Lambda$ without non-trivial elements of order $p$, all irreducible polynomials $f\in \F_q[t]$ of degree $\geq 2$ and all nice sets of $T$ disjoint colors $\cI \in \binom{[d]}{r+1}^{[T]}$, we have
    \[h^1(F^{\cI}X, \Lambda) \geq \beta_1 \text{ and } H^1(F^{\cI}X, \Lambda)=1\]
    where $X = X(A_d,\F_q,f)$.
\end{proposition}
We start by showing the vanishing of the first cohomology group.
\begin{lemma} \label{lem: vanishing of cohom for colored faces complex}
     Let $r,d,T \in \N, d,T \geq 3$. For all prime powers $q =p^m \geq 5$, all groups $\Lambda$ without non-trivial elements of order $p$, all irreducible polynomials $f\in \F_q[t]$ of degree $\geq 2$ and all nice sets of disjoint colors $\cI \in \binom{[d]}{r+1}^{[T]}$, we have
     \[H^1(F^{\cI}X(A_d,\F_q,f),\Lambda) =1.\]
\end{lemma}

\begin{proof}
Let $X = X(A_d,\F_q,f)$. Recall that 
\[X = \CC(G;(H_i^f)_{i \in [d]})\]
where $G = \Sl_{d+1}(\F_q[t]/(f))$ and $H_i^f = \pi_f(U_{[d] \setminus \{i\}})$ with $\pi_f: \Sl_{d+1}(\F_q[t]) \to \Sl_{d+1}(\F_q[t]/(f))$, $U_i=\langle e_{i,i+1}(\lambda) \mid \lambda \in \F_q \rangle \leq \Sl_{d+1}(\F_q[t]), U_0 = \langle e_{d+1,1}(\lambda t) \mid \lambda \in \F_q \rangle$ and $U_J = \langle U_i \mid i \in J\rangle$ for al $J \subseteq [d]$.
Let $\cI \in \binom{[d]}{r+1}^{[T]}$ be the set of disjoint colors from above. For $J \in \cI$ we set $H_J^f = \bigcap_{i \in J} H_i^f$. Then by \pref{lem: colored faces coset complex} we can describe $F^{\cI}X$ as
\[F^{\cI}X = \CC(G;(H_J)_{J\in \cI}). \]
We want to apply basically the same argument as in \cite{OppenheimV2025}, thus we also refer to that work for more details on the terminology.

We can see each such complex (i.e. for each choice of $f$) as a quotient of the same underlying infinite complex, which we describe now. Let $\KMS{d}{q} \leq \Sl_{d+1}(\F_q[t])$ be the corresponding concrete KMS group, denote $H_i =U_{[d] \setminus \{i\}}$, and consider the infinte complex
\[Y = \CC(\KMS{d}{q};(H_i)_{i \in [d]}).\]
We refer to this as fundamental KMS complex. Let $\varphi_f$ denote the map $\pi_f$ restricted to $\KMS{q}{d}$ and let $N_f = \ker \varphi_f$. In \cite[Lemma 3.5]{OppenheimV2025}, it was shown that $X \cong N_f \backslash Y$. Analogously, it holds that $F^{\cI}X \cong N_f \backslash F^{\cI}Y$.
We want to apply the following result. 
\begin{theorem}[\cite{KaufmanOW2025}] \label{thm: homology of quotients}
Let $X$ be a connected pure 2-dimensional simplicial complex and $N$ be a group acting on $X$. Assume $p: X \to N\backslash X$ is rigid, i.e. it is simplicial and for all $\sigma \in X$ we have $\lvert \sigma \rvert = \lvert p(\sigma)\rvert$. For every group $\Lambda$, if $H^1(X,\Lambda) = 0$ and $H^1(N,\Lambda) = 0$ then $H^1(N \backslash X, \Lambda)=0$ (the cohomology of $N$ is taken with the trivial action of $N$ on $\Lambda$). 
\end{theorem}

In \cite[Cororllary 3.7]{OppenheimV2025} it was shown that if $\Lambda$ has no non-trivial elements of order $p$ (recall $q=p^m$), then $H^1(N_f,\Lambda) =1$. Thus, it remains to show that $H^1(F^{\cI}Y, \Lambda) = 1$ which we will achieve by showing that $F^{\cI}Y$ is 1-connected using  \pref{prop: CC simply connected if G is amalgam of subgroups}. We know by the work of \cite{devillers_simple_2007,Abramenko-TwinBuildings} (see also \cite[Proposition 3.8]{GraveDePeraltaV2025}) that $\KMS{d}{q}$
    is the free product of the subgroups $U_{\{i,j\}}, \{i,j\} \subseteq [d]$ amalgamated along their intersections, given  $d \geq 3, q \geq 5$. 
Hence, it remains to show that for each $\left\{ i,j \right\} \subseteq [d]$ there exists an $L \in \cI $ such that $\left\{ i,j \right\} \subseteq [d]\setminus L$ (and hence $U_{\{i,j\}} \subseteq H_L = \cap_{i \in L}U_{[d]\setminus \{i\}} = U_{[d]\setminus L} $. This clearly holds if $\vert \cI \vert \geq 3$, as all sets in $\cI$ are disjoint.    
\end{proof}
We further reduce the task of bounding the cocycle expansion constant in the next step. 
\subsection{Step 2: Reduction to the links via the local-to-glocal machinery} \label{subsec:overview-links}
To show coboundary expansion of the nice colored faces complexes, we reduce the question to showing coboundary expansion for all vertex links via the local-to-global result by Dikstein-Dinur.
\begin{theorem}[{\cite[Theorem 8]{DiksteinD2023cbdry}}]
\label{thm: local to global}
Let $0 \leq \lambda <1,  \beta >0$ be constants and $\Lambda$ be a group.  For a finite pure $n$-dimensional simplicial complex $Y$ with $n \geq 3$,  if
\begin{itemize}
    \item for every vertex $v \in Y(0)$,  $h^1 (Y_{v}, \Lambda) \geq \beta$ and $H^1(Y_{v},\Lambda)$ is trivial,
    \item and $Y$ is a $\lambda$-local spectral expander,  
\end{itemize}
    then
\[h^1 (Y, \Lambda) \geq \frac{(1-\lambda) \beta}{24} - e \lambda,\]
where $e \approx 2.71$ is the Euler constant.
\end{theorem}
To apply this theorem, we first of all need to show local spectral expansion for $F^{\cI}X$.
% \begin{lemma} \label{lem: colored faces is loc spec}
%      Let $r,d \in \N, d \geq 3r$ and $T\in \N, T\geq 4$. For every $\gamma_2 \in (0,1)$ there exists $q_0 \in \N$ such that for all prime powers $q =p^m \geq q_0$, all irreducible polynomials $f\in \F_q[t]$ and all nice sets of $T$ disjoint colors $\cI \in \binom{[d]}{r+1}^{[T]}$, we have that $F^{\cI}X(A_d,\F_q,f)$ is a $\gamma_2$-local spectral expander. 
% \end{lemma}
\ydnote{Modified proof of the following lemma, old proof is commented out:}
\begin{lemma} \label{lem: colored faces is loc spec}
Let $X=X(A_d,\F_q,f)$, where $q>16$ and $f$ is irreducible of degree at least two. Let $\cI$ consist of $T\geq 3$ pairwise disjoint color sets, each of size $r+1$. Then $F^{\cI}X$ is a $\frac{2(r+1)}{\sqrt q-2}$-local spectral expander.
\end{lemma}

\begin{proof}
We first bound the bipartite walk between any two specified colors in any link $X_\tau$. Write $D=d-|\tau|$ and $\eta_q=2/(\sqrt q-2)$. If this bipartite walk has a nontrivial singular value $s$, extending its centered singular vectors by zero to the other colors gives a mean-zero function on $X_\tau$ with Rayleigh quotient $s/D$. By \pref{thm: KMS complexes are local spectral}, this is at most $\eta_q/D$, so $s\leq\eta_q$.

Thus every two-color walk in every link has nontrivial singular values at most $\eta_q$. By \cite[Lemma 3.3]{GurLL2022}, the walk between faces on two disjoint color sets $I,J$ has nontrivial singular values at most $\sqrt{|I||J|}\eta_q$. For two members of $\cI$, this is $(r+1)\eta_q=\lambda$.

The same bound holds after fixing any face of $F^{\cI}X$, since this amounts to passing to a link of $X$. In each such link, decompose a function into its mean on each remaining color class and a function with mean zero on every class. The first subspace carries the walk on a complete graph, whose nonconstant eigenvalues are negative. On the second subspace, averaging the bipartite walks between classes gives operator norm at most $\lambda$. Hence every link has second eigenvalue at most $\lambda$.
\end{proof}

% The fact that $F^{\cI}X$ is a local spectral expander follows from the fact the KMS complex $X$ is a local spectral expander and the following lemma, which was proven in \cite[Theorem 7.1]{DiksteinD2019} and similar statements were independently shown in \cite{AlevJT2019} and \cite[Lemma 3.3]{GurLL2022}.
% \begin{lemma}
%     Let $X$ be a finite, pure $d$-dimensional simplicial complex, $1 \leq r \leq \frac{d}{3}, T\geq 2$, $\cI \subset \binom{[d]}{r+1}^{[T]}$ be a set of $T$ disjoint colors and $\gamma \in (0,1)$. If $X$ is a $\gamma$-local spectral expander then $F^{\cI}X$ is a $r\gamma$-local spectral expander.
% \end{lemma}
% By \pref{thm: KMS complexes are local spectral}, we have that for every $d\in \N$, $\gamma \in (0,1)$ there exists a $q_0 \in \N$ such that for all prime powers $q \geq q_0$ and for every irreducible polynomial $f\in \F_q[t]$ the KMS complex $X(A_d,\F_q,f)$ is a $\gamma$-local spectral expander. This proves \pref{lem: colored faces is loc spec}.

Finally, we need to show that the links are coboundary expanders. 
\begin{proposition} \label{prop: cbe of links}
     Let $r,d \in \N, d \geq 3$ and $T\in \N, T\geq 4$. There exist $\beta_2 >0, q_0 \in \N$ such that for all prime powers $q =p^m \geq q_0$, all groups $\Lambda$ without non-trivial elements of order $p$, all irreducible polynomials $f\in \F_q[t]$ of degree $\geq 2$ and all nice sets of $T$ disjoint colors $\cI \in \binom{[d]}{r+1}^{[T]}$, we have for all vertices $v\in F^{\cI}X(0)$
     \[h^1(F^{\cI}X_{v},\Lambda) \geq \beta_2 \text{ and } H^1(F^{\cI}X_{v},\Lambda) =1\]
     where $X = X(A_d,\F_q,f)$.
\end{proposition}
In the next step, we further reduce the statement of this proposition using the cone method.

\subsection{Step 3: Coboundary expansion of the links via the cone method} \label{subsec:overview-cones-dehn}
The next step is to show that each link $F^{\cI}X_{v}$ admits a 1-cone function and bounding its radius. Then the following provides a lower bound on the cocycle expansion constant. 
\begin{proposition}[{\cite[Lemma 1.6,  Lemma 4.3]{DiksteinD2023swap}}] \label{prop: cone implies CBE}
Let $Y$ be a pure $d$-dimensional finite simplicial complex, $d\geq 2$.  Assume that there is a 1-cone  on $Y$ such that $\Rad_i (\Cone_Y) < \infty, i=0,1$ and that the automorphisms group $\Aut(Y)$ acts transitively on $Y(d)$. Then for every group $\Lambda$,
$$h^1(Y,  \Lambda) \geq \frac{1}{\binom{d+1}{3} \Rad_1 (\Cone_Y)}.$$
\end{proposition}
Note that any coset complex satisfies the symmetry assumption about the action of the automorphism group. %Hence, also 

Thus, our goal is to show the following.
\begin{proposition} \label{prop: restatement bound on 1-cone}
     Let $r,d \in \N, d \geq 3$, and $T\in \N, T\geq 4$. There exist $\beta_3 >0, q_0 \in \N$ such that for all prime powers $q =p^m \geq q_0$, all irreducible polynomials $f\in \F_q[t]$ of degree $\geq 2$of degree $\geq 2$ and all nice sets of $T$ disjoint colors $\cI \subseteq \binom{d+1}{r+1}$, we have for all vertices $v\in F^{\cI}X(0)$ that $(F^{\cI}X)_{v}$ admits a 1-cone function $\Cone_v$ with $\Rad_1(\Cone_v) \leq \beta_3$ (again $X = X(A_d,\F_q,f)$).
\end{proposition}

To show the existence of a cone function, we have a closer look at the structure of links. 

Let $X = X(A_d,\F_q,f)$. For every vertex $v \in F^{\cI}X(0)$ we can describe $F^{\cI}X_{v}$ explicitly as coset complex. Recall that 
\[X = \CC(G;(H_i^f)_{i \in [d]})\]
where $G = \Sl_{d+1}(\F_q[t]/(f))$ and $H_i^f = \pi_f(U_{[d] \setminus \{i\}})$ where $\pi_f: \Sl_{d+1}(\F_q[t]) \to \Sl_{d+1}(\F_q[t]/(f))$, $U_i=\langle e_{i,i+1}(\lambda) \mid \lambda \in \F_q \rangle \leq \Sl_{d+1}(\F_q[t]), U_0 = \langle e_{d+1,1}(\lambda t) \mid \lambda \in \F_q \rangle$ and $U_J = \langle U_i \mid i \in J\rangle$ for al $J \subseteq [d]$.
Let $\cI \in \binom{[d]}{r+1}^{[T]}$ be the set of disjoint colors from above. For $J \in \cI$ we set $H_J=H_J^f = \bigcap_{i \in J} H_i^f$. Then by \pref{lem: colored faces coset complex} we can describe $F^{\cI}X$ as
\[F^{\cI}X = \CC(G;(H_J)_{J\in \cI}). \]
Thus, for a vertex $gH_J \in F^{\cI}X$, we can describe its link again as coset complex using \pref{lem: link in coco} as follows.
\[F^{\cI}X_{gH_J} \cong \CC(H_J;(H_J \cap H_I)_{I \in \cI\setminus \{J\}}). \]
We can now use the fact (\cite[Lemma 4.2,4.6]{GraveDePeraltaV2025}) that $H_J = \pi_f(U_{[d]\setminus J}) \cong U_{[d]\setminus J}$ and get
\[(F^{\cI}X)_{gH_J} \cong \CC(U_{[d]\setminus J};(U_{[d]\setminus (J\cup I)})_{I \in \cI\setminus \{J\}}). \]

\ydnote{added this explanaition:}
For each fixed \(J\), we now use the upper triangular matrix model obtained by first passing to \(X_u\), with \(\col(u)=j_*\in J\), as explained before \pref{def:nice-colors}. We relabel the color sets and root groups in the path order chosen there, retaining the notation \(U_i\) for the corresponding simple root groups. Thus \(0\in J\), and for \(1\leq i\leq d\) the group \(U_i\) consists of the matrices \(e_{i,i+1}(\lambda)\), \(\lambda\in\F_q\). We use these coordinates whenever \(J\) is fixed in the rest of this section and in \pref{sec:dirprod} and \pref{sec: bounding the Dehn function}.

\begin{corollary}\label{cor: bound on cone radius gives cbe our case}
    Every link \((F^{\cI}X)_{v}\) is a (\(\vert \cI \setminus\{J\} \vert -1 = T-2 \))-dimensional coset complex. Hence its automorphism group acts transitively on the maximal, i.e. $(T-2)$-dimensional, simplices. Thus, if $\Rad_1(\Cone_v) \leq \beta_3$ then $h^1((F^{\cI}X)_{v},\Lambda) \geq \frac{1}{\binom{T-1}{3}\beta_3}$. 
\end{corollary}

But before showing the bound on the cone radius, we show that the links are 1-connected and hence in particular have vanishing first cohomology group. 

\begin{lemma} \label{lem: 1-connected of colored faces}
    Assume that the size of the underlying field $q \geq 5$. Then the complex $\CC(U_{[d]\setminus J};(U_{[d]\setminus (J\cup I)})_{I \in \cI\setminus \{J\}})$ is connected if $\vert \cI \vert \geq 3$ and it is 1-connected if $\vert \cI \vert \geq 4$. Thus in particular $H^1(F^{\cI}X_{gH_J},\Lambda) =1$ if $\vert \cI\vert \geq 4$.
\end{lemma} 
\begin{proof}
    Let $Y = \CC(U_{[d]\setminus J};(U_{[d]\setminus (J\cup I)})_{I \in \cI\setminus \{J\}})$. Then by \pref{prop: CC connected if subgroups generate G}, $Y$ is connected if and only if $U_{[d]\setminus J}= \langle U_{[d]\setminus (J\cup I)}\mid I \in \cI\setminus \{J\} \rangle$. The later is equivalent to 
    \[\forall i \in [d] \setminus J \ \exists L \in \cI \setminus \{J\}: i \in [d] \setminus (L \cup J) \iff \bigcup_{L \in \cI \setminus \{J\}} [d]\setminus(L \cup J) = [d] \setminus J.\]
    Now, we have
    \[\bigcup_{L \in \cI \setminus \{J\}} [d]\setminus(L \cup J)= \left( \bigcup_{L \in \cI \setminus \{J\}} [d]\setminus L\right) \cap [d] \setminus J= \left([d] \setminus \bigcap_{L \in \cI \setminus \{J\}}L\right) \cap [d] \setminus J.\]
    Hence $Y$ is connected if and only if $\bigcap_{L \in \cI \setminus \{J\}}L = \emptyset$ which, since the color sets in $\cI$ are distinct, holds if and only if $\lvert \cI \rvert \geq 3$.

    For the 1-connectedness, we want to use \pref{prop: CC simply connected if G is amalgam of subgroups}. We know by the work of \cite{tits_ensembles_1986,Abramenko-TwinBuildings} (see also \cite[Proposition 3.8]{GraveDePeraltaV2025}) that $U_{[d]\setminus J}$
    is the free product of the subgroups $U_{\{i,j\}}, \{i,j\} \subseteq [d]\setminus J$ amalgamated along their intersections. Hence, it remains to show that for each $\left\{ i,j \right\} \subseteq [d]\setminus J$ there exists an $L \in \cI \setminus \{J\}$ such that $\left\{ i,j \right\} \subseteq [d]\setminus L$. 
    To show that this holds if $\vert \cI \vert \geq 4$, we fix $\left\{ i,j \right\} \subseteq [d] \setminus J, i \neq j$. 
    Both $i,j$ can be contained in at most one set from $\cI$, call them $i \in L_i, j \in L_j$ respectively. Since $\vert \cI \vert \geq 4$, there will be at least one set in $L \in \cI \setminus \{J,L_i,L_j\}$. Then $i,j \notin L$ and hence $\{i,j\} \subseteq [d] \setminus L$ as required.
\end{proof}

Thus we know that, given $ \vert \cI \vert \geq 4$, every link $F^{\cI}X_{gH_J}$ is 1-connected and thus admits a cone function, see Remark \pref{rem: connectednes implies cone}. 
\begin{remark}
    The statement above holds in a more general setting. Consider the coset complex $X = \CC(G;(H_i)_{i \in [d]})$, set for $J \subseteq [d], H_J = \cap_{j \in J}H_j$ and assume $\langle H_J,H_L \rangle = H_{J\cap L}$ for all $J,L \subseteq [d]$ (this for example holds if $G,H_i, i\in [d]$ form a subgroup geometric system in the sense of \cite{kaufman_high_2023}). Let $\cI \subseteq \binom{[d]}{r+1}$ be a set of disjoint colors. Assume that the chamber system generated by $G, H_i,i\in [d]$ is $m$-simply connected, in the sense of e.g. \cite[Definition 2.13]{GraveDePeraltaV2025}. Then $F^{\cI}X$ is 1-connected if $\vert \cI \vert \geq m+1$. The general proof is similar to the one stated above. Since we do not need the general result here, we will not discuss the detailed proof here. 
\end{remark}

It remains to bound the radius of the cone function. We will do that by considering the Dehn function of a specific presentation of $U_{[d]\setminus J}$ in the next step. 

\subsubsection{Bounding the cone radius using the Dehn function} \label{subsec:overview-bins}
We bound the radius of the cone function using the Dehn function of a specific presentation using the following.

\begin{theorem} \label{thm: dehn implies cone}
Let $G$ be a finite group with subgroups $H_i$, $i\in I$, where $|I|\geq 3$, and assume that the coset complex $Y=\CC(G;(H_i)_{i\in I})$ is $1$-connected. Let $\Dehn_Y$ be as in \pref{def: dehn function of coset complex}, and write $D=\diam(Y)$. Then $Y$ admits a $1$-cone $\Cone_Y$ satisfying
\[
\Rad_0(\Cone_Y)\leq D,\qquad
\Rad_1(\Cone_Y)\leq 15\bigl(\Dehn_Y(2D+1)+2D+1\bigr).
\]
\end{theorem}
The relation between Dehn functions and fillings of coset complexes was previously studied in \cite[Theorem 6.6]{KaufmanO2021} where this theorem was first discovered. We extend their theorem in a rather straightforward manner in  \pref{app:dehn-equivalence}, where we also prove that the converse direction is also true --- a cone can also serve as a bound to the Dehn function. 

Thus, in order to bound the cone radius, it suffices to bound the diameter and the Dehn function of the relevant coset complex. We will prove the following theorem.

Note that the polynomial $p$ is increasing in both variables. Thus, we have reduced the original problem to showing the following theorem.

\begin{theorem} \label{thm: dehn and diameter bound for link}
    Let $r,d \in \N, d \geq 3$, and $T\in \N, T\geq 4$. There exist $\beta_4, \beta_5>0,\beta_6, \beta_7\geq 0$ and $q_0 \in \N$ such that for every prime power $q \geq q_0$ and every nice set of $T$ disjoint colors $\cI \in \binom{[d]}{r+1}^{[T]}$ for every $J \in \cI$ the complex
    \[Y = \CC(U_{[d]\setminus J};(U_{[d]\setminus (J\cup I)})_{I \in \cI\setminus \{J\}})\]
    has $\diam(Y) \leq \beta_4$ and $\Dehn_Y(m) \leq \beta_5 m^2 + \beta_6m+\beta_7$ for all $m  \in \N$.
\end{theorem}

Thus, it remains to bound the diameter and the Dehn function of $Y$. For the diameter, we get the following bound, which will be proven in \pref{sec: bounding the Dehn function}, see \pref{prop: bounding diameter}.

\begin{lemma} \label{lem: diameter Y}
    It holds that \(\diam(Y) \leq 6(r+1)+1.\)
\end{lemma}

We reduce the task of bounding the Dehn function further below.

\subsection{Step 4: Decomposing into bins} \label{subsec:overview-main-proof}
Let \(\mathcal{I}\in\binom{[d]}{r+1}^{[T]}\) be a nice set of \(T\) disjoint color sets, fix \(J\in\cI\), and let \(B_0,\dots,B_r\) be its bins in the path order of \pref{def:nice-colors}. In the upper triangular matrix model, a nonempty bin \(B=[a,b]\) acts on the coordinate block \(\{a,\dots,b+1\}\), so \(U_B\cong\Unip_{b-a+2}(\F_q)\). The blocks corresponding to different bins are disjoint, and every simple root indexed by \([d]\setminus J\) belongs to one of the bins. Therefore the group decomposes as the internal direct product
\[U_{[d]\setminus J}=\prod_{\ell=0}^r U_{B_\ell}.\]
Let \(Y_\ell=\CC(U_{B_\ell};(U_{B_\ell\setminus I})_{I\in\cI\setminus\{J\}})\).

The following lemma is a specific case of \cite[Lemma 16]{Abramenko-TwinBuildings}.
\begin{lemma} \label{lem: oppo decomposes into join}
    \[\CC(U_{[d]\setminus J};(U_{[d]\setminus (J \cup \{i\})})_{i \in [d]\setminus J}) = *_{\ell=0}^r \CC(U_{B_\ell};(U_{B_\ell \setminus \{i\}})_{i \in B_\ell}) \]
\end{lemma}

With this in mind, we investigate the colored faces complex of a join.

\begin{lemma} \label{lem: colored faces of join}
Let $X$ be a finite, pure $d$-dimensional, $(d+1)$-partite simplicial complex. Assume that $X=Y_1*Y_2$ is the join of two complexes $Y_1,Y_2$, and let $\cI$ be a collection of pairwise disjoint nonempty color sets. Write $C_i=\col(Y_i)$ and $\cI_i=(J\cap C_i)_{J\in\cI}$ for $i=1,2$, retaining the original label $J$ on each restricted color set, including empty ones.

Define the \emph{colored join} $F^{\cI_1}Y_1 \overset{*}{\times} F^{\cI_2}Y_2$ as follows. Its vertices are pairs $(\sigma_1,\sigma_2)$ such that $\sigma_i\in Y_i$ and, for some $J\in\cI$, $\col(\sigma_i)=J\cap C_i$ for $i=1,2$. We assign this vertex color $J$. Empty faces are allowed in either coordinate. A collection $\{(\sigma_{1,j},\sigma_{2,j})\}_{j=0}^k$ forms a face if its vertices have distinct colors and $\bigcup_{j=0}^k\sigma_{i,j}\in Y_i$ for each $i=1,2$.

Then the map $(\sigma_1,\sigma_2)\mapsto \sigma_1\sqcup\sigma_2$ induces an isomorphism
\[F^{\cI_1}Y_1 \overset{*}{\times} F^{\cI_2}Y_2 \cong F^{\cI}X.\]
\end{lemma}

\begin{proof}
Every face $\sigma\in X=Y_1*Y_2$ decomposes uniquely as $\sigma=\sigma_1\sqcup\sigma_2$ with $\sigma_i\in Y_i$. If $\col(\sigma)=J$, then $\col(\sigma_i)=J\cap C_i$. Thus the proposed map is a color-preserving bijection between the vertex sets.

Let $\{(\sigma_{1,j},\sigma_{2,j})\}_{j=0}^k$ be a collection of vertices, and write $\sigma_j=\sigma_{1,j}\sqcup\sigma_{2,j}$. Their images form a face of $F^{\cI}X$ precisely when their colors are distinct and $\bigcup_{j=0}^k\sigma_j\in X$. By the definition of the join,
\[\bigcup_{j=0}^k\sigma_j\in X \iff \bigcup_{j=0}^k\sigma_{i,j}\in Y_i \text{ for each } i=1,2.\]
These are exactly the conditions defining a face of the `colored join', proving the isomorphism.
\end{proof}

We can iterate the result and get the following corollary.
\begin{corollary} \label{cor:faces-of-join}
    Let $X= Y_1* \dots * Y_k$ be a finite, pure $d$-dimensional, $(d+1)$-partite simplicial complex. Further, let $\cI$ be a set of distinct colors. Then 
    \[F^{\cI}X = \prod_{1\leq i \leq k}^* F^{\cI_i}Y_i.\]
    Here $\displaystyle \prod^*$ denotes the straightforward generalization of $\overset{*}{\times}$ to multiple complexes.
\end{corollary}

Let \(X = \CC(U_{[d]\setminus J};(U_{[d]\setminus (J \cup \{i\})})_{i \in [d]\setminus J}), Y = F^{\cI}X \). By the above results, we have
 \[Y =  \prod_{0\leq \ell \leq r}^*\CC(U_{B_\ell};(U_{B_\ell \setminus I})_{I \in \cI \setminus \set{J}}).\]
Denote by \(Y_\ell = \CC(U_{B_\ell};(U_{B_\ell \setminus I})_{I \in \cI \setminus \set{J}})\) the \(\ell\)-th component of this product.
\begin{lemma}\label{lem:direct-product-of-groups in overview}
    Assume that every \(Y_\ell\) has a bound of \(Dehn_{Y_\ell}(t)\leq At\) on its Dehn function. Then \(Y\) has a bound of \(Dehn_Y(t) \leq O((t(r+1))^2)+(r+1)At\) on its Dehn function.
\end{lemma}
The idea is that we take a word written with generators from \(K_{I_0 \cup I_m}\), and use \(\poly(r)\) commuting relations in order to write the word as \(r+1\) words that live only in the respective \(K^{B_{\ell}}\), and then to contract each word separately. We prove this lemma in \pref{sec:dirprod}.

Thus, as a last step, we need to bound the Dehn function of each component $Y_\ell$.

\begin{lemma} \label{lem: Yl Dehn function}
    Let $Y_\ell$ be as above. Then there exists a universal constant \(A\) such that $\Dehn_{Y_\ell}(t) \leq At$.
\end{lemma}
We prove this lemma in \pref{sec: bounding the Dehn function}, see \pref{lem:dehn-bound-on-bin}.

\subsection{Proof of the main theorem}
We now put all the pieces together to prove \pref{thm: main thm ch 3}.
We fix the parameters $r,d,\gamma_1,n$ as in the statement of the theorem and set $T=r^{10}$. Let $q=p^m$ be a prime power and $\Lambda$ a group without non-trivial elements of order $p$. Fix an irreducible polynomial $f\in \F_q[t]$ of degree at least 2. Let $X = X(A_d,\F_q,f)= \CC(G,(H_i^f)_{i\in [d]})$ be the corresponding KMS complex. Let $\cI \in \binom{[d]}{r+1}^{[T]}$ be a nice set of disjoint colors. Then 
\(F^{\cI}X \cong \CC(G,(H_J^f)_{J\in \cI})\) (where $H_J^f = \cap_{j \in J} H_j^f$) and for any vertex $gH_J \in F^{\cI}X$ its link has the form 
\[(F^{\cI}X)_{gH_J} \cong \CC(U_{[d]\setminus J};(U_{[d]\setminus (J\cup I)})_{I \in \cI\setminus \{J\}}). \]
We fix one such link and denote it by \(Y\). We use the path order associated with \(J\) in \pref{def:nice-colors}, and let \(B_0,\dots,B_r\) be the resulting \(J\)-bins. Then \(Y\) decomposes as restricted join
\[Y =  \prod_{0\leq \ell \leq r}^*\CC(U_{B_\ell};(U_{B_\ell \setminus I})_{I \in \cI \setminus \set{J}}),\]
and we denote by \(Y_\ell = \CC(U_{B_\ell};(U_{B_\ell \setminus I})_{I \in \cI \setminus \set{J}})\) the \(\ell\)-th component of this decomposition. By \pref{lem: Yl Dehn function}, we know that there exists a universal constant such that $\Dehn_{Y_\ell}(t) \leq At$. From \pref{lem:direct-product-of-groups in overview} we deduce that \(\Dehn_Y(t) \leq \beta_5t^2 + \beta_6 t + \beta_7= O((t(r+1))^2)\) with $\beta_5,\beta_6,\beta_7$ at most quadratic in $r$ and independent of all other parameters. By \pref{lem: diameter Y}, we know that $\diam(Y) \leq 6(r+1)+1=\beta_4$. Thus we have established \pref{thm: dehn and diameter bound for link}. 

Next, we want to apply \pref{thm: dehn implies cone} to $Y$. By \pref{lem: 1-connected of colored faces}, since $\vert \cI\vert = T >4$ and assuming $q \geq 5$, we know that $Y$ is simply connected. Thus we get 
\[\Rad_0(\Cone_Y) \leq \beta_4, \qquad \Rad_1(\Cone_Y) \leq 15(\beta_5(2\beta_4+1)^2 + \beta_6 (2\beta_4+1) + \beta_7 +2\beta_4+1) =:\beta_3.\]
Since $\beta_4$ is of order $r$, $\beta_5,\beta_6,\beta_7$ are of order $r^2$, we get that $\beta_3$ is of order $r^4$ and independent of all other parameters. Thus we have established \pref{prop: restatement bound on 1-cone} with $q_0=5$. Using \pref{cor: bound on cone radius gives cbe our case}, we get 
\[h^1((F^{\cI}X)_{gH_J},\Lambda) \geq \frac{1}{\binom{T-1}{3}\beta_3}=:\beta_2.\]
Since $T = r^{10}$, $\binom{T-1}{3}$ is of order $T^3=r^{30}$. In total, $\beta_2$ is of order $\frac{1}{r^{34}}$. Together with \pref{lem: 1-connected of colored faces}, we have shown \pref{prop: cbe of links}. 

In the following step, we want to combine coboundary expansion of the links with local spectral expansion to deduce cocycle expansion of $F^{\cI}X$ via the local-to-global theorem. By \pref{lem: colored faces is loc spec} we know that for each $\gamma_2$ there exists $q_0$ such that if $q\geq q_0$ then $F^{\cI}X$ is a $\gamma_2$-local spectral expansion. Given $\beta_2$ we choose $\gamma_2$ small enough such that $\frac{(1-\gamma_2) \beta_2}{24}-e\gamma_2=: \beta_1$ is positive. This choice of $\gamma_2$ then determines a lower bound for $q_0$. After fixing $q_0$, $\beta_1$ will grow with the same order as $\beta_2$, i.e. order $r^{-34}$. Together with \pref{lem: vanishing of cohom for colored faces complex} we have shown \pref{prop: nice colored faces are CBE}. 
Finally, we will apply \pref{thm: DD-reduction to colored faces} to deduce coboundary expansion of $F^rX$. By our choice of $n$ and $T$, we have $F^{\cI}X = F^{\cI}(X^{\leq n})$, in particular $F^{\cI}(X^{\leq n})$ is a $\beta_1$-coboundary expander for any nice set of colors $\cI$. Our choice of $n$ ensures that we can apply \pref{lem: LSE for faces complex of skeleton} which gives a lower bound on $q_0$ and shows that $F^r(X^{\leq n})$ is a $\gamma_1$-local spectral expander. Since
\[\Prob[\cI \in \binom{[d]}{r+1}^{[T]}]{\cI \text{ is nice}} \geq 0.9\]
by \pref{lem:most-colors-are-nice}, we have 
\[\Prob[\cI \in \binom{[d]}{r+1}^{[T]}]{h^1(F^I(X^{\leq n}),\Lambda) \geq \beta_1 \land H^1(F^{\cI}(X^{\leq n}),\Lambda)=1} \geq 0.9\]
and thus get $H^1(F^rX,\Lambda)=1$ and $h^1(F^rX,\Lambda) \geq \Omega(r^{-34})$ which finishes the proof.

\ivnote{Check what you think about the new version (old one in comments).}
\section{The color sampling lemma} \label{sec: color sampling lemma}
In this subsection we prove \pref{lem:most-colors-are-nice}. We start with the following observation. Cutting the cyclic order at one of the \(r+1\) colors of a color set leaves a path of \(d\) positions with \(r\) selected colors. These divide the remaining \(d-r\) coordinates into \(r+1\) bins, so the average bin size is \(\frac{d-r}{r+1}\approx\frac{d+1}{r+1}\). We first show that, with high probability, no bin is larger than this average by more than a logarithmic factor (and to cover the first item, also that no bin is smaller than \(101\)). More precisely, for most collections \(\cI\), all bins of all \(I\in\cI\) have size \(\lesssim\frac{\ln T}{r+1}(d+1)\). This is the content of \eqref{eq:bin-is-small}, up to a union bound over the color sets.

Next, fix a color set defining a small bin \(B\), and a forbidden subset \(S\subseteq B\) occupying at most four fifths of the bin, up to rounding. We show that a random color set \(I'\) avoids \(S\) with probability \(\gtrsim(Tr)^{-\alpha}\), where \(\alpha=\frac{4(r+1)}{5r}<0.81\). Since we have \(T-1\) other color sets to try, the probability that all of them intersect \(S\) is exponentially small in \(T(Tr)^{-\alpha}\). For \(T\geq r^{10}\), this is small enough to take a union bound over all the requirements in the definition of niceness.\footnote{The meticulous reader may object and say that the color sets in \(\cI\) must be disjoint, so they are not independent: why should having many candidates make it so unlikely that they all intersect \(S\)? The point is that, under our lower bound on \(d\), independently sampled color sets are already disjoint with high probability. Thus requiring disjointness changes the failure probability only slightly.} The same argument applies when \(S\) is contained in two bins of this same color set, provided its total size satisfies the same bound; fixing that color set determines both bins before we sample the candidates.

Finally, we apply this observation to the forbidden sets appearing in the definition of niceness. For example, to find a color set \(I'\) whose intersection with \(B\) lies in its first fifth, we take \(S\) to be the remaining four fifths of \(B\). Avoiding \(S\) is exactly the required condition. The last-fifth and middle-fifth conditions are handled similarly. For the fourth condition, \(S\) is the union of two intervals, each occupying approximately two fifths of its respective bin. Since both bins are small, the total size of \(S\) is again at most four fifths of our bin-size cutoff, up to rounding. Thus the observation applies to every requirement, and a union bound shows that they all hold simultaneously with high probability.

\begin{proof}[Proof of \pref{lem:most-colors-are-nice}]
Before proving the lemma, we will prove that all sizes of all bins are small with high probability. Then we will show for each bin, and for every pair of bins, that there exists a color set that avoids the respective intervals appearing in the definition of niceness. This will make use of the fact that the intervals we need to avoid occupy only a constant fraction of the relevant bins. Note that different requirements in the definition may have different color sets witnessing them.

We first sample \(I_0,\dots,I_{T-1}\) independently and uniformly from \(\binom{[d]}{r+1}\). These sets need not be disjoint. Let \(D\) be the event that they are pairwise disjoint. We will estimate the probability that the niceness requirements hold in this independent sampling, and only at the end condition on \(D\).

Indeed, the probability that two independently sampled color sets intersect is at most \((r+1)^2/(d+1)\). To see this, fix the first set. Each of its \(r+1\) colors belongs to the second set with probability \((r+1)/(d+1)\), and we take a union bound over these colors. Taking another union bound over all pairs of color sets gives \[\Prob{D^c}\leq \binom{T}{2}\frac{(r+1)^2}{d+1}\leq \frac{1}{200}.\] Conditioned on \(D\), all ordered collections of pairwise disjoint color sets are equally likely. Forgetting the ordering gives precisely the distribution in the lemma, since each collection has the same number \(T!\) of orderings.

Set \(A\coloneqq100T(r+1)\), \(m\coloneqq d(1-A^{-1/r})\), and \(\alpha\coloneqq\frac{4(r+1)}{5r}\). We claim that
\begin{equation} \label{eq:bin-is-small}
    \Prob[I_0]{\exists\text{ a bin }B\text{ of }I_0\text{ with }|B|>m}\leq \frac{1}{100T}.
\end{equation}
This probability is small enough so that we could later on union bound on all \(T\) sets \(I_j\).

\ydnote{changed explanation to accomodate for cyclicity:}To prove the claim, choose \(j_*\) uniformly from \(I_0\) and condition on its value. Conditional on \(j_*\), the other \(r\) colors of \(I_0\) form a uniformly random \(r\)-subset of \([d]\setminus\{j_*\}\). Rotate the colors so that \(j_*=0\). The remaining selected colors are then a uniformly random \(r\)-subset of the path \(\{1,\dots,d\}\), whose \(r+1\) gaps are exactly the bins of \pref{def:nice-colors}. Changing \(j_*\) only permutes these bins, so the event in \eqref{eq:bin-is-small} is independent of the choice of \(j_*\). There is no additional union bound over possible choices of \(j_*\).

Fix one of the \(r+1\) bin positions, and let \(B\) be its bin. For an integer \(0\leq s\leq d-r\), under this conditional distribution we have \[\Prob{|B|\geq s}=\frac{\binom{d-s}{r}}{\binom dr}.\] The denominator counts all possible choices of the \(r\) selected positions. For the numerator, take a choice whose specified bin has at least \(s\) unused positions, delete the first \(s\) of these positions, and shift all subsequent positions down by \(s\). This produces an \(r\)-subset of \(\{1,\dots,d-s\}\). Conversely, inserting \(s\) unused positions at the start of the specified bin recovers a unique choice with a bin of size at least \(s\). This bijection proves the equality.

Canceling factorials gives \[\frac{\binom{d-s}{r}}{\binom dr}=\prod_{j=0}^{r-1}\left(1-\frac{s}{d-j}\right)\leq\left(1-\frac{s}{d}\right)^r.\] Taking \(s=\lfloor m\rfloor+1\), we obtain \[\Prob{|B|>m}\leq\left(1-\frac{m}{d}\right)^r=\frac1A.\] If \(s>d-r\), the probability is zero, since there are only \(d-r\) unselected positions. A union bound over the \(r+1\) bins gives \((r+1)/A=1/(100T)\). This holds for every value of \(j_*\), so averaging over \(j_*\) proves \eqref{eq:bin-is-small}. A further union bound over the \(T\) color sets shows that the probability that any bin has size greater than \(m\) is at most \(1/100\).

We also need to rule out bins of size at most \(100\). Let \(E\) be the event that such a bin occurs in one of the sampled color sets. The two selected colors bounding such a bin have clockwise distance at most \(101\). View each color set as sampled in a uniformly random order. For two fixed, distinct sampling positions, conditioned on the first color, the second is uniform among the \(d\) remaining colors, so it lies in the next \(101\) positions with probability \(101/d\). A union bound over the \(r(r+1)\) ordered pairs in each of the \(T\) color sets gives \[\Prob{E}\leq\frac{101Tr(r+1)}{d}\leq\frac2T<10^{-5},\] where we used \(d+1\geq100(r+1)^2T^2\), \(T\geq r^{10}\), and \(r\geq10^4\).

The next observation estimates the probability of finding a color set that avoids a prescribed set \(S\). Later, \(S\) will consist of the colors that are forbidden to a witness for one particular niceness requirement. The exponent \(\alpha\) accounts for the fact that the bin cutoff came from \(r\) random cut points, whereas a candidate color set contains \(r+1\) colors.
\begin{observation}~
\begin{enumerate}
    \item Let \(S\subseteq[d]\) be a fixed set of size at most \(\frac45m+20\). For a uniformly sampled color set \(I_1\), \[\Prob[I_1]{I_1\cap S=\emptyset}\geq\frac{1}{2A^\alpha}.\]
    \item If \(t\) color sets are sampled independently of one another and of \(S\), then the probability that none of them avoids \(S\) is at most \[\left(1-\frac{1}{2A^\alpha}\right)^t\leq\exp\left(-\frac{t}{2A^\alpha}\right).\]
\end{enumerate}
\end{observation}

To prove the first item, we sample the \(r+1\) colors of \(I_1\) one at a time without repetition. Suppose that the first \(j\) sampled colors were all outside \(S\). There are then \(d+1-j\) colors left to choose from, of which \(d+1-|S|-j\) are outside \(S\). Therefore, conditioned on having avoided \(S\) so far, the probability that the next color also avoids \(S\) is \((d+1-|S|-j)/(d+1-j)\). Multiplying these conditional probabilities for \(j=0,\dots,r\) gives \[\Prob[I_1]{I_1\cap S=\emptyset}=\prod_{j=0}^{r}\frac{d+1-|S|-j}{d+1-j}.\] This uses conditional probabilities; the individual sampled colors are not independent.

Put \(\theta=\frac15+\frac45A^{-1/r}\). Since \(m=d(1-A^{-1/r})\leq(d+1)(1-A^{-1/r})\), the bound on \(|S|\) gives \[d+1-|S|\geq d+1-\frac45m-20\geq(d+1)\theta-20.\] Since \(j\leq r\), each numerator in the product above is at least \((d+1)\theta-(r+20)\), while its denominator is at most \(d+1\). This lower bound on the numerator is positive under our assumptions. Hence \[\Prob[I_1]{I_1\cap S=\emptyset}\geq\left(\theta-\frac{r+20}{d+1}\right)^{r+1}\geq\theta^{r+1}\left(1-\frac{5(r+20)}{d+1}\right)^{r+1},\] where the last inequality uses \(\theta\geq1/5\).

By the elementary inequality \((1-x)^{r+1}\geq1-(r+1)x\) for \(0\leq x\leq1\), the last factor is at least \[1-\frac{5(r+1)(r+20)}{d+1}\geq\frac12.\] Here we used \(d+1\geq100(r+1)^2T^2\geq10(r+1)(r+20)\). Finally, the weighted arithmetic--geometric mean inequality gives \[\theta=\frac15\cdot1+\frac45\cdot A^{-1/r}\geq\left(A^{-1/r}\right)^{4/5},\] and therefore \(\theta^{r+1}\geq A^{-4(r+1)/(5r)}=A^{-\alpha}\). Combining these inequalities proves the first item of the observation. For the second item, each candidate fails to avoid \(S\) with probability at most \(1-1/(2A^\alpha)\). Independence allows us to multiply these failure probabilities, and the last inequality follows from \(1-x\leq e^{-x}\).

We now explain how to apply the observation to each of the four items in the definition of niceness. For each item, we must identify a forbidden set \(S\) such that \(I'\cap S=\emptyset\) guarantees the required property of \(I'\). We must also check that \(|S|\leq\frac45m+20\) when the relevant bins have size at most \(m\). Only colors in \([d]\) can belong to \(I'\), so any forbidden set below is intersected with \([d]\) if necessary; this can only decrease its size. For each fixed color set \(I\), we use its path order from \pref{def:nice-colors} for all color sets and intervals below.
\begin{enumerate}
    \item[item 1] Fix a bin \(B=[a,b]\) of \(I\). We want a color set \(I'\ne I\) whose colors in \(B\) all lie in the first fifth of the bin. Thus we take \(S=(a+\lfloor\frac{b-a}{5}\rfloor,b]\). Avoiding \(S\) means precisely that \(I'\cap B\subseteq[a,a+\lfloor\frac{b-a}{5}\rfloor]\). The number of forbidden colors is \[|S|=(b-a)-\left\lfloor\frac{b-a}{5}\right\rfloor\leq\frac45(b-a+1)=\frac45|B|.\] Here we used \(\lfloor\frac{b-a}{5}\rfloor\geq\frac{b-a-4}{5}\). In particular, \(|S|\leq\frac45m\) if \(|B|\leq m\).

    \item[item 2] We want the colors of \(I'\) in \(B\) to lie in its last fifth. Accordingly, we take \(S=[a,b-\lfloor\frac{b-a}{5}\rfloor)\). Avoiding \(S\) means precisely that \(I'\cap B\subseteq[b-\lfloor\frac{b-a}{5}\rfloor,b]\). This forbidden set has the same size as in item 1, so again \(|S|\leq\frac45m\) whenever \(|B|\leq m\).

    \item[item 3] We want the colors of \(I'\) in \(B\) to lie in the specified middle interval. We therefore take \(S=[a,b]\setminus[a+2\lceil\frac{b-a}{5}\rceil,a+3\lfloor\frac{b-a}{5}\rfloor]\). Avoiding \(S\) is exactly the required condition. The middle interval contains \[\max\left\{0,3\left\lfloor\frac{b-a}{5}\right\rfloor-2\left\lceil\frac{b-a}{5}\right\rceil+1\right\}\] colors, where the maximum with zero accounts for the possibility that the interval is empty. Since \(\lceil\frac{b-a}{5}\rceil\leq\lfloor\frac{b-a}{5}\rfloor+1\), this number is at least \[\left\lfloor\frac{b-a}{5}\right\rfloor-1\geq\frac{|B|}{5}-2.\] Hence the number of colors outside the middle interval is at most \(\frac45|B|+2\), and \(|S|\leq\frac45m+2\) when \(|B|\leq m\).

    \item[item 4] Fix two bins \(B=[a,b]\), \(B'=[a',b']\) of the same color set \(I\) and a pair \((\ell,\ell')\in\set{0,1,2,3}^2\). This time we need the same color set \(I'\) to avoid both intervals in the fourth item. Thus the forbidden set is their union: \[\begin{aligned}
        S={}&[a+\ell\lfloor\tfrac{b-a}{5}\rfloor,\min\{a+(\ell+2)\lceil\tfrac{b-a}{5}\rceil+4,b\}]\\
        &\quad\cup[a'+\ell'\lfloor\tfrac{b'-a'}{5}\rfloor,\min\{a'+(\ell'+2)\lceil\tfrac{b'-a'}{5}\rceil+4,b'\}].
    \end{aligned}\] Avoiding this union means avoiding both intervals simultaneously, as required. The first interval has at most \[(\ell+2)\left\lceil\frac{b-a}{5}\right\rceil-\ell\left\lfloor\frac{b-a}{5}\right\rfloor+5\leq2\left\lfloor\frac{b-a}{5}\right\rfloor+\ell+7\leq\frac25|B|+10\] colors. The \(+5\) in the first expression accounts for the \(+4\) in the upper endpoint and for counting both endpoints. The first inequality uses \(\lceil\frac{b-a}{5}\rceil\leq\lfloor\frac{b-a}{5}\rfloor+1\), and the second uses \(\ell\leq3\). The same argument bounds the second interval by \(\frac25|B'|+10\). Consequently, \[|S|\leq\frac25|B|+\frac25|B'|+20\leq\frac45m+20\] whenever both bins have size at most \(m\).
\end{enumerate}

We have now checked that every individual niceness requirement can be satisfied by finding a color set that avoids a set \(S\) of the size allowed in the observation. We next bound the probability that such a color set does not exist. First fix which requirement we are considering: the excluded color set \(I\), the position of its bin in items 1--3, or the two bin positions and \((\ell,\ell')\) in item 4. At this point the bins, and hence \(S\), are still random.

Reveal the one or two color sets defining these bins. For example, if the two bins in item 4 belong to \(I_j\) and \(I_{j'}\), we reveal these sets. This determines the bins and \(S\). The other color sets are still independently and uniformly sampled. We look for a witness among them, also excluding \(I\), since the definition requires \(I'\ne I\). In total we exclude at most three color sets, so at least \(T-3\) independent candidates remain. Whenever the revealed bins have size at most \(m\), the observation shows that the probability that all these candidates fail is at most \[\exp\left(-\frac{T-3}{2A^{4/5}}\right).\]

This bound holds for every choice of the revealed color sets with the required bin sizes. Averaging over those choices, the probability that the defining bins have size at most \(m\) and this particular requirement has no witness is also at most the displayed bound. Note that we have not conditioned on all bins of all color sets being small. We fixed only the sets defining the requirement, so the candidate witnesses retained their original independent distributions.

It remains to make all requirements hold simultaneously. For items 1--3 there are \(T\) choices of \(I\), \(r+2\) bin positions, and three items to check. Thus these items give at most \(3T(r+2)\) requirements. For item 4 there are at most \(T\) choices of the excluded set \(I\), at most \(T(r+2)\) choices for each of the two bin positions, and \(16\) choices of \((\ell,\ell')\). This gives at most \(16T^3(r+2)^2\) requirements. This count allows the two bins to belong to arbitrary color sets, as in the definition.

Let \(F\) be the event that some niceness requirement fails in the independent sampling. Either there is a bin larger than \(m\), or there is a failed requirement whose defining bins all have size at most \(m\). The first event has probability at most \(1/100\), and we have just bounded the probability of the second event for each individual requirement. Taking a union bound gives \[\Prob{F}\leq\frac1{100}+\bigl(3T(r+2)+16T^3(r+2)^2\bigr)\exp\left(-\frac{T-3}{2A^{4/5}}\right).\]

We finish by checking explicitly that this error bound is small for the stated values of \(r\) and \(T\). The number of requirements satisfies \[3T(r+1)+16T(r+1)^2\leq20T(r+1)^2\leq21Tr^2.\] Also, \(T-1\geq T/2\), \(r+1\leq1.01r\), and \(\alpha\leq0.81\), so \[\frac{T-1}{2A^\alpha}\geq\frac{T^{1-\alpha}}{4[100(r+1)]^\alpha}\geq\frac{T^{1-\alpha}}{200r^\alpha},\] where we used \(4\cdot101^{0.81}<200\). Consequently, \[\Prob{F}\leq\frac1{100}+21Tr^2\exp\left(-\frac{T^{1-\alpha}}{200r^\alpha}\right).\]

To bound the second term, first hold \(r\) fixed and increase \(T\). Since \(1-\alpha\geq0.19\) and \[\frac{T^{1-\alpha}}{r^\alpha}\geq r^{10-11\alpha}=r^{6/5-44/(5r)}\geq r^{11/10},\] the derivative of the logarithm of this term with respect to \(T\) is \[\frac1T\left(1-\frac{(1-\alpha)T^{1-\alpha}}{200r^\alpha}\right)<0\qquad\text{for }T\geq r^{10},\ r\geq10^4.\] Thus the term is largest at \(T=r^{10}\), where it is at most \(21r^{12}\exp(-r^{11/10}/200)\). The logarithmic derivative of this upper bound with respect to \(r\) is \[\frac1r\left(12-\frac{11r^{11/10}}{2000}\right)<0\qquad\text{for }r\geq10^4.\] It is therefore enough to evaluate at \(r=10^4\). At this value its logarithm is \[\ln21+48\ln10-\frac{10^{22/5}}{200}<-12,\] so the term is less than \(10^{-5}\). We conclude that \(\Prob{F}\leq1/100+10^{-5}\).

Finally, we return to the distribution in the lemma by conditioning on pairwise disjointness. We proved at the start that \(\Prob{D}\geq1-1/200\). Thus \[\Prob[D]{F\cup E}\leq\frac{\Prob{F}+\Prob{E}}{\Prob{D}}\leq\frac{1/100+10^{-5}+2/T}{1-1/200}<0.011<0.1.\] On \(D\), the collection is nice whenever neither \(F\) nor \(E\) occurs. Hence a uniformly sampled collection of \(T\) pairwise disjoint color sets is nice with probability at least \(0.9\), as claimed.
\end{proof}

\section[Bound on the Dehn function for the group of upper triangular matrices]{Bound on the Dehn function for \(\Unip_n (\mathcal{R})\)}
\label{sec: Unip sec}

In \pref{lem: Yl Dehn function} in the proof outline above, we bound the Dehn function of coset complex $Y_{\ell}$ based on the relations of the vertex stabilizers under the niceness assumption. We note that $Y_{\ell}$ is a coset complex of an upper-triangular matrix group and bounding the Dehn function requires studying the Dehn function of such groups under group presentations related to the niceness assumption. Roughly speaking, the niceness assumption is meant to ensure that the vertex stabilizers in $Y_{\ell}$ will have ``enough'' relations to yield a good bound on the Dehn function. In this section, we study presentations of upper-triangular matrix groups and derive a presentation that on one hand has enough relations to ensure a linear bound on the Dehn function and on the other hand its set of relations in ``small'' enough to be derived from the niceness assumption. 

Throughout let $\mathcal{R}$ be a commutative unital ring in which we denote the unit as $1$. For our application below, $\mathcal{R}$ will be a finite field, but the results of the section are general and do not depend on the field or finiteness properties. 

We denote $\Unip_n (\mathcal{R})$ to be the group of $n \times n$ upper-triangular matrices with entries in $\mathcal{R}$ and $1$'s along the main diagonal with the binary operation of usual matrix multiplication. Below, we will abbriviate and refer to such matrices as upper-triangular (without explicitly mentioning the $1$'s on the main diagonal or the $\mathcal{R}$ entries). 

The key result of this section is to find a bound for the Dehn function of $\Unip_n (\mathcal{R})$ for specific presentation that will be useful in the sequel.

\subsection{Generalities regarding commutators in groups}

Here we cover some general definitions and facts regarding commutators in groups that will be used below.

Let $G$ be a group. For two elements $g,h \in G$, the \textit{commutator of g and h} is the group element defined as $[g,h] = g^{-1} h^{-1} g h$. We say that $g$ and $h$ \textit{commute} if $[g,h] = e$ or equivalently $gh = hg$. The following identities are well-known (and straightforward to prove):
\begin{claim}
\label{claim:comm ident claim}
Let $G$ be a group and $g,h,g_1,g_2 \in G$ be elements. Then
\begin{itemize}
    \item $[g,h]^{-1} = [h,g]$
    \item $[g,e] = [e,g] =e$
    \item If $[g,h] = e$, then  $[g^{-1},h] = [g^{-1},h^{-1}] = [g,h^{-1}]=e$
    \item If $[g_1,h] = [g_2,h] = e$, then $[[g_1,g_2],h] =e$.
\end{itemize}
\end{claim}

Using the notation $g^h = h^{-1} g h$, the following identity is known as the Hall–Witt identity: For every $g_1, g_2, g_3 \in G$ it always hold that
$$[[g_1, g_2^{-1}],g_3]^{g_2} [[g_2, g_3^{-1}],g_1]^{g_3} [[g_3, g_1^{-1}],g_2]^{g_1} =e.$$

In the sequel, we will need the following implication of this identity:
\begin{proposition}
\label{HW implication prop} 
Let $G$ be a group and $g_1,g_2,g_3$ be elements of $G$. If $[g_3,g_1^{-1}] =e$ and $[[g_1, g_2^{-1}],g_3]$ commutes with $g_2$ and $g_3^{-1}$, then
$$[[g_1, g_2^{-1}],g_3] = [g_1,[g_2,g_3^{-1}]].$$
\end{proposition}

\begin{proof}
We note that the assumption that $[g_3,g_1^{-1}] =e$ implies that $[[g_3, g_1^{-1}],g_2]^{g_1} =e$. Also, by the assumption that $[[g_1, g_2^{-1}],g_3]$ commutes with $g_2$ it follows that $[[g_1, g_2^{-1}],g_3]^{g_2} = [[g_1, g_2^{-1}],g_3]$. Thus, the Hall–Witt identity is reduced to
$$[[g_1, g_2^{-1}],g_3] [[g_2, g_3^{-1}],g_1]^{g_3} =e.$$
Using the assumption that $[[g_1, g_2^{-1}],g_3]$ commutes with $g_3^{-1}$ yields that 
$$e=[[g_1, g_2^{-1}],g_3] [[g_2, g_3^{-1}],g_1]^{g_3} = [[g_1, g_2^{-1}],g_3] g_3^{-1} [[g_2, g_3^{-1}],g_1] g_3 = g_3^{-1} [[g_1, g_2^{-1}],g_3] [[g_2, g_3^{-1}],g_1] g_3,$$
and thus 
$$[[g_1, g_2^{-1}],g_3] [[g_2, g_3^{-1}],g_1] = e$$
or equivalently
$$[[g_1, g_2^{-1}],g_3]  = [[g_2, g_3^{-1}],g_1]^{-1} = [g_1,[g_2, g_3^{-1}]],$$
as needed.
\end{proof}

We will also need the following result
\begin{proposition}[{\cite[Lemma 1]{BissD}}]
\label{BissD result prop}
For a group $G$ and elements $g_1, g_2, h \in G$. Assume that $g_2$ commutes with $[g_1,h]$. Then $[g_1,h] [g_2,h] = [g_1 g_2, h]$.
\end{proposition}

\subsection{Basic presentation and block matrices}

 For $1 \leq i<j \leq n$ and $a \in \mathcal{R}$, we define the elementary matrix $e_{i,j} (a) \in \Unip_n (\mathcal{R})$ to be the matrix with $1$'s along the main diagonal, $a$ in the $(i,j)$ entry and $0$'s in all other entries. The following fact is well-known (see for instance \cite[Theorem 2]{BissD} where a much smaller presentation is shown to suffice):
\begin{theorem}
\label{thm: very basic presentation}
The group $\Unip_n (\mathcal{R})$ can be presented in terms of generators and relations as follows: The set $S_{\elm} = \lbrace e_{i,j} (a) : 1 \leq i<j \leq n, a \in \mathcal{R} \rbrace$ is a generating set of $\Unip_n (\mathcal{R})$ and $\Unip_n (\mathcal{R}) = \langle S_{\elm} \vert R_{\elm} \rangle$ where $R_{\elm}$ are the following relations:
\begin{enumerate}
    \item For every $1 \leq i<j \leq n$ and every $a_1, a_2 \in \mathcal{R}$, $e_{i,j} (a_1) e_{i,j} (a_2) =  e_{i,j} (a_1 + a_2)$.
    \item For every $1 \leq i,j, k,l \leq n$ with $i<j$ and $k<l$ and every $a_1, a_2 \in \mathcal{R}$, 
    $$[e_{i,j} (a_1), e_{k,l} (a_2)] = \begin{cases}
    e_{i,k} (a_1 a_2) & j=k \\
    e_{k,j} (- a_1 a_2) & i =l \\
    e & \text{otherwise}
    \end{cases}.$$
    
\end{enumerate}
\end{theorem}

The basic presentation above is well suited for describing $\Unip_n(\mathcal R)$, but it is less convenient for bounding the Dehn function. The goal of this section is to introduce a block presentation that is adapted to a recursive analysis of the Dehn function.

A key idea is to partition matrices into five large blocks and to regard each off-diagonal block as a single generator carrying an entire matrix rather than a single ring element. This yields a presentation whose generators are naturally indexed by pairs $(i,j)$ of block positions. Explicitly, assuming that $n \geq 5$, we denote $k_0 = \lfloor \frac{n}{5} \rfloor$ and $k_1 = n - 4k_0$ (note that $n \geq 5$ implies that $k_0, k_1$ are positive integers). A matrix $A \in \Unip_n (\mathcal{R})$ can always be written in the following block form:
$$A = 
\begin{pmatrix}
    A_{1,1} & A_{1,2} & A_{1,3} & A_{1,4} & A_{1,5} \\
    0_{k_0} & A_{2,2} & A_{2,3} & A_{2,4} & A_{2,5} \\
    0_{k_1 \times k_0} & 0_{k_1 \times k_0} & A_{3,3} & A_{3,4} & A_{3,5} \\
    0_{k_0} & 0_{k_0} & 0_{k_0 \times k_1} & A_{4,4} & A_{4,5} \\
    0_{k_0} & 0_{k_0} & 0_{k_0 \times k_1} & 0_{k_0} & A_{5,5}
\end{pmatrix}$$
where
\begin{itemize}
    \item $0_{k_i \times k_j}$ denotes $k_i \times k_j$ 0 matrix and $0_{k_i}$ denotes the $k_i \times k_i$ 0 matrix
    \item $A_{1,1}, A_{2,2}, A_{4,4}, A_{5,5} \in \Unip_{k_0} (\mathcal{R}) $
    \item $A_{3,3} \in \Unip_{k_1} (\mathcal{R})$
    \item For $i \neq 3, j \neq 3$, $A_{i,j} \in M_{k_0} (\mathcal{R})$
    \item $A_{1,3}, A_{2,3} \in M_{k_0 \times k_1} (\mathcal{R})$
    \item $A_{3,4}, A_{3,5} \in M_{k_1 \times k_0} (\mathcal{R})$
\end{itemize}
With this block form in mind, we introduce the following notation: 
\begin{enumerate}
    \item For $1 \leq i < j \leq 5$, and a matrix $A$ (of the suitable dimensions), we denote $E_{i,j} (A)$ to be the block matrix in the form above, where $A_{1,1},...,A_{5,5}$ are identity matrices (of the suitable dimensions), $A_{i',j'} = 0$ for every $(i',j') \neq (i,j)$ and $A_{i,j} = A$. E.g., 
    $$E_{1,2} (A) = 
    \begin{pmatrix}
    I & A  & 0_{k_0 \times k_1} & 0_{k_0} & 0_{k_0} \\
    0_{k_0} & I & 0_{k_0 \times k_1} & 0_{k_0} & 0_{k_0} \\
    0_{k_1 \times k_0} & 0_{k_1 \times k_0} & I & 0_{k_1 \times k_0} & 0_{k_1 \times k_0} \\
    0_{k_0} & 0_{k_0} & 0_{k_0 \times k_1} & I & 0_{k_0} \\
    0_{k_0} & 0_{k_0} & 0_{k_0 \times k_1} & 0_{k_0} & I
\end{pmatrix},$$
where $A \in M_{k_0} (\mathcal{R})$.
    \item For $1 \leq i \leq 5$, and an upper triangular matrix $A$ (of the suitable dimensions), we denote $E_{i,i} (A)$ to be the block matrix of the form above, where $A_{i,i} = A$, all the other $A_{i',i'}$'s are identity matrices and the $A_{i,j}$'s are zero matrices. E.g.,
    $$E_{1,1} (A) = 
    \begin{pmatrix}
    A & 0_{k_0}  & 0_{k_0 \times k_1} & 0_{k_0} & 0_{k_0} \\
    0_{k_0} & I & 0_{k_0 \times k_1} & 0_{k_0} & 0_{k_0} \\
    0_{k_1 \times k_0} & 0_{k_1 \times k_0} & I & 0_{k_1 \times k_0} & 0_{k_1 \times k_0} \\
    0_{k_0} & 0_{k_0} & 0_{k_0 \times k_1} & I & 0_{k_0} \\
    0_{k_0} & 0_{k_0} & 0_{k_0 \times k_1} & 0_{k_0} & I
\end{pmatrix},$$
where $A \in \Unip_{k_0} (\mathcal{R})$.
\end{enumerate}

Above, we used the imprecise phrase ``of suitable dimension'', in order to make it precise, we define the following notion: A matrix $A$ will be called \textit{$(i,j)$-admissible}, if $E_{i,j} (A)$ is well-defined. Explicitly, $A$ is $(i,j)$-admissible is equivalent to:
$$\begin{cases}
    A \in M_{k_0} (\mathcal{R}) & i < j, i \neq 3, j \neq 3 \\
    A \in M_{k_0 \times k_1} (\mathcal{R}) & i< 3, j=3 \\
    A \in M_{k_1 \times k_0} (\mathcal{R}) & 3 < j, i=3 \\
    A \in \Unip_{k_0} (\mathcal{R}) & i=j, i \neq 3 \\
    A \in \Unip_{k_1} (\mathcal{R}) & i=j= 3 
 \end{cases}$$

With this notation, we will define a presentation of $\Unip_{n} (\mathcal{R})$ and prove a bound for the Dehn function of this presentation. 

\subsection{Main result - a bound on the Dehn function}
\label{subsection: Main result}

We define the generating set of our presentation to be 
$$S_{\Elm} = \lbrace E_{i,j} (A) : 1 \leq i \leq j \leq 5, (i,j) \neq (1,5), A \text{ is } (i,j) \text{-admissible} \rbrace.$$

We also define the following set of relations $R_{\Elm}$ on the alphabet $S_{\Elm}$:
\begin{enumerate}
    \item[\customlabel{Elm:1}{(Elm1)}]  For every $1 \leq i \leq j \leq 5, (i,j) \neq (1,5)$ and all matrices $A_1, A_2$ that are $(i,j)$-admissible, 
    $$E_{i,j} (A_1) E_{i,j} (A_2) =  
    \begin{cases}
    E_{i,j} (A_1 + A_2) & i <j \\
    E_{i,j} (A_1 A_2) & i = j
    \end{cases}.
    $$
    \item[\customlabel{Elm:2}{(Elm2)}] For every $1 \leq i,j, k,l \leq 5$ such that one of the following holds:
    \begin{itemize}
        \item $1 \leq i < j \leq 4$ and $1 \leq k < l \leq 4$
        \item $2 \leq i < j \leq 5$ and $2 \leq k < l \leq 5$
        \item $\lbrace (i,j), (k,l) \rbrace = \lbrace (1,2), (4,5) \rbrace$
    \end{itemize}
     and all matrices $A_1, A_2$ that are $(i,j)$ and $(k,l)$ admissible correspondingly,
    $$[E_{i,j} (A_1), E_{k,l} (A_2)] = \begin{cases}
    E_{i,k} (A_1 A_2) & j=k \\
    E_{k,j} (- A_2 A_1 ) & i =l \\
    e & \text{otherwise}
    \end{cases}.$$    
     \item[\customlabel{Elm:3}{(Elm3)}] For every $1 \leq i,j, l \leq 5, i \leq j$ such that one of the following holds:
    \begin{itemize}
        \item $1 \leq i \leq j \leq 4$ and $1 \leq l \leq 4$
        \item $2 \leq i \leq j \leq 5$ and $2 \leq l \leq 5$
        \item  $1 \leq i \leq j \leq 2$ and $l= 5$
        \item $4 \leq i \leq j \leq 5$ and $l=1$
    \end{itemize}
     and all matrices $A, B$ that are $(i,j)$ and $(l,l)$ admissible correspondingly,
     $$(E_{l,l} (B))^{-1} E_{i,j} (A) E_{l,l} (B) = 
     \begin{cases}
     E_{i,j} (A) &  i \neq l \text{ and } j \neq l \\
     E_{i,j} (B^{-1} A) & i=l \text{ and } j \neq l \\
     E_{i,j} (AB) & j =l \text{ and } i \neq l \\
     E_{l,l} (B^{-1} AB) & i=j=l
     \end{cases}.$$
\end{enumerate}

We note that currently, we have not yet shown that $S_{\Elm}$ is a generating set of $\Unip_n (\mathcal{R})$ or that $\langle S_{\Elm} \vert \mathcal{R}_{\Elm} \rangle$ is a presentation of $\Unip_n (\mathcal{R})$. Our main result in this section is the following theorem:

\begin{theorem}
\label{thm:unip-presentation}
The group $\Unip_n (\mathcal{R})$ can be presented as $\langle S_{\Elm} \vert R_{\Elm} \rangle$. Moreover, there is a universal constant $C'$ (independent on $n$ and $\mathcal{R}$) such that
the Dehn function with respect to this presentation can be bounded as follows:
\[\Dehn_{ S_{\Elm}, R_{\Elm} } (m) \leq C' m.\]
\end{theorem}

The strategy of proving this theorem goes as follows: we start with a different ``larger'' presentation of $\Unip_n (\mathcal{R})$, which we denote by $\langle S_{\Elm} ' \vert R_{\Elm} ' \rangle$. For this larger presentation, we bound the Dehn function linearly by a universal constant. We then show that moving to our ``smaller'' presentation $\langle S_{\Elm} \vert R_{\Elm} \rangle$ defined above still has almost the same Dehn function - the only thing that changes is the universal constant in the linear bound which increases. 

\subsection{Bounding the Dehn function for a larger presentation}

We define the following set 
\begin{align*}
    & S_{\Elm} ' = \lbrace E_{i,j} (A) : 1 \leq i \leq j \leq 5, A \text{ is } (i,j) \text{-admissible} \rbrace
\end{align*}

By block matrix multiplication, we derive the following relations which we denote by $R_{\Elm} '$:
\begin{enumerate}
    \item[\customlabel{Elm':1}{(Elm1')}]  For every $1 \leq i \leq j \leq 5$ and all matrices $A_1, A_2$ that are $(i,j)$-admissible, 
    $$E_{i,j} (A_1) E_{i,j} (A_2) =  
    \begin{cases}
    E_{i,j} (A_1 + A_2) & i <j \\
    E_{i,j} (A_1 A_2) & i = j
    \end{cases}.
    $$
    \item[\customlabel{Elm':2}{(Elm2')}] For every $1 \leq i,j, k,l \leq 5$ with $i<j$ and $k<l$ and every all matrices $A_1, A_2$ that are $(i,j)$ and $(k,l)$ admissible correspondingly,
    $$[E_{i,j} (A_1), E_{k,l} (A_2)] = \begin{cases}
    E_{i,k} (A_1 A_2) & j=k \\
    E_{k,j} (- A_2 A_1 ) & i =l \\
    e & \text{otherwise}
    \end{cases}.$$    
     \item[\customlabel{Elm':3}{(Elm3')}] For every $1 \leq i,j, l \leq 5, i \leq j$ and all matrices $A, B$ that are $(i,j)$ and $(l,l)$ admissible correspondingly,
     $$(E_{l,l} (B))^{-1} E_{i,j} (A) E_{l,l} (B) = 
     \begin{cases}
     E_{i,j} (A) &  i \neq l \text{ and } j \neq l \\
     E_{i,j} (B^{-1} A) & i=l \text{ and } j \neq l \\
     E_{i,j} (AB) & j =l \text{ and } i \neq l \\
     E_{l,l} (B^{-1} A B & i=j=l
     \end{cases}.$$
\end{enumerate}

\ionote{Added because of Inga's remark}
\begin{proposition}
\label{prop: S' R' pres}
It holds that  $\langle S_{\Elm} ' \vert R_{\Elm} ' \rangle$ is a presentation of $Unip_n (\mathcal{R})$.
\end{proposition}

\begin{proof}
We note that for every $1 \leq i < j \leq n$ and every $a \in \mathcal{R}$,  $e_{i,j} (a) \in S_{\Elm} '$. We also note that all the relations between the  $e_{i,j} (a)$'s specified in \pref{thm: very basic presentation} appear in $R_{\Elm} '$.  Thus By \pref{thm: very basic presentation},  $\langle S_{\Elm} ' \vert R_{\Elm} ' \rangle$ is a presentation of $Unip_n (\mathcal{R})$.
\end{proof}

As we explained above, we will bound the Dehn function of $Unip_n (\mathcal{R})$ with respect to the presentation $\langle S_{\Elm} ' \vert R_{\Elm} ' \rangle$ by a universal linear bound (this bound will be independent of $n$ and $\mathcal{R}$). 

Our method of bounding the Dehn function for the presentation $\langle S_{\Elm} ' \vert R_{\Elm} ' \rangle$ is to decompose the group into a finite nested sequence of groups (where the last group is the complete group $\Unip_n (\mathcal{R})$) and then bound the Dehn function of each of these subgroups recursively (thus achieving a bound for  the Dehn function of  $\Unip_n (\mathcal{R})$).

We define the following order on the set $\lbrace (i,j) : 1 \leq i \leq j \leq 5 \rbrace$: For $(i,j)$ and $(i',j')$ we define $(i,j) < (i',j')$ as follows:
$$(i,j) < (i',j') \Leftrightarrow
 \begin{cases}
 j' - i' > j-i \\
 \text{or}
 \\
j' - i' = j-i \text{ and } i' >i
 \end{cases}.$$
We denote  $(i,j) \leq (i',j')$ if either $(i,j) < (i',j')$ or $(i,j) = (i',j')$. 

With this order, we define the following subgroups of $\Unip_{n} (\mathcal{R})$ for $1 \leq i_0 \leq j_0 \leq 5$:
$$G_{\geq (i_0,j_0)} = \langle E_{i,j} (A_{i,j}) : (i_0,j_0) \leq (i,j), A_{i,j} \text{ is } (i,j) \text{-admissible} \rangle,$$
$$G_{> (i_0,j_0)} = \langle E_{i,j} (A_{i,j}) : (i_0,j_0) < (i,j), A_{i,j} \text{ is } (i,j) \text{-admissible} \rangle,$$
$$G_{(i_0,j_0)} = \langle E_{i_0,j_0} (A_{i_0,j_0}) :  A_{i_0,j_0} \text{ is } (i_0,j_0) \text{-admissible} \rangle.$$

For $1 \leq i \leq j \leq 5$, we further define $S_{\geq (i,j)} = S_{\Elm} ' \cap G_{\geq (i,j)}$ and $R_{\geq (i,j)}$ to be the subset of the relations in $R_{\Elm} '$ that involve only elements in $S_{\geq (i,j)}$. We similarly define $S_{> (i,j)}, S_{(i,j)}, R_{> (i,j)}, R_{(i,j)}$. 

\begin{observation}
\label{obs:trivial-intersect-obs}
We note is that the $(i_0,j_0)$-block in $G_{> (i_0,j_0)}$ is trivial, i.e., it is the $0$ matrix when $i < j$ and the identity matrix if $i=j$. For example,
$$G_{>(2,3)} = 
\left\lbrace
\begin{pmatrix}
    I & 0 & A_{1,3} & A_{1,4} & A_{1,5} \\
    0 & I & 0 & A_{2,4} & A_{2,5} \\
    0 & 0 & I & A_{3,4} & A_{3,5} \\
    0 & 0 & 0 & I & A_{4,5} \\
    0 & 0 & 0 & 0 & I
\end{pmatrix} : \text{the non-zero } A_{i,j} \text{ are } (i,j) \text{-admissible}
\right\rbrace
$$
and the $(2,3)$-block is indeed trivial.
\end{observation}

The following notion will play an important role in the sequel. A group
$G$ is called a semidirect product of subgroups $H$ and $N$, denoted
$G=H\ltimes N$, if every element $g\in G$ admits a unique decomposition
\[
g=hn,
\]
where $h\in H$ and $n\in N$, and multiplication of elements is compatible
with this decomposition. Equivalently, $N$ is normal in $G$, $G=HN$ and
$H\cap N=\{e\}$.

One may think of a semidirect product as a way of representing every
element of $G$ by two pieces of information: an $H$-component and an
$N$-component. The main advantage of such a decomposition is that it
allows one to separate a problem on $G$ into a problem on $H$ and a
problem on $N$. In particular, if
\[
g=hn=e,
\]
where $h\in H$ and $n\in N$, then necessarily $h=e$ and $n=e$.

The next proposition shows that every group $G_{\geq(i_0,j_0)}$
admits such a decomposition, where $G_{(i_0,j_0)}$ records the
$(i_0,j_0)$-block and $G_{>(i_0,j_0)}$ records all higher-order blocks. 

\begin{proposition}
\label{prop:semidirect-prop}
For every $1 \leq i_0 \leq j_0 \leq 5$, it holds that $G_{\geq (i_0,j_0)} = G_{(i_0,j_0)} \ltimes G_{> (i_0,j_0)}$. In particular, for every $g \in G_{\geq (i_0,j_0)}$ there are unique $A_{i_0,j_0}$ and $g' \in G_{> (i_0,j_0)}$ such that $g = E_{i_0, j_0} (A_{i_0,j_0}) g'$.
\end{proposition}

\begin{proof}
We note that $G_{(i_0,j_0)}$ normalizes $G_{> (i_0,j_0)}$ either by \ref{Elm':2}, if $i_0 < j_0$ or by \ref{Elm':3} if $i_0 = j_0$. 

By the definition of $G_{\geq (i_0,j_0)}$ this group is generated by $G_{> (i_0,j_0)}$ and $G_{(i_0,j_0)}$. It follows that $G_{> (i_0,j_0)}$ is normal in $G_{\geq (i_0,j_0)}$ and that $G_{\geq (i_0,j_0)} = G_{ (i_0,j_0)} G_{> (i_0,j_0)}$. 

Last, by \pref{obs:trivial-intersect-obs}, the groups $G_{ (i_0,j_0)}$ and $ G_{> (i_0,j_0)}$ intersect trivially, and thus it follows that $G_{\geq (i_0,j_0)} = G_{(i_0,j_0)} \ltimes G_{> (i_0,j_0)}$.
\end{proof}

\ionote{This corollary was changed due to Inga's remark - now it also states the $S_{\geq (i,j)}$ generates $G_{\geq (i,j)}$ for every $1 \leq i \leq j \leq 5$. The previous version was just for $G_{\geq (1,1)}$ which is now the "in particular case"}

\begin{corollary}
\label{g as product of E's coro}
Every element $g \in \Unip_n (\mathcal{R})$ can be written as a product 
$$g = \prod_{1 \leq i \leq j \leq 5} E_{i,j} (A_{i,j}),$$
where the order of the elements in the product is according to the order $<$ on $\lbrace (i,j) : 1 \leq i \leq j \leq 5 \rbrace$ defined above.
\end{corollary}

\begin{proof}
We note that $\Unip_n (\mathcal{R}) = G_{\geq (1,1)}$. Fix $g \in G_{\geq (1,1)}$. By \pref{prop:semidirect-prop}, there are $g' \in G_{>(1,1)} = G_{\geq (2,2)}$ and $A_{1,1}$ such that $g = E_{1,1} (A_{1,1}) g'$. Continuing in the same manner, using \pref{prop:semidirect-prop} gives the needed result.
\end{proof}

\ydnote{Is \(|g|\) and \(|g|_{S}\) defined somewhere? I saw some explanation in the proof but maybe it should be in prelims to this section. Also, reading the proof made me think I don't fully understand the proposition itself. Is \(g'=g\) in the first case?}\ionote{The notation $\vert g \vert_S$ was defined aobve in the Prelim. in the subsection about Dehn functions, but maybe it was easy to miss. I made it into a proper definition in that subsection. Is it enough? Also, the proposition had a mistake - now corrected (I wrote $s^{-1} g' s$ where it should have been just $g'$. See if it makes sense now. } 

\begin{proposition}
\label{prop:change-length-in-Dehn}
Let $1 \leq i \leq j \leq 5, (i,j) \neq (1,5)$, $g \in G_{> (i,j)}$ and $s \in S_{(i,j)}$. The word $s^{-1} g s$ can be transformed to a word $g'$ in the alphabet $S_{>(i,j)}$ by applying at most $\vert g \vert_{S_{>(i,j)}} $ relations in $R_{\geq (i,j)}$ such that
$$\vert  g '  \vert_{S_{>(i,j)}} \leq 
\begin{cases}
    \vert g \vert_{S_{>(i,j)}} &  i =j \text{ or } (i,j) \geq (2,4) \\
    2  \vert g \vert_{S_{>(i,j)}} & i <j \text{ and } (i,j) \leq (1,3)
\end{cases}.$$
\end{proposition}

\begin{proof}
We denote $s = E_{i,j} (A)$ and $g = s_1 ... s_k$ where $s_1,...,s_k \in S_{>(i,j)}$ and $\vert g \vert_{S_{>(i,j)}} = k$. By free expansion moves, we have that 
$$s^{-1} g s = s^{-1} (s_1 s_2 ... s_k) s =(s^{-1} s_1 s) (s^{-1} s_2 s)  ... (s^{-1} s_k s).$$

Assume first that $(i,j) \geq (2,4)$. By \ref{Elm':2}, for every $s' = E_{i',j'} (A') \in S_{> (i,j)}$, it holds that 
$[E_{i',j'} (A'), E_{i,j} (A)]=e$ or equivalently, 
$$s^{-1} s' s = E_{i,j} (-A) E_{i',j'} (A') E_{i,j} (A) = E_{i',j'} (A ') [E_{i',j'} (A'), E_{i,j} (A)] = E_{i',j'} (A ').$$
Thus, by applying one relation for every $s^{-1} s_i s$, we deduce that $s^{-1} g s = g$ and in particular $\vert s^{-1} g s \vert_{S_{>(i,j)}} = \vert g \vert_{S_{>(i,j)}}$.

Next, we assume that $i=j$. By \ref{Elm':3}, for every $s' = E_{i',j'} (A') \in S_{> (i,j)}$ it holds that $s^{-1} s' s \in S_{> (i,j)}$. Thus, by applying one relation for every $s^{-1} s_i s$, we get a word in the alphabet $S_{> (i,j)}$ of length $\vert g \vert_{S_{>(i,j)}} $.

Last, we assume that $i <j$ and $(i,j) \leq (1,3)$. In that case, for $s' = E_{i',j'} (A') \in S_{> (i,j)}$ there are three possible cases. If $i' \neq j$, then by \ref{Elm':2}, $s^{-1} E_{i',j'} (A') s = E_{i',j'} (A')$. If $i' =j$, then by \ref{Elm':2}, $s^{-1} E_{i',j'} (A') s = E_{i',j'} (A') E_{i,j'} (-A A')$. If $j' =i$, then  by \ref{Elm':2}, $s^{-1} E_{i',j'} (A') s = E_{i',j'} (A') E_{i',j} ( A' A)$. Thus, by applying one relation for every $s^{-1} s_i s$, we get a word in the alphabet $S_{> (i,j)}$ with a length that is $\leq 2\vert g \vert_{S_{>(i,j)}} $.
\end{proof}

We define $\Dehn^{G_{\geq (i,j)}}$ to be the Dehn function of $G_{\geq (i,j)}$ with respect to the presentation $\langle S_{\geq (i,j)} \vert R_{\geq (i,j)} \rangle$ and similarly define $\Dehn^{G_{> (i,j)}}$, $\Dehn^{G_{(i,j)}}$. 

\begin{proposition}
\label{prop: recursive bound prop}
 We have the following recursive bounds on the Dehn functions:
 \begin{itemize}
    \item $\Dehn^{G_{(1,5)} (m)} \leq m$.
    \item For every $1 \leq i \leq j \leq 5$, $(i,j) < (1,5)$,
    $$\Dehn^{G_{\geq (i,j)}} (m)\leq 
    \begin{cases}
    2 m + \Dehn^{G_{> (i,j)}} (m) & i =j \text{ or } (i,j) \geq (2,4) \\
    2 m + \Dehn^{G_{> (i,j)}} (2m) & i < j, (i,j) \leq (1,3)
    \end{cases} .
    $$
\end{itemize}
\end{proposition}

\begin{proof}
We start by proving a generalized version of the first bound: Let $g \in G_{(i,j)}$, then by applying at most $\vert g \vert_{S_{(i,j)}} -1$ relations in $R_{(i,j)}$, $g$ can be reduced to an element in $S_{(i,j)}$. Indeed, if $\vert g \vert_{S_{(i,j)}} =l$, then $g = E_{i,j} (A_1) ... E_{(i,j)} (A_l)$. By applying \ref{Elm':1} $l-1$ times, $g$ can be reduced to 
$$\begin{cases}
    E_{i,j} (A_1 + ...+A_l) \in S_{(i,j)} & i <j \\
    E_{i,j} (A_1 ... A_l) \in S_{(i,j)} & i =j
\end{cases}.$$

As a special case, when $(i,j) = (1,5)$ and $w$ is a relation in $G_{(i,j)}$ with $\vert w \vert_{S_{(1,5)}} =m$, it can be reduced by applying $m-1$ relations to an element in $S_{(1,5)}$, i.e., to $e$. Thus, $\Dehn^{G_{(1,5)} (m)} \leq m$.

Assume that $(i,j) \neq (1,5)$ and let $w$ be a relation in $G_{\geq (i,j)}$ with $\vert w \vert_{S_{\geq (i,j)}} = m$. Such $w$ can be written as a product 
$$w = g_1 h_1 ... g_l h_l g_{l+1}$$
where $g_1,...,g_{l+1} \in G_{> (i,j)}$ (where $g_1, g_{l+1}$ may be trivial) and $h_1,...,h_l \in G_{(i,j)}$ and 
$$\vert g_1 \vert_{S_{> (i,j)}} + \vert h_1 \vert_{S_{(i,j)}} +... + \vert g_{l+1} \vert_{S_{> (i,j)}} = m.$$
By the previous bound, every $h_i$ can be reduced to $s_i \in S_{(i,j)}$ by applying at most $\vert h_i \vert_{S_{(i,j)}}-1$ relations. Thus, by applying $\leq \vert h_1 \vert_{S_{(i,j)}} + ... + \vert h_l \vert_{S_{(i,j)}}$ relations, $w$ can be reduced to 
$$w = g_1 s_1 ... g_l s_l g_{l+1}.$$
By a free expansion, we get that $w = g_1 s_1 ... s_{l-1} s_l (s_l^{-1} g_l s_l) g_{l+1} $ and by applying a relation of type \ref{Elm':1}, we can reduce $s_{l-1} s_l$ to $s_{l-1} ' \in S_{(i,j)}$. Thus, we get that 
$g_1 s_1 ... g_{l-1} s_{l-1} ' (s_l^{-1} g_l s_l) g_{l+1}$. Continuing in the same manner (applying another $l-1$ relations of type \ref{Elm':1}), we reduce the word to $s_1' ((s_1')^{-1} g_1 s_1') ((s_2 ')^{-1} g_2 s_2 ' ) ...((s_{l-1}')^{-1} g_{l-1} s_{l-1}') (s_l^{-1} g_l s_l) g_{l+1}$. 
By \pref{prop:change-length-in-Dehn}, for each $(s_i')^{-1} g_i s_i'$ we can apply $\vert g_i \vert_{S_{> (i,j)}}$ relations to change it into $g_i' \in G_{> (i,j)}$. To conclude, we applied at most
$$\vert h_1 \vert_{S_{(i,j)}} + ... + \vert h_l \vert_{S_{(i,j)}} + l + \vert g_1 \vert_{S_{> (i,j)}} + ... +\vert g_{l+1} \vert_{S_{> (i,j)}} \leq 2m$$
relations and changed the relation $w$ into the relation
$s_1' (g_1 '... g_{l+1} ')$. We note that $g_1 ' ... g_{l+1} ' \in G_{>(i,j)}$ and since $s_1' (g_1 '... g_{l+1} ')$ is a relation, it follows from \pref{prop:semidirect-prop}, that $s_1 ' =e$ and that $w' = g_1 ' ... g_{l+1} ' $ is a relation in $G_{>(i,j)}$. Therefore we deduce that
$$\Area_{R_{\geq (i,j)}} (w) \leq 2m + \Area_{R_{>(i,j)}} (w').$$
By \pref{prop:change-length-in-Dehn}, 
\begin{align*}
   & \vert w' \vert_{S_{>(i,j)}} = \vert g_1 ' ... g_{l+1} '  \vert_{S_{>(i,j)}} \leq \\
 &  \begin{cases}
   \sum_{i=1}^{l+1} \vert g_i \vert_{S_{>(i,j)}} &  i =j \text{ or } (i,j) \geq (2,4) \\
   \sum_{i=1}^{l+1} 2  \vert g_i \vert_{S_{>(i,j)}} & i <j \text{ and } (i,j) \leq (1,3)
\end{cases} \leq \\
& \begin{cases}
   m &  i =j \text{ or } (i,j) \geq (2,4) \\
   2m & i <j \text{ and } (i,j) \leq (1,3)
\end{cases}.
\end{align*}

The above bounds were for any arbitrary relation $w$ of length $\leq m$ and therefore
$$\Dehn^{G_{\geq (i,j)}} (m)\leq 
    \begin{cases}
    2 m + \Dehn^{G_{> (i,j)}} (m) & i =j \text{ or } (i,j) \geq (2,4) \\
    2 m + \Dehn^{G_{> (i,j)}} (2m) & i < j, (i,j) \leq (1,3)
    \end{cases} ,
    $$
as needed.
\end{proof}

\begin{corollary}
\label{Elm' dehn bound coro}
The Dehn function of $\Unip_n (\mathcal{R})$ with respect to the presentation $\langle S_{\Elm} ' \vert R_{\Elm} ' \rangle$ can be bounded as follows:
$$\Dehn_{S_{\Elm} ', R_{\Elm} '} (m) \leq 360 m .$$
\end{corollary}

\begin{proof}
We note that by definition $G_{\geq (1,1)} = \Unip_n (\mathcal{R})$. 

We spell out the order on the $(i,j)$:
\( (1,1)<(2,2)<(3,3)<(4,4)<(5,5)<(1,2)<(2,3)<(3,4)<(4,5)<(1,3)<(2,4)<(3,5)<(1,4)<(2,5)<(1,5).\)
Also, observe that for two consecutive $(i,j)<(i',j')$ in the list, we have $G_{>(i,j)}=G_{\geq (i',j')}$. Thus, iteratively using \pref{prop: recursive bound prop}, we get
\begin{align*}
    &\Dehn_{S_{\Elm} ', R_{\Elm} '} (m) = \Dehn^{G_{\geq (1,1)}} (m)  \\
    &\leq 2m + \Dehn^{G_{\geq (2,2)}} (m) \leq \ldots \leq 8m + \Dehn^{G_{\geq (5,5)}} (m) \leq 10m + \Dehn^{G_{\geq (1,2)}} (m) \\
    &\leq 10 m + 2m +  \Dehn^{G_{\geq (2,3)}} (2m) \leq 12m + 4m +  \Dehn^{G_{\geq (3,4)}} (4m) \leq 16m + 8m + \Dehn^{G_{\geq (4,5)}} (8m) \\
    &\leq 24m + 16m +\Dehn^{G_{\geq (1,3)}} (16m) \leq 40m + 32m + \Dehn^{G_{\geq (2,4)}} (32m) \leq 72m + 64m + \Dehn^{G_{\geq (3,5)}} (32m) \\
    &\leq 136m + 64m + \Dehn^{G_{\geq (1,4)}} (32m) \leq 200m + 64m + \Dehn^{G_{\geq (2,5)}} (32m) \\
    &\leq 264m+ 64m + \Dehn^{G_{\geq (1,5)}} (32m) \leq 328m + 32m = 360m.
\end{align*}

\end{proof}

\subsection{Proof of our main Dehn bound}

Corollary \ref{Elm' dehn bound coro} provides a bound for the Dehn function of $\Unip_n (\mathcal{R})$ with respect to the presentation  $\langle S_{\Elm} ' \vert R_{\Elm} ' \rangle$. As explained above, we will now show how this bound yields a bound for the Dehn function for the ``smaller'' presentation  $\langle S_{\Elm} \vert R_{\Elm} \rangle$ defined in \ref{subsection: Main result} above.   

We start by showing that $S_{\Elm}$ defined in \ref{subsection: Main result} is a generating set:

\begin{claim}
The set $S_{\Elm}$ generates $\Unip_n (\mathcal{R})$. Moreover, every element $g \in \Unip_n (\mathcal{R})$ can be written as a product 
$$g = \left(\prod_{1 \leq i \leq j \leq 5, (i,j) \neq (1,5)} E_{i,j} (A_{i,j}) \right) [E_{1,4} (A_{1,5}), E_{4,5} (I)],$$
where the order of the elements in the product is according to the order $<$ on $\lbrace (i,j) : 1 \leq i \leq j \leq 5 \rbrace$ defined above and $I$ is the $k_0 \times k_0$ identity matrix.
\end{claim}

\begin{proof}
We note that for every $(1,5)$-admissible $A_{1,5}$, it holds that $[E_{1,4} (A_{1,5}), E_{4,5} (I)] = E_{1,5} (A_{1,5})$ and use Corollary \ref{g as product of E's coro}.
\end{proof}

After that we can prove our main theorem in the section which gives a linear bound on the Dehn function of $\Unip_n (\mathcal{R})$ with respect to the presentation  $\langle S_{\Elm} \vert R_{\Elm} \rangle$:

\begin{proof}[Proof of \pref{thm:unip-presentation}]
Whenever we write $E_{i,j} (A)$, we mean that $A$ is $(i,j)$-admissible (without mentioning it explicitly). We also note that for $1 \leq i < j \leq 5$, if $i \neq 3, j \neq 3$, an $(i,j)$-admissible matrix is a square matrix and when denoting $E_{i,j} (I)$, we mean that $I$ is the identity matrix.

By Corollary \ref{Elm' dehn bound coro}, it is enough to show that there is a universal constant $C''$ such that $\Dehn_{ S_{\Elm}, R_{\Elm} } (m) \leq C'' \Dehn_{ S_{\Elm} ', R_{\Elm} ' } ( m)$. In order to do so, we will use \pref{prop: comparing dehn for dif pres prop}. Since $S_{\Elm} \subseteq S_{\Elm '}$ it follows that $\Free (S_{\Elm}) \subseteq \Free (S_{\Elm '})$ and we define $\phi : \Free (S_{\Elm}) \rightarrow \Free (S_{\Elm '})$ via this inclusion. Explicitly, for every $1 \leq i \leq j \leq 5$, $(i,j) \neq (1,5)$ and every $A$ that is $(i,j)$-admissible, we define $\phi (E_{i,j} (A)) = E_{i,j} (A)$. The map in the other direction $\phi ' :\Free (S_{\Elm} ') \rightarrow \Free (S_{\Elm })$ is induced from the map
$$\phi ' (E_{i,j} (A)) = 
\begin{cases}
E_{i,j} (A) & (i,j) \neq (1,5) \\
[E_{1,4} (A), E_{4,5} (I)] & (i,j) = (1,5)
\end{cases}.$$
We will bound the constants of \pref{prop: comparing dehn for dif pres prop} given by this map. Two bounds are very easy: first 
$$C_1 = \max_{s \in S_{\Elm}} \vert \phi (s) \vert_{S_{\Elm '}} = 1$$
and second
$$C_3 = \max_{s \in S_{\Elm}} \Area_{R_{\Elm}} (s (\phi ' (\phi (s)))^{-1}) = \max_{s \in S_{\Elm}} \Area_{R_{\Elm}} (s s^{-1}) =0.$$
Thus by \pref{prop: comparing dehn for dif pres prop} it is enough to bound $C_2 = \max_{\rho ' \in R_{\Elm '}} \Area_{R_{\Elm}} (\phi ' (\rho '))$, and if we show that $C_2$ can be bounded universally, then for $C '' = C_2$ it will hold that $\Dehn_{ S_{\Elm}, R_{\Elm} } (m) \leq C'' \Dehn_{ S_{\Elm} ', R_{\Elm} ' } ( m) \leq C'' 360 m$. 

Next, we will go systematically over all the relations $\rho ' \in R_{\Elm '}$ and show that for every such $\rho '$ the $\Area_{R_{\Elm}} (\phi ' (\rho '))$ can be bounded independently of $n, \mathcal{R}$. We recall that $\Area_{R_{\Elm}} (\phi ' (\rho '))$ is the number of relations in $R_{\Elm}$ one needs to to make $\phi ' (\rho ')$ into the trivial word. This will also show that $\langle S_{\Elm} \vert R_{\Elm} \rangle$ is indeed a presentation of $\Unip_n (\mathcal{R})$.

We will not keep track of the actual area book-keeping, but only show that every relation $\phi ' (\rho ')$ can be deduced from the relations in $R_{\Elm}$ (and that the number of relations in $R_{\Elm}$ needed is independent of $n, \mathcal{R}$). 

We note that every $\rho ' \in R_{\Elm} \subseteq R_{\Elm '}$ is trivially deduced from $R_{\Elm}$.

\textbf{Step 1:} We will show how to deduce the relations in $ R_{\Elm '}$ that involve only $E_{i,j} (A)$'s in which $i < j$ and $(i,j) \neq (1,5)$ (and are not included in $R_{\Elm}$). 

\paragraph*{$[E_{1,2} (A), E_{3,5} (B)] =e$:} By \ref{Elm:2}, $E_{3,5} (B) = [E_{3,4} (B), E_{4,5} (I)]$ and also by \ref{Elm:2}, $[E_{1,2} (A), E_{3,4} (B)] = [E_{1,2} (A), E_{4,5} (I)] =e$. It follows from \pref{claim:comm ident claim} that $[E_{1,2} (A), E_{3,5} (B)] =e$.

\paragraph*{$[E_{1,3} (A), E_{4,5} (B)] =e$:} By \ref{Elm:2}, $E_{1,3} (A) = [E_{1,2} (I), E_{2,3} (A)]$ and also by \ref{Elm:2}, $[E_{1,2} (I), E_{4,5} (B)] = [E_{2,3} (A), E_{4,5} (B)] =e$. It follows from \pref{claim:comm ident claim} that $[E_{1,3} (A), E_{4,5} (B)] =e$.

\paragraph*{$[E_{1,3} (A), E_{2,5} (B)] =e$:} By \ref{Elm:2}, $E_{2,5} (B) = [E_{2,4} (I), E_{4,5} (B)]$ and also  by \ref{Elm:2}, $[E_{1,3} (A), E_{2,4} (I)] =e$. We have shown above that $[E_{1,3} (A), E_{4,5} (B)]=e$ and thus it follows from \pref{claim:comm ident claim} that $[E_{1,3} (A), E_{2,5} (B)] =e$.

\paragraph*{$[E_{1,4} (A), E_{3,5} (B)] =e$:} By \ref{Elm:2}, $E_{1,4} (A) = [E_{1,2} (I), E_{2,4} (A)]$ and also by \ref{Elm:2}, $[E_{2,4} (A), E_{3,5} (B)] =e$.  We have shown above that $[E_{1,2} (I), E_{3,5} (B)]=e$ and thus it follows from \pref{claim:comm ident claim} that $[E_{1,4} (A), E_{3,5} (B)] =e$.

\paragraph*{$[E_{1,4} (A), E_{2,5} (B)] =e$:} We define $B_1$ to be the $k_0 \times k_1$ matrix given in block form $B_1 = \begin{pmatrix}
    B & 0_{k_0 \times (k_1 - k_0)}
\end{pmatrix}$. We also define $B_2$ to be the $k_1 \times k_0$ matrix given in block form $B_2 = \begin{pmatrix}
    I_{k_0 \times k_0} \\
    0_{(k_1 - k_0) \times k_0}
\end{pmatrix}$. By \ref{Elm:2}, $E_{2,5} (B) = [E_{2,3} (B_1), E_{3,5} (B_2)]$ and also by \ref{Elm:2}, $[E_{1,4} (A), E_{2,3} (B_1)] = e$. We have shown above that $[E_{1,4} (A), E_{3,5} (B_2)]=e$ and thus it follows from \pref{claim:comm ident claim} that $[E_{1,4} (A), E_{2,5} (B)] =e$.

\textbf{Step 2:} We will show how to deduce the relations in $ R_{\Elm '}$ of the form $\phi ' ([E_{1,5} (A), E_{i,j} (B)])$ for every $1 \leq i < j \leq 5, (i,j) \neq (1,5), (4,5)$ and every $A,B$, i.e., by our definition of $\phi '$, we need to show how to deduce the relations of the form $[[E_{1,4} (A), E_{4,5} (I)],E_{i,j} (B)] =e$ for every $1 \leq i < j \leq 5, (i,j) \neq (1,5), (4,5)$ and every $A,B$. We will actually show more: we will show how to deduce the relations of the form $[[E_{1,4} (A_1), E_{4,5} (A_2)],E_{i,j} (B)] =e$ for every $1 \leq i < j \leq 5, (i,j) \neq (1,5), (4,5)$ and every $A_1, A_2, B$. 

\paragraph*{$[[E_{1,4} (A_1), E_{4,5} (A_2)],E_{i,j} (B)] =e$ for $(i,j) \in \lbrace (1,2), (1,3), (2,3), (2,5), (3,5)$:} For every $(i,j) \in \lbrace (1,2), (1,3), (2,3), (2,5), (3,5)$, we have either by \ref{Elm:2} or by the previous step that $[E_{1,4} (A_1), E_{i,j} (B)] = [E_{4,5} (A_2), E_{i,j} (B)] =e$ and thus it follows from \pref{claim:comm ident claim} that $[[E_{1,4} (A_1), E_{4,5} (A_2)],E_{i,j} (B)] =e$.

\paragraph*{$[[E_{1,4} (A_1), E_{4,5} (A_2)],E_{3,4} (B)] =e$:} By Hall-Witt identity:

\begin{align*}
   & [[E_{1,4} (A_1), E_{4,5} (-A_2)],E_{3,4} (B)]^{E_{4,5} (A_2)} \cdot \\
   & \cdot [[E_{4,5} (A_2),E_{3,4} (-B)], E_{1,4} (A_1)]^{E_{3,4} (B)} \cdot \\
   & \cdot [[E_{3,4} (B), E_{1,4} (-A_1)], E_{4,5} (A_2)]^{E_{1,4} (A_1)} =e.
\end{align*}

By \ref{Elm:2}, $[E_{3,4} (B), E_{1,4} (-A_1)] =e$ and thus $[[E_{3,4} (B), E_{1,4} (-A_1)], E_{4,5} (A_2)]^{E_{1,4} (A_1)} =e$. By \ref{Elm:2}, $[E_{4,5} (A_2),E_{3,4} (-B)] = E_{3,5} (B A_2 )$. By the previous step, $[E_{3,5} (B A_2 ), E_{1,4} (A_1)]=e$ and thus
$$[[E_{4,5} (A_2),E_{3,4} (-B)], E_{1,4} (A_1)]^{E_{3,4} (B)} = [E_{3,5} (A_2 B), E_{1,4} (A_1)]^{E_{3,4} (B)} =e. $$
It follows from the Hall-Witt identity that $[[E_{1,4} (A_1), E_{4,5} (-A_2)],E_{3,4} (B)]^{E_{4,5} (A_2)} =e$ which implies that $[[E_{1,4} (A_1), E_{4,5} (-A_2)],E_{3,4} (B)] =e$. Since this is true for every $A_1, A_2, B$, we can replace $A_2$ with $- A_2$ and get $[[E_{1,4} (A_1), E_{4,5} (A_2)],E_{3,4} (B)] =e$ as needed.

\paragraph*{$[[E_{1,4} (A_1), E_{4,5} (A_2)],E_{i,4} (B)] =e$ for $i=1,2$:} We define $B_1$ to be the $k_0 \times k_1$ matrix given in block form $B_1 = \begin{pmatrix}
    B & 0_{k_0 \times (k_1 - k_0)}
\end{pmatrix}$. We also define $B_2$ to be the $k_1 \times k_0$ matrix given in block form $B_2 = \begin{pmatrix}
    I_{k_0 \times k_0} \\
    0_{(k_1 - k_0) \times k_0}
\end{pmatrix}$.
By \ref{Elm:2}, $[E_{i,3} (B_1), E_{3,4} (B_2)] = E_{i,4} (B)$. We showed above that $[[E_{1,4} (A_1), E_{4,5} (A_2)],E_{i,3} (B_1)] = [[E_{1,4} (A_1), E_{4,5} (A_2)],E_{3,4} (B_2)] =e$ and thus by \pref{claim:comm ident claim} it follows that $[[E_{1,4} (A_1), E_{4,5} (A_2)],E_{i,4} (B)] =e$.

\textbf{Step 3:}  We will show how to deduce the relations in $ R_{\Elm '}$ of the form $\phi ' ([E_{1,j} (A), E_{j,5} (B)] (E_{1,5} (AB))^{-1})$ for every $2 \leq j \leq 4$ and every $A,B$. By the definition of $\phi '$, we need to show how to deduce $[E_{1,4} (AB), E_{4,5} (I)] = [E_{1,j} (A), E_{j,5} (B)]$ for every $2 \leq j \leq 4$ and every $A,B$.

\paragraph*{$[E_{1,4} (A_1 A_2), E_{4,5} (A_3)] = [E_{1,j} (A_1), E_{j,5} (A_2 A_3)] $ for every $j =2,3$ and every $A_1,A_2,A_3$}: We denote $g_1 = E_{1,j} (A_1), g_2 = E_{j,4} (-A_2)$ and $g_3 = E_{4,5} (A_3)$ and use the Hall-Witt identity:
$$[[g_1, g_2^{-1}],g_3]^{g_2} [[g_2, g_3^{-1}],g_1]^{g_3} [[g_3, g_1^{-1}],g_2]^{g_1} =e.$$
By \ref{Elm:2} and our choice of $g_1, g_2, g_3$ it holds that: $[g_1, g_2^{-1}] = E_{1,4} (A_1 A_2)$, $[g_2, g_3^{-1}] = E_{j,5} (A_2 A_3)$ and $[g_3, g_1^{-1}] = e$. Thus, the Hall-Witt identity gives:
$$[E_{1,4} (A_1 A_2),E_{4,5} (A_3)]^{E_{j,4} (-A_2)} [E_{j,5} (A_2 A_3),E_{1,j} (A_1)]^{E_{4,5} (A_3)} = e.$$
We showed above that $[[E_{1,4} (A_1 A_2)],E_{4,5} (A_3)]$ commutes with $E_{j,4} (-A_2)$ and thus $[[E_{1,4} (A_1 A_2)],E_{4,5} (A_3)]^{E_{j,4} (-A_2)} = [[E_{1,4} (A_1 A_2)],E_{4,5} (A_3)]$. Also, by \ref{Elm:2}, $[E_{j,5} (A_2 A_3), E_{4,5} (A_3)] = [E_{1,j} (A_1), E_{4,5} (A_3)] =e$ and thus $[E_{j,5} (A_2 A_3),E_{1,j} (A_1)]^{E_{4,5} (A_3)}= [E_{j,5} (A_2 A_3),E_{1,j} (A_1)]$. Therefore we get 
$$[E_{1,4} (A_1 A_2),E_{4,5} (A_3)] [E_{j,5} (A_2 A_3),E_{1,j} (A_1)] = e$$
or equivalently,
$[E_{1,4} (A_1 A_2), E_{4,5} (A_3)] = [E_{1,j} (A_1), E_{j,5} (A_2 A_3)]$.

\paragraph*{$[E_{1,4} (AB), E_{4,5} (I)] = [E_{1,j} (A), E_{j,5} (B)]$ for every $j=2,3$ and every $A,B$}: We showed that $[E_{1,4} (A_1 A_2), E_{4,5} (A_3)] = [E_{1,j} (A_1), E_{j,5} (A_2 A_3)] $. Take $A_1 = A, A_2 =B, A_3 =I$.

\paragraph*{$[E_{1,4} (AB), E_{4,5} (I)] = [E_{1,4} (A), E_{4,5} (B)]$ for every $A,B$}: We showed that $[E_{1,4} (A_1 A_2), E_{4,5} (A_3)] = [E_{1,2} (A_1), E_{2,5} (A_2 A_3)] $. For $A_1 = I, A_2 =A, A_3 =B$, then $[E_{1,4} (A), E_{4,5} (B)] = [E_{1,2} (I), E_{2,5} (AB)] $. In order to conclude the proof, we use the fact that we already shown that $[E_{1,2} (I), E_{2,5} (AB)] = [E_{1,4} (AB), E_{4,5} (I)]$.

\textbf{Step 4:} We will show how to deduce the relations in $ R_{\Elm '}$ that involve only $E_{i,j}$'s with $i <j$. We have two such relations: First, $\phi '  ([E_{1,5} (A), E_{4,5} (B)])$ or equivalently by the definition of $\phi '$, $[[E_{1,4} (A), E_{4,5} (I)], E_{4,5} (B)] =e$. Second, $\phi '  ((E_{1,5} (A) E_{1,5} (B)) (E_{1,5} (A+B) )^{-1})$ or equivalently by the definition of $\phi '$, $[E_{1,4} (A), E_{4,5} (I)] [E_{1,4} (B), E_{4,5} (I)]  =[E_{1,4} (A+B), E_{4,5} (I)]$.

\paragraph*{$[[E_{1,4} (A), E_{4,5} (I)], E_{4,5} (B)] =e$:} We showed above that $[E_{1,4} (A), E_{4,5} (I)] = [E_{1,2} (A), E_{2,5} (I)]$. By \ref{Elm:2}, $[E_{1,2} (A), E_{4,5} (B)] = [E_{2,5} (I), E_{4,5} (B)]=e$. It follows from \pref{claim:comm ident claim} that $[[E_{1,4} (A), E_{4,5} (I)], E_{4,5} (B)] =e$.

\paragraph*{$[E_{1,4} (A), E_{4,5} (I)] [E_{1,4} (B), E_{4,5} (I)]  =[E_{1,4} (A+B), E_{4,5} (I)]$:} We showed above that $E_{1,4} (B)$ commutes with $[E_{1,4} (A), E_{4,5} (I)]$ and therefore by Proposition \ref{BissD result prop} and \ref{Elm:1},
$$ [E_{1,4} (A), E_{4,5} (I)] [E_{1,4} (B), E_{4,5} (I)] =  [E_{1,4} (A) E_{1,4} (B), E_{4,5} (I)] = [E_{1,4} (A+B), E_{4,5} (I)]$$
as needed.

\textbf{Step 5:} We will show how to deduce the relations in $ R_{\Elm '}$ that involve $E_{l,l}$'s and $E_{i,j}$'s with $(i,j) \neq (1,5)$ and do not appear in $ R_{\Elm}$.

\paragraph*{$(E_{1,1} (B))^{-1} E_{i,5} (A) E_{1,1} (B) = E_{i,5} (A)$ for $i=2,3$:} By \ref{Elm:2}, $E_{i,5} (A) = [E_{i,4} (A), E_{4,5} (I)]$. By free expansion, 
$$(E_{1,1} (B))^{-1} [E_{i,4} (A), E_{4,5} (I)] E_{1,1} (B) = [(E_{1,1} (B))^{-1} E_{i,4} (A) E_{1,1} (B), (E_{1,1} (B))^{-1} E_{4,5} (I) E_{1,1} (B)].$$
By \ref{Elm:3}, $(E_{1,1} (B))^{-1} E_{i,4} (A) E_{1,1} (B) = E_{i,4} (A)$ and $ (E_{1,1} (B))^{-1} E_{4,5} (I) E_{1,1} (B) = E_{4,5} (I)$ and thus
$$[(E_{1,1} (B))^{-1} E_{i,4} (A) E_{1,1} (B), (E_{1,1} (B))^{-1} E_{4,5} (I) E_{1,1} (B)] = [E_{i,4} (A), E_{4,5} (I)]$$
as needed.

\paragraph*{$(E_{5,5} (B))^{-1} E_{1,j} (A) E_{5,5} (B) = E_{1,j} (A)$ for $j=3,4$:} By \ref{Elm:2}, $E_{1,j} (A) = [E_{1,2} (I), E_{2,j} (A)]$. By free expansion, 
$$(E_{5,5} (B))^{-1} [E_{1,2} (I), E_{2,j} (A)] E_{5,5} (B) = [(E_{5,5} (B))^{-1} E_{1,2} (I) E_{5,5} (B), (E_{5,5} (B))^{-1} E_{2,j} (A) E_{5,5} (B)].$$
By \ref{Elm:3}, $(E_{5,5} (B))^{-1} E_{1,2} (I) E_{5,5} (B) = E_{1,2} (I)$ and $  (E_{5,5} (B))^{-1} E_{2,j} (A) E_{5,5} (B)= E_{2,j} (A)$ and thus
$$ [(E_{5,5} (B))^{-1} E_{1,2} (I) E_{5,5} (B), (E_{5,5} (B))^{-1} E_{2,j} (A) E_{5,5} (B)] =  [ E_{1,2} (I)),  E_{2,j} (A)]$$
as needed.

\textbf{Step 6:}  We will show how to deduce the relations in $ R_{\Elm '}$ that involve $\phi ' (E_{1,5} (A))$. By our definition of $\phi '$, we need to show for every $A,B$ that 
$$(E_{l,l} (B))^{-1} [E_{1,4} (A), E_{4,5} (I)] E_{l,l} (B) = 
\begin{cases}
[E_{1,4} (A), E_{4,5} (I)] & 2 \leq l \leq 4 \\
[E_{1,4} (B^{-1} A), E_{4,5} (I)] & l=1 \\
[E_{1,4} (AB), E_{4,5} (I)] & l=5 
\end{cases}.$$
We note that by free expansion, 
$$(E_{l,l} (B))^{-1} [E_{1,4} (A), E_{4,5} (I)] E_{l,l} (B) =  [(E_{l,l} (B))^{-1} E_{1,4} (A) E_{l,l} (B),(E_{l,l} (B))^{-1} E_{4,5} (I) E_{l,l} (B)].$$

To finish, we will go case by case: 

For $l=2,3$, we have by \ref{Elm:3}, that $(E_{l,l} (B))^{-1} E_{1,4} (A) E_{l,l} (B) = E_{1,4} (A)$ and $(E_{l,l} (B))^{-1} E_{4,5} (I) E_{l,l} (B) = E_{4,5} (I)$ and thus we are done. 

For $l=1$, we have by \ref{Elm:3}, that $(E_{1,1} (B))^{-1} E_{1,4} (A) E_{1,1} (B) = E_{1,4} (B^{-1} A)$ and $(E_{1,1} (B))^{-1} E_{4,5} (I) E_{1,1} (B) = E_{4,5} (I)$ and thus we get $[E_{1,4} (B^{-1} A), E_{4,5} (I)]$ as needed.

For $l=5$, we have by \ref{Elm:3}, that $(E_{5,5} (B))^{-1} E_{1,4} (A) E_{5,5} (B) = E_{1,4} (A)$ and $(E_{5,5} (B))^{-1} E_{4,5} (I) E_{5,5} (B) = E_{4,5} (B)$ and thus we get $[E_{1,4} (A), E_{4,5} (B)]$ and by Step 3, this is equal to $[E_{1,4} (AB), E_{4,5} (I)]$ as needed.

For $l=4$, we have by \ref{Elm:3}, that $(E_{4,4} (B))^{-1} E_{1,4} (A) E_{4,4} (B) = E_{1,4} (AB)$ and $(E_{4,4} (B))^{-1} E_{4,5} (I) E_{4,4} (B) = E_{4,5} (B^{-1})$ and thus we get $[E_{1,4} (AB), E_{4,5} (B^{-1})]$ and by Step 3, this is equal to $[E_{1,4} (AB B^{-1}), E_{4,5} (I)] = [E_{1,4} (A), E_{4,5} (I)]$ as needed.
\ivnote{Is there a strategic way to see that these are all relations? Does this implicitly show that $\langle S_{\Elm}\mid R_{\Elm}\rangle \cong \Unip$? Because I don't think we mention anywhere that we haven proven this, only state at some point that we still need to prove it.}\ionote{Added Proposition \ref{prop: S' R' pres} above}
\end{proof}

\subsection[Element decomposition in the unipotent group]{Element decomposition in $\Unip_n (\mathcal{R})$}
Below, we will also need the following result (which is not about the Dehn function):

\begin{proposition}\label{prop:diameter-of-large-triangle-matrix-groups}
Let $H_1, H_2 < \Unip_n (\mathcal{R})$ be the following groups:
$$H_1 = \langle E_{i,j} (A) : 1 \leq i \leq j \leq 4, A \text{ is } (i,j) \text{-admissible} \rangle,$$
$$H_2 = \langle E_{i,5} (A) : 2 \leq i \leq 5, A \text{ is } (i,j) \text{-admissible} \rangle.$$
Then every element $g \in \Unip_n (\mathcal{R})$ can be written as a product of at most $6$ elements of $H_1 \cup H_2$.
\end{proposition}

\begin{proof}
We note that by definition, $H_1$  is composed of all the matrices of the form:
$$\begin{pmatrix}
    A_{1,1} & A_{1,2} & A_{1,3} & A_{1,4} & 0_{k_0}\\
    0_{k_0} & A_{2,2} & A_{2,3} & A_{2,4} & 0_{k_0} \\
    0_{k_1 \times k_0} & 0_{k_1 \times k_0} & A_{3,3} & A_{3,4} & 0_{k_1 \times k_0} \\
    0_{k_0} & 0_{k_0} & 0_{k_0 \times k_1} & A_{4,4} & 0_{k_0} \\
    0_{k_0} & 0_{k_0} & 0_{k_0 \times k_1} & 0_{k_0} & I
\end{pmatrix}$$
and $H_2$ is composed of all the matrices of the form:
$$\begin{pmatrix}
    I &  0_{k_0} & 0_{k_1 \times k_0} &  0_{k_0} & 0_{k_0}\\
    0_{k_0} & I & 0_{k_1 \times k_0} &  0_{k_0} & A_{2,5}\\
    0_{k_1 \times k_0} & 0_{k_1 \times k_0} & I & 0_{k_1 \times k_0} & A_{3,5} \\
    0_{k_0} & 0_{k_0} & 0_{k_0 \times k_1} & I & A_{4,5} \\
    0_{k_0} & 0_{k_0} & 0_{k_0 \times k_1} & 0_{k_0} & A_{5,5}
\end{pmatrix}.$$
Therefore, for every element $g \in \Unip_n (\mathcal{R})$ there are elements $h_1 \in H_1, h_2 \in H_2$ and a matrix $A$ such that $g = h_1 h_2 E_{1,5} (A)$. We finish by $E_{1,5} = [E_{1,2} (I), E_{2,5} (A)]$ and noting that $E_{1,2} (I) \in H_1$ and $E_{2,5} (A) \in H_2$. 
\end{proof}

\section{The direct product lemma}\label{sec:dirprod}
For \(B,I \subseteq [d]\), we use the notation
\[H_I^B \coloneqq U_{B \setminus I}.\]
Thus, \(H_I^B\) is the subgroup generated by the elementary unipotent matrices indexed by \(B \setminus I\). We use this notation throughout this section.

For the rest of this section, let \(\cI\) be a nice set of indexes and let \(J \in \cI\). The complex we analyze is
\[Y \coloneqq \CC\left(U_{[d]\setminus J};\left(H_I^{[d]\setminus J}\right)_{I \in \cI \setminus \set{J}}\right).\]

Let \(B_0,B_1,\dots,B_r\) be the \(J\)-bins. The block decomposition gives the internal direct product decompositions
\[U_{[d]\setminus J} = \prod_{\ell=0}^r U_{B_\ell} \qquad \text{and} \qquad H_I^{[d]\setminus J} = \prod_{\ell=0}^r H_I^{B_\ell} \quad \text{for every } I \in \cI \setminus \set{J}.\]
For every \(0 \leq \ell \leq r\), denote by
\[Y_\ell \coloneqq \CC\left(U_{B_\ell};\left(H_I^{B_\ell}\right)_{I \in \cI \setminus \set{J}}\right)\]
the coset complex corresponding to the \(\ell\)-th factor of the direct product.

We recall which presentations are used in the Dehn functions below. For every \(I \in \cI \setminus \set{J}\), let \(R_I\) be the multiplication-table relations of \(H_I^{[d]\setminus J}\), and let \(R_I^\ell\) be the multiplication-table relations of \(H_I^{B_\ell}\). In the notation of \pref{def: dehn function of coset complex}, let
\[\Dehn_Y = \Dehn_{S,R}, \qquad S = \bigcup_{I \in \cI \setminus \set{J}} H_I^{[d]\setminus J}, \qquad R = \bigcup_{I \in \cI \setminus \set{J}} R_I,\]
whereas
\[\Dehn_{Y_\ell} = \Dehn_{S_\ell,R_\ell}, \qquad S_\ell = \bigcup_{I \in \cI \setminus \set{J}} H_I^{B_\ell}, \qquad R_\ell = \bigcup_{I \in \cI \setminus \set{J}} R_I^\ell.\]
Thus, each \(\Dehn_{Y_\ell}\) is computed using the presentation associated with the component complex \(Y_\ell\), while \(\Dehn_Y\) is computed using the presentation associated with \(Y\).

%\ionote{There might be a small issue here: $Y$ is defined as a coset complex and $\Dehn_{Y}$ was defined as the Dehn function of the group acting on $Y$ using the presentation coming from the vertex stabilizers. However, we need the direct product decomposition at the group level, i.e., at some point we need to say/prove that the direct product decomposition of the coset complex is equivalent to the direct product decomposition at the group level}

\begin{lemma}\label{lem:direct-product-of-groups}
    Assume that there is a constant \(A\) such that
    \[\Dehn_{Y_\ell}(s) \leq A s\]
    for every \(0 \leq \ell \leq r\) and every \(s \in \N\), where each Dehn function is taken with respect to the coset-complex presentation of \(Y_\ell\) described above. Then, with respect to the coset-complex presentation of \(Y\),
    \[\Dehn_Y(t) \leq O((t(r+1))^2) + (r+1) A t\]
    for every \(t \in \N\).
\end{lemma}

We will start with the following more basic claim that we use as a building block. We will say that an upper-triangular matrix $g$ is \textit{supported on block size $k$}, if 
$$g = 
\begin{pmatrix}
        I_{k'} & 0 & 0 \\
        0 & A & 0 \\
        0 & 0 & I_{k''}
    \end{pmatrix},$$
where $A$ is a $k \times k$ upper-triangular matrix. We will call $A$ \textit{the support of $g$}.
\begin{claim}\label{claim:diameter-bound-for-good-groups}
    Let $0 \leq \ell \leq r$ and denote $B_{\ell} =[a,b]$. There is a universal polynomial $p$ (independent of $\ell, r, [a,b]$) such that for every $C \in \mathbb{N}, 20 \geq C \geq 1$, every \(I\in\cI\setminus\set{J}\) and every $g \in H_I^{B_\ell}$, there are $g_1,...,g_{s} \in H_I^{B_\ell}$ such that $s \leq p (C)$, $g = g_1 ... g_{s}$, and each $g_i$ is supported on block size $\leq \frac{b-a}{C} $. Moreover, $g$ can be transformed into $g_1 ... g_{s}$ by applying at most $p(C)$ relations in $R_I^\ell$.  
\end{claim}
\ydnote{Consider changing the statement so that \(g,g_1,g_2,\dots,g_s \in H_I^{B_\ell}\). That way it is clear the number of operations is \(O(s)\) since we can just use the multiplication table of \(H_I^{B_\ell}\) which is in \(R_\ell\).}
\ivnote{Sounds good to me, since any $g\in S_\ell$ will be in one of the $H_I^{B_\ell}$ and then the factors will also be in there. We might just need basically what is now the first sentence of the proof to later in Lemma 7.2 when we use this claim again. } \ionote{See corrected version}

\begin{proof}
Fix \(0 \leq \ell \leq r\), denote \(B_\ell=[a,b]\). Also fix \(I\in\cI\setminus\set{J}\) and \(g\in H_I^{B_\ell}\).

We first observe that if \(g\in H_I^{B_\ell}\) is decomposed into a product \(g=g_1\cdots g_6\) of matrices in \(U_{B_\ell}\), we can assume that all six matrices belong to \(H_I^{B_\ell}\). Indeed, replace each \(g_i\) by the matrix obtained by setting to zero all entries outside the diagonal blocks determined by \(I\cap B_\ell\). Since the diagonal blocks of a product of upper-triangular matrices are the products of the corresponding diagonal blocks, the resulting matrices still multiply to \(g\). Moreover, this operation yields matrices in \(H_I^{B_\ell}\).

We show by induction on \(m\geq 0\) that \(g\) can be transformed into a product of at most \(6^m\) elements of \(H_I^{B_\ell}\), each supported on block size at most \(\max\{(b-a+2)(9/10)^m,5\}\), using at most \(6^m-1\) relations in \(R_I^\ell\).

For \(m=0\), we simply take \(g\) itself. Assume that the statement holds for \(m\). Consider one of the resulting matrices \(h\), supported on block size \(k\). If \(k\leq 5\), we leave it unchanged. Otherwise, by \pref{prop:diameter-of-large-triangle-matrix-groups}, \(h\) can be decomposed into at most six matrices supported on block size at most
\[k-\lfloor k/5\rfloor\leq 9k/10.\]

By the observation above, we can take all these matrices to belong to \(H_I^{B_\ell}\). Thus, replacing \(h\) by their product requires at most five relations from \(R_I^\ell\).

Applying this to every factor, we obtain a product of at most \(6^{m+1}\) elements, each supported on block size at most \(\max\{(b-a+2)(9/10)^{m+1},5\}\). The total number of relations used is at most
\[(6^m-1)+5\cdot 6^m=6^{m+1}-1.\]
This completes the induction.

Finally, take $m=\lceil\log_{10/9}C\rceil +1 $. The resulting factors have block size at most \(\max\{\frac{9}{10}(b-a+2)/C,5\} \leq^{b-a \geq 100} \frac{b-a}{C} \), and both their number and the number of relations used are bounded by
\[6^m\leq 6^2 C^{18}.\]
Thus we can take \(p(C)=6^2 C^{18}\).
\end{proof}

\begin{proof}[Proof of \pref{lem:direct-product-of-groups}]
    Let \(w=w_1 w_2 \dots w_s\) be a word over \(S\) such that \(s \leq t\) and \(w=1\) in \(U_{[d]\setminus J}\). For every \(1 \leq j \leq s\), choose \(I_j \in \cI \setminus \set{J}\) such that
    \[w_j \in H_{I_j}^{[d]\setminus J}.\]
    By the direct product decomposition of \(H_{I_j}^{[d]\setminus J}\), we can write
    \[w_j = w_j^0 w_j^1 \dots w_j^r\]
    with
    \[w_j^\ell \in H_{I_j}^{B_\ell}\]
    for every \(0 \leq \ell \leq r\).

    Every partial product of the factors \(w_j^0,w_j^1,\dots,w_j^r\) belongs to \(H_{I_j}^{[d]\setminus J}\). Therefore, using the multiplication-table relations of \(H_{I_j}^{[d]\setminus J}\), we can replace \(w_j\) by \(w_j^0 w_j^1 \dots w_j^r\) in \(O(r)\) relation steps. Applying this to every letter transforms \(w\), in \(O(sr)\) steps, into
    \[w_1^0 w_1^1 \dots w_1^r w_2^0 w_2^1 \dots w_2^r \dots w_s^0 w_s^1 \dots w_s^r.\]

    We next sort the letters according to their superscripts, obtaining
    \[(w_1^0 w_2^0 \dots w_s^0)(w_1^1 w_2^1 \dots w_s^1) \dots (w_1^r w_2^r \dots w_s^r).\]
    We use a bubble-sort procedure, repeatedly commuting adjacent letters whose superscripts are in the wrong order. The following claim shows that each such commutation costs only a constant number of relation steps (its proof is in the end of this section).

    \begin{claim}\label{claim:commuting-a-pair}
        Let \(x=w_i^j\) and \(y=w_{i'}^{j'}\) for \(j \ne j'\). Then one can transform \(xy\) into \(yx\) using \(O(1)\) relations from \(R\).
    \end{claim}

    There are \(s(r+1)\) letters, and hence at most \(O((s(r+1))^2)\) adjacent commutations. By \pref{claim:commuting-a-pair}, the sorting therefore uses \(O((s(r+1))^2)\) relation steps.

    For every \(0 \leq \ell \leq r\), denote
    \[w^{(\ell)} \coloneqq w_1^\ell w_2^\ell \dots w_s^\ell.\]
    Since
    \[U_{[d]\setminus J} = \prod_{\ell=0}^r U_{B_\ell}\]
    and \(w=1\) in \(U_{[d]\setminus J}\), each \(w^{(\ell)}\) is equal to the identity in \(U_{B_\ell}\). It is a word of length at most \(s \leq t\) over \(S_\ell\), and therefore
    \[\Area_{R_\ell}(w^{(\ell)}) \leq \Dehn_{Y_\ell}(t) \leq A t.\]

    Every relation in \(R_\ell\) is a multiplication-table relation inside some
    \[H_I^{B_\ell} \leq H_I^{[d]\setminus J},\]
    and hence is also a relation in \(R\). Thus, the same reduction contracts \(w^{(\ell)}\) inside the presentation defining \(\Dehn_Y\). Contracting all \(r+1\) component words costs at most \((r+1)At\) additional relation steps. Altogether,
    \[\Area_R(w) \leq O((t(r+1))^2) + (r+1)At.\]
    Taking the maximum over all trivial words of length at most \(t\) proves the lemma.

    It remains to prove \pref{claim:commuting-a-pair}.

    \begin{proof}[Proof of \pref{claim:commuting-a-pair}]
        Write
        \[B_j=[a,b], \qquad B_{j'}=[a',b'], \qquad c=b-a, \qquad c'=b'-a'.\]
        By \pref{claim:diameter-bound-for-good-groups}, we can decompose
        \[x=x_1 x_2 \dots x_{C_0} \qquad \text{and} \qquad y=y_1 y_2 \dots y_{C_1},\]
        where \(C_0,C_1=O(1)\), every \(x_m\) is supported on a block of size at most \(\frac{c}{20}\), and every \(y_{m'}\) is supported on a block of size at most \(\frac{c'}{20}\). These decompositions can be performed using \(O(1)\) multiplication-table relations in \(H_{I_i}^{[d]\setminus J}\) and \(H_{I_{i'}}^{[d]\setminus J}\), respectively.

        For every \(x_m\), there is some \(q \in \set{0,1,2,3}\) such that the indexes of its support are contained in
        \[[a + q \lfloor \tfrac{c}{5} \rfloor,\min \set{a+(q+2)\lceil \tfrac{c}{5} \rceil+4,b}].\]
        Similarly, for every \(y_{m'}\), there is some \(q' \in \set{0,1,2,3}\) such that the indexes of its support are contained in the analogous interval inside \(B_{j'}\).

        By niceness of \(\cI\), for every pair of these intervals there exists \(K \in \cI \setminus \set{J}\) that avoids both intervals. Consequently, the corresponding pair \(x_m,y_{m'}\) lies in the same subgroup
        \[H_K^{[d]\setminus J} = U_{[d]\setminus(J\cup K)}.\]
        Since \(x_m\) and \(y_{m'}\) belong to distinct direct factors, they commute in \(U_{[d]\setminus J}\). We can therefore transform \(x_m y_{m'}\) into \(y_{m'}x_m\) using two multiplication-table relations of \(H_K^{[d]\setminus J}\): first replace \(x_m y_{m'}\) by their product, and then expand this product as \(y_{m'}x_m\).

        Commuting all pairs \(x_m,y_{m'}\), and then recombining the factors, transforms \(xy\) into \(yx\) using \(O(C_0 C_1)=O(1)\) relations from \(R\).
    \end{proof}
\end{proof}

\section{Bounding the diameter and the Dehn function} \label{sec: bounding the Dehn function}

The aim of this section is to prove \pref{lem: diameter Y} and \pref{lem: Yl Dehn function} and thus conculding the proof of \pref{thm: dehn and diameter bound for link}.

Let \(\cI\) is a nice set of indices, let \(J \in \cI\) and let \[Y \coloneqq \CC\left(U_{[d]\setminus J};\left(H_I^{[d]\setminus J}\right)_{I \in \cI \setminus \set{J}}\right)\] and for $0 \leq \ell \leq r$ 
\[Y_\ell \coloneqq \CC\left(U_{B_\ell};\left(H_I^{B_\ell}\right)_{I \in \cI \setminus \set{J}}\right).\]

In light of \pref{thm: dehn implies cone}, in order to bound the cone radii of $Y$, it is enough to bound the diameter of $Y$ and the Dehn function of  $U_{[d]\setminus J}$.

We will start by bounding the diameter of $Y$ (proving \pref{lem: diameter Y}):
\begin{proposition} \label{prop: bounding diameter}
It holds that $\diam (Y)\leq 6(r+1) +1$. 
\end{proposition}

\begin{proof}
By \cite[Proposition 6.3]{KaufmanO2021}, it is enough to show that every $g \in U_{[d]\setminus J}$ can be written as a product of at most $6(r+1)$ elements in $\bigcup_{I \in \cI \setminus \set{J}} H_I^{[d]\setminus J}$. Since $U_{[d]\setminus J}$ is a direct product of the groups $U_{B_0},...,U_{B_r}$ it follows that there are $g_\ell \in U_{B_{\ell}}, 0 \leq \ell \leq r$ such that $g = g_0 ... g_r$. Thus it is enough to show that for every $0 \leq \ell \leq r$ and every $g_{\ell} \in U_{B_{\ell}}$, $g_{\ell}$ can be written as at most $6$ elements of $\bigcup_{I \in \cI \setminus \set{J}} H_I^{B_{\ell}}$. 

Fix $0 \leq \ell \leq r$. By the assumption that \(\cI\) is nice, it holds that there are $I_1,I_2 \in \cI$ such that $I_1 \cap B_{\ell} \subseteq [a, a+ \lfloor \frac{b-a}{5} \rfloor]$, $I_{2} \cap B_{\ell} \subseteq [b- \lfloor \frac{b-a}{5} \rfloor, b]$. Thus by \pref{prop:diameter-of-large-triangle-matrix-groups}, each element in $U_{B_{\ell}}$ can be written as a product of at most $6$ elements in $H_{I_1}^{B_{\ell}} \cup H_{I_2}^{B_{\ell}}$ as needed.
\end{proof}

In light of \pref{lem:direct-product-of-groups}, in order to prove a bound on the Dehn function of $U_{[d]\setminus J}$, it is enough to prove a bound on the Dehn function of $\Dehn_{Y_\ell}$ for every $0 \leq \ell \leq r$, where we recall that 
 \(R_I^\ell\) is the multiplication-table relations of \(H_I^{B_\ell}\) and
\[\Dehn_{Y_\ell} = \Dehn_{S_\ell,R_\ell}, \qquad S_\ell = \bigcup_{I \in \cI \setminus \set{J}} H_I^{B_\ell}, \qquad R_\ell = \bigcup_{I \in \cI \setminus \set{J}} R_I^\ell.\]

Thus, we can restate and prove \pref{lem: Yl Dehn function}:

\begin{lemma}
\label{lem:dehn-bound-on-bin}
There is a constant $A$ such that for every $0 \leq \ell \leq r$ and every $s \in \mathbb{N}$,
$\Dehn_{S_\ell,R_\ell} (s) \leq A s$.
\end{lemma}

\begin{proof}
Fix $0 \leq \ell \leq r$ and denote $B_{\ell} = [a,b], n = b-a$. In the notation of \pref{sec: Unip sec}, $U_{B_{\ell}}$ is isomorphic to $\Unip_n (\mathbb{F}_q)$. 

\begin{comment}
We define $I_{\min}, I_{\max} \in \mathcal{I} \setminus \{ J \}$ to be the sets such that
$$\max (I_{\min} \cap B_{\ell}) \leq \min_{I \in \mathcal{I} \setminus \{ J \}} \max (I \cap B_{\ell}),$$
$$\min (I_{\max} \cap B_{\ell}) \geq \max_{I \in \mathcal{I} \setminus \{ J \}} \min (I \cap B_{\ell}).$$ 
\ydnote{reverse inequality?}
\end{comment}

By the assumption that \(\cI\) is nice, it holds that there are $I_1,I_2, I_3 \in \cI$ such that $I_1 \cap B_{\ell} \subseteq [a, a+ \lfloor \frac{b-a}{5} \rfloor]$, $I_{2} \cap B_{\ell} \subseteq [b- \lfloor \frac{b-a}{5} \rfloor, b]$ and $I_{3} \cap B_{\ell} \subseteq [a+  \lceil \frac{b-a}{5} \rceil, a+ 3\lfloor \frac{b-a}{5} \rfloor]$.

Thus, in the notion of \pref{sec: Unip sec}, $S_{\Elm} \subseteq S_{\ell}$ and every relation in $R_{\Elm}$ can be derived from $R_{I_1}^{\ell} \cup R_{I_2}^{\ell} \cup R_{I_3}^{\ell} \subseteq R_{\ell}$ by applying at most $4$ relations in $R_{I_1}^{\ell} \cup R_{I_2}^{\ell} \cup R_{I_3}^{\ell} \subseteq R_{\ell}$. 

By \pref{thm:unip-presentation}, there is a universal constant $C'$ such that $\Dehn_{S_{\Elm}, R_{\Elm}} (m) \leq C' m$. 

Thus it is enough to show that there are universal constants $C_0, C_1$ such that given a word $g \in S_{\ell}$, $g$ can be transformed into $w_1 ... w_m$ where $m \leq C_0, w_1,...,w_m \in S_{\Elm}$ and this requires $C_1$ relations in $R_{\ell}$. Indeed, assume we showed this for every $g \in S_{\ell}$. Then given a trivial word $g_1 ... g_s$ in the alphabet $S_{\ell}$, it can be transformed into a word $w_1...w_t$ in the alphabet $S_{\Elm}$ with $t \leq C_0 s$ using $C_1 s$ relations in $R_{\ell}$. Then using \pref{thm:unip-presentation}, $w_1...w_t$ can be reduced to the trivial word using at most $C' m$ relations in $R_{\Elm}$ and thus by at most $4 C' m$ relations in  $R_{\ell}$. It follows that 
$$\Dehn_{S_{\ell}, R_{\ell}} (s) \leq (4C_0 C' + C_1)s,$$
as needed.

Let $g \in S_{\ell}$. By \pref{claim:diameter-bound-for-good-groups}, there are constants $C_0 ', C_1'$ such that $g$ can be transformed into $w_1 ' ... w_{m'}'$ where $m' \leq C_0 '$, and $w_i '$ are of block size $\leq \frac{n}{5}$. Thus applying at most $2$ relations from $R_\ell$, each $w_i '$ can be transformed to a product of at most $3$ elements in $S_{\Elm}$. It follows that we can take $C_0 = 3 C_0'$ and $C_1 = C_1 ' + 2C_0'$ (all are universal constants). 
\end{proof}

With the above Lemma, we can finally restate and prove \pref{thm: dehn and diameter bound for link}:

\begin{theorem}
With respect to the coset-complex presentation of \(Y\),
    \[\Dehn_Y(t) \leq O((t(r+1))^2) \]
    for every \(t \in \N\).
\end{theorem}

\begin{proof}
Combine \pref{lem:direct-product-of-groups} and \pref{lem:dehn-bound-on-bin}.
\end{proof}

\section{AI Methodology} \label{sec:ai}
All ideas and proofs are independent work of the authors. We used Open AI's Astra 6 to simplify arguments and optimize constants in the proof of \pref{lem:most-colors-are-nice} as our original proof had a bound of \(r \geq 10^{2700}\). We also used various other LLMs to proof read, edit text, create diagrams and search the literature for references. The authors verified the correctness and originality of all content, including references.

\printbibliography

\appendix
\section{PCPs with improved proof-length dependence}
\label{app:pcp}

We show that substituting \pref{thm:main} into the construction of Bafna, Minzer, Vyas and Yun improves its dependence on the soundness parameter. We follow their construction and only record the parameter changes. 

% \subsection{Preliminaries for this appendix}
% \ydnote{TODO: preliminaries of what a PCP is, what a label cover is. The agreement test.}

\subsection[Adapting Theorem 1.2 in Bafna et.\ al.\ ]{Adapting {\cite[Theorem 1.2]{BafnaMVY2025}}}
Following \cite{BafnaMVY2025} we prove the following theorem.
\restatetheorem{thm:pcp}

% \ydnote{The alphabet bound stated here is the bound obtained directly by converting the final robust PCP into label cover. \cite[Theorem 8.9]{BafnaMVY2025}, states the stronger bound \(\poly(1/\varepsilon)\); check the additional accounting needed for that bound before claiming it here.}

We first establish the agreement theorem and the quantitative complex construction, and then apply the embedding and composition arguments of \cite{BafnaMVY2025}.

\subsubsection{Complexes for the PCP construction}
\label{subsec:pcp-complexes}

In this subsection we give a quantitative variant of \cite[Theorem 2.13]{BafnaMVY2025}. The improvement comes from substituting \pref{thm:main} into the agreement-testing argument of \cite{DiksteinD2023agr}. Throughout this section, all polynomial dependencies have absolute constant exponents, independent of the soundness parameter. %We use our convention that \(X(k)\) consists of faces of size \(k+1\).

\paragraph{Agreement from swap coboundary expansion}
The \(V\)-test on an ensemble \(F=\set{F_s:s\to\Sigma}_{s\in X(k)}\) samples two \(k\)-faces \(s,s'\) whose intersection has size \(\lfloor\sqrt{k+1}\rfloor\), according to the natural distribution induced by \(X\), and accepts if \(F_s|_{s\cap s'}=F_{s'}|_{s\cap s'}\). Equivalently, we first sample their union according to the natural face distribution and then sample the pair uniformly subject to the prescribed intersection size.

We use the following straightforward consequence of the proof of \cite[Theorem 3.1]{DiksteinD2023agr} (a similar theorem with slightly different guarantees was proven in \cite{BafnaM2023}). Here and below, taking a parameter to be a sufficiently large polynomial means that the polynomial is fixed independently of the complex and the alphabet.

\begin{theorem}[Agreement from swap coboundary expansion]
    \label{thm:agreement-from-swap}
    Let \(\varepsilon\in(0,1/2)\), let \(k\geq\exp(\poly(1/\varepsilon))\), and let \(d\) be a sufficiently large polynomial in \(k\). Let \(X\) be a \(d\)-dimensional clique complex satisfying the following:
    \begin{enumerate}
        \item \(X\) is a \(1/d^2\)-one-sided local spectral expander.
        \item For a large intermediate dimension \(r=\poly(k)\), with \(d\)  polynomially in \(r\), the faces complex \(\mathcal{F}^rX\) is a \((\poly(1/r),m)\)-coboundary expander for \(m=\poly(1/\varepsilon)\).
    \end{enumerate}
    
    Let \(\Sigma\) be a finite alphabet and let \(F=\set{F_s:s\to\Sigma}_{s\in X(k)}\). If the \(V\)-test accepts \(F\) with probability at least \(\varepsilon\), then there exists \(f:X(0)\to\Sigma\) such that
    \[\Prob[s\in X(k)]{\dist(F_s,f|_s)\leq\gamma}\geq\varepsilon^{O(1)},\]
    where \(\gamma=\exp(\poly(1/\varepsilon))\poly(1/k)\).
\end{theorem}

%The displayed dependence of \(\gamma\) uses the refinement discussed here. The statement in \cite{DiksteinD2023agr} gives an \(\exp(\poly(1/\varepsilon))\) loss.

\begin{remark}
    Note that \cite{DiksteinD2023agr} additionally requires a property called well connectedness, which is used to descend a cover of the faces complex to a cover of \(X\). This requirement can be removed under our stronger assumption of coboundary expansion: the corrected cochain is a coboundary, so we can synchronize the lists directly without constructing or descending a cover. See also \cite{BafnaM2023}.    
\end{remark}

\begin{proof}[Proof sketch]
    We explain the quantitative improvement over the argument of \cite[Section 3]{DiksteinD2023agr}. In that proof, a cochain representing the agreement problem is constructed on a sufficiently large faces complex of \(X\). Coboundary expansion is then used to correct this cochain to a nearby coboundary, from which a global assignment is constructed. We refer the reader to \cite[Section 3]{DiksteinD2023agr} for the details of this argument.

    The fraction of triangles violated by the constructed cochain is at most \(\mu = \exp(-\poly(\varepsilon)r/k) + O(r/d)\). In \cite{DiksteinD2023agr}, the available coboundary expansion was \(\exp(-O(\sqrt{r}))\), and therefore the argument required \(d \gg r\exp(O(\sqrt{r}))\) to ensure that \(\mu\) is sufficiently small for the correction to be close to the original cochain. Since \(r > \poly(\varepsilon^{-1},k)\), this resulted in a superpolynomial dependence of \(d\) on \(k\). Using our improved coboundary expansion of \(1/\poly(r)\), it is enough to take \(d = \poly(r) = \poly(k)\), while the rest of the argument remains unchanged.

    The bound \(d=\poly(k)\) follows by tracking the dimension requirements throughout the argument of \cite[Section 3]{DiksteinD2023agr}; we omit the full parameter bookkeeping for brevity.
\end{proof}

\paragraph{The quantitative complex construction}
We now give the variant of \cite[Theorem 2.13]{BafnaMVY2025} used for the PCP construction. For a target size \(n\), we construct a complex with between \(n\) and \(nq^{O(d^2)}\) vertices, where \(q\) is the field size. Since \(q\) is polylogarithmic in \(n\), this additional factor is absorbed into the claimed proof-length bound.

We will use the following fact about the links of KMS complexes. We sketch its proof after the proof of \pref{thm:pcp-complexes}.

\begin{fact}\label{fact:pcp-link-expansion-diameter}
    There is an absolute constant \(C_0>0\) such that the following holds. Let \(X=X(A_d,\F_q,f)\), where \(f\) is irreducible of degree at least two and \(q\geq C_0(d+1)\).
    \begin{enumerate}
        \item For every face link \(L\), including \(X\) itself, the bipartite graph between any two distinct colors is uniformly weighted over its edges and has second singular value at most \(4/\sqrt q\).
        \item For every vertex link \(L\), order its colors as \(1,\dots,d\) according to the subspace dimensions in its opposition-complex description. For \(1\leq a<b\leq d\), the bipartite graph between colors \(a,b\) has diameter \(O(d/(b-a))\), with an absolute implied constant. In particular, this diameter is constant when \(b-a=\Omega(d)\).
    \end{enumerate}
\end{fact}

%We now give the variant of \cite[Theorem 2.13]{BafnaMVY2025} used for the PCP construction. We allow a size bound of \(nq^{O(d^2)}\) (here \(n\) is the number of vertices and \(q\), the fieldsize, is a polynomial in the degree of a vertex), rather than a constant multiple of \(n\) independent of \(q\). This is sufficient because \(q\) is polylogarithmic in \(n\).

\begin{theorem}
    \label{thm:pcp-complexes}
    For every \(\varepsilon\in(0,1/2)\), there are parameters
    \[k=\exp(\poly(1/\varepsilon)),\qquad d=\poly(k),\qquad q_0=\exp(\poly(1/\varepsilon))\]
    with the following property. For every fixed \(C>0\), all sufficiently large \(n\), and every prime \(q=\Theta((\log n)^C)\) with \(q\geq q_0\), one can construct in polynomial time a \(d\)-dimensional KMS complex \(X\) such that
    \[n\leq |X(0)|\leq nq^{O(d^2)}\]
    and: \begin{enumerate}
    \item \(X\) is a \((d+1)\)-partite clique complex.
    \item Every vertex belongs to at most \(q^{O(d^2)}\) top-dimensional faces.
    \item For every vertex link \(L\), there is a group of automorphisms acting transitively on \(L(d-1)\).
        \item For every vertex link \(L\), order its colors as \(1,\dots,d\) according to the subspace dimensions in its opposition-complex description. For distinct colors \(i,j\) in this ordering, the bipartite graph between \(L[i]\) and \(L[j]\) is uniformly weighted over its edges and has diameter \(O(d/|i-j|)\).
        \item For every face link \(L\), including \(X\) itself, the bipartite graph between any two distinct colors has second singular value at most \(4/\sqrt q\).
        \item The normalized adjacency operator of the graph \((X(0),X(1))\) has second singular value at most \[\frac{1}{d}+\frac{4}{\sqrt q}\leq\frac{2}{d}.\]
        \item Let \(\Sigma\) be a finite alphabet and let \(F=\set{F_s:s\to\Sigma}_{s\in X(k)}\). Consider the following test:
        \begin{enumerate}
            \item Sample \(t\in X(d)\).
            \item Sample a uniform subset \(b\subseteq t\) of size \(\lfloor\sqrt{k+1}\rfloor\).
            \item Independently sample two uniform \(k\)-faces \(s,s'\subseteq t\) containing \(b\).
            \item Accept if \(F_s|_b=F_{s'}|_b\).
        \end{enumerate}
        If this test accepts with probability at least \(\varepsilon\), then there exists \(f:X(0)\to\Sigma\) such that
        \[\Prob[s\in X(k)]{\dist(F_s,f|_s)\leq\varepsilon}\geq\varepsilon^{O(1)}.\]
    \end{enumerate}
    All constants in the displayed parameter bounds are independent of \(\varepsilon\).
\end{theorem}

\begin{proof}[Proof sketch]
    Choose \(k=\exp(\poly(1/\varepsilon))\), \(d=\poly(k)\), and \(q_0=\exp(\poly(1/\varepsilon))\) sufficiently large that \pref{thm:agreement-from-swap} applies to the family of KMS complexes for every prime \(q\geq q_0\), with acceptance threshold \(\varepsilon/2\) and error \(\gamma\leq\varepsilon\). We also take \(q_0\) large enough to satisfy the hypotheses of \pref{fact:pcp-link-expansion-diameter} and the numerical bound in the sixth item.

    The first three items and the sixth follow from basic properties of the KMS complexes and their links \cite{GraveDePeraltaV2025,DiksteinHIP2026}. For the clique-complex property, we take \(\deg(f)\geq4\) and use the opposition-complex description in \cite[Section 3.2]{GraveDePeraltaV2025}.\footnote{Before taking the congruence quotient, the complex is a full opposition subcomplex of a building and hence a clique complex. Products of three local subgroup matrices have entries of degree at most three, so reduction modulo \(f\) creates no new triangles when \(\deg(f)\geq4\). Every triangle in the quotient therefore spans a face. Since its vertex links are clique complexes, the same holds for every larger clique.} The fourth and fifth items follow from \pref{fact:pcp-link-expansion-diameter}, which we prove below. We omit the other standard verifications.

    For size selection, vary \(s=\deg(f)\geq4\). Each vertex stabilizer has order \(q^{\binom{d+1}{2}}\), so the number of vertices is \[N_s=(d+1)\frac{|\Sl_{d+1}(\F_{q^s})|}{q^{\binom{d+1}{2}}}.\] The standard order formula for \(\Sl_{d+1}\) gives \(N_{s+1}/N_s\leq2q^{(d+1)^2-1}=q^{O(d^2)}\). Choose the least \(s\geq4\) for which \(N_s\geq n\). For sufficiently large \(n\), we have \(N_4<n\), since \(d\) is fixed and \(q\) is polylogarithmic in \(n\). Minimality then gives \(n\leq N_s\leq nq^{O(d^2)}\). The resulting complex can be constructed in polynomial time by \cite{DiksteinHIP2026}.

    Note that the test in the last item differs slightly from the \(V\)-test, since the sampled faces may intersect outside \(b\). However, this occurs with probability at most \(O(k^2/d)\). Taking \(d\) sufficiently large ensures that acceptance probability \(\varepsilon\) in our test implies acceptance probability at least \(\varepsilon/2\) in the \(V\)-test, so the agreement theorem applies.
\end{proof}

\begin{proof}[Proof sketch of \pref{fact:pcp-link-expansion-diameter}]
    By \cite{GraveDePeraltaV2025,DiksteinHIP2026}, the two-color graphs in vertex links are symmetric induced subgraphs of the corresponding subspace-incidence graphs, retaining a \(1-O(1/q)\) fraction of the vertices and neighbors. Symmetry gives uniform edge weights. The full incidence graphs have second singular value at most \(1/\sqrt q\) and diameter \(O(d/(b-a))\). For sufficiently large \(q\) relative to \(d\), the standard comparisons for these dense symmetric restrictions increase both bounds by at most a factor of two.

    The spectral bound also holds in every proper link, since these are joins of smaller opposition complexes. Finally, partite trickle-down \cite[Lemma 2.10]{BafnaMVY2025} gives the bound \(4/\sqrt q\) for \(X\) itself.
\end{proof}

For these parameters, the number of top-dimensional faces is at most
\[|X(d)|\leq |X(0)|q^{O(d^2)}\leq nq^{O(d^2)}.\]
In particular, when \(q=\Theta((\log n)^C)\) for fixed \(C\), this is at most \(n(\log n)^{\exp(\poly(1/\varepsilon))}\), for sufficiently large \(n\). This gives the quantitative complex input for the \cite{BafnaMVY2025} construction; the remaining parameter calculation follows their embedding, soundness amplification and alphabet reduction.

\subsection{The PCP construction and its parameters}

\begin{proof}[Proof of \pref{thm:pcp}]
    We follow the construction of \cite[Section 8]{BafnaMVY2025}, replacing their Theorem 2.13 by \pref{thm:pcp-complexes}. We explain the resulting parameter improvement and the minor adjustment required in the embedding step.

    As in \cite[Lemma 8.8 and Theorem 8.9]{BafnaMVY2025}, it suffices to start with an agreement test of soundness \(\delta^2\), where \(\delta=\varepsilon^{O(1)}\). Applying \pref{thm:pcp-complexes}, we obtain a family of complexes with parameters \(k=\exp(\poly(1/\varepsilon))\), \(d=\poly(k)\), and \(q=\Theta((\log n)^C)\), for a sufficiently large absolute constant \(C>0\). In particular, we may assume that \(q\) is larger than the threshold required by \pref{thm:pcp-complexes} by taking \(n\) sufficiently large.

    We apply the embedding construction of \cite[Lemma 8.6]{BafnaMVY2025}, starting with a constant-gap CSP. To accommodate the size of the chosen complex, we take sufficiently many disjoint copies of the initial CSP so that its number of vertices matches that of the auxiliary routing graph up to fewer than one copy, as in \cite[Lemma 5.2]{BafnaMVY2025}. This preserves a constant soundness gap. In each vertex link, we use the color ordering from \pref{fact:pcp-link-expansion-diameter} when applying the local routing argument; the remaining robust routing argument is unchanged. We then apply \cite[Lemma 6.2]{BafnaMVY2025} to obtain a generalized label-cover instance with perfect completeness and soundness at most \(\delta\).

    The quantitative difference is in the size of this instance. By \pref{thm:pcp-complexes}, the chosen complex has at most \(n(\log n)^{\exp(\poly(1/\varepsilon))}\) vertices, and every vertex belongs to at most \(q^{O(d^2)}\) top-dimensional faces. Consequently, the number of top-dimensional faces, and hence the size of the resulting generalized label-cover instance, is at most \(n(\log n)^E\), where \(E=\exp(O(\varepsilon^{-a}))\) for a sufficiently large absolute constant \(a>0\). Its logarithmic alphabet sizes and decision complexity are also bounded by \((\log n)^E\). Passing to a robust PCP using \cite[Lemma 7.3]{BafnaMVY2025}, as in \cite[Lemma 8.8]{BafnaMVY2025}, gives proof length at most \(n(\log n)^E\), randomness complexity at most \(\log_2 n+E\log_2\log n\), and decision complexity at most \((\log n)^E\), after increasing the absolute constants in \(E\).

    The remaining construction is exactly that of \cite[Lemma 8.8 and Theorem 8.9]{BafnaMVY2025}: two compositions with the Reed--Muller based decodable PCP, followed by the Hadamard based composition and the conversion back to label cover. Their parameter analysis shows that these operations preserve the proof-length bound \(n(\log n)^{O(E)}\), while achieving perfect completeness and soundness at most \(\varepsilon\). The final robust PCP has alphabet size \((1/\varepsilon)^{O(1)}\) and query complexity \((1/\varepsilon)^{O(\log(1/\varepsilon))}\). Thus the conversion to a two-query PCP gives a proof alphabet of size at most \(\exp((1/\varepsilon)^{O(\log(1/\varepsilon))})\), as required. Increasing the absolute constant in the definition of \(E\) completes the proof.
\end{proof}
\section{Equivalence between cones and bounding the Dehn function} \label{app:dehn-equivalence}

In this subsection we establish a quantitative connection between the geometry of a coset complex and the presentation of its underlying group. More precisely, we show that the diameter of a coset complex is comparable to the word diameter of the group, and that filling loops in the complex is equivalent, up to a linear distortion, to reducing words using the multiplication tables of its subgroups. In particular, in the end of this section we prove \pref{thm: dehn implies cone}.

The implication from the Dehn function to fillings of the coset complex was established in a different quantitative form in \cite[Theorem 6.6]{KaufmanO2021}. We give a proof of both directions. Our proof is based on a general comparison between a simplicial complex and the nerve of a cover by subcomplexes.

For a loop $P$ in a simplicial complex $Z$, let $\Area_Z(P)$ be the minimum number of $\TR$ moves needed to contract $P$, where $\BT$ moves are free, as in \pref{subsec:dehn}. We set $\Area_Z(P)=\infty$ if $P$ is not contractible, and write
\[
\Area_Z(t)=\sup\set{\Area_Z(P):P\text{ is a loop of length at most }t}.
\]
In other words, we count the triangles used in a filling, with backtracking canceled freely.

\begin{definition}[Nerve]
    Let $\mathcal Y$ be an indexed family of nonempty subcomplexes of $Z$ that covers every vertex and every edge. Its \emph{nerve} $\mathcal A$ has one vertex for each member of $\mathcal Y$, and $\set{Y_0,\dots,Y_j}$ is a face if these subcomplexes have a common vertex. We retain the indices: two members of the family may be equal as subcomplexes and still give distinct vertices of $\mathcal A$.
\end{definition}

The application we have in mind is a coset complex $Y=\CC(G;(H_i)_{i\in I})$. Define an auxiliary simplicial complex $Z$ whose vertices are the elements of $G$, and whose nonempty faces are the finite subsets contained in a single coset $gH_i$. Its $1$-skeleton is the Cayley graph of $G$ with generating set $\bigcup_{i\in I}H_i\setminus\set{1}$. Note that this is not necessarily a clique complex.

\ydnote{write down as a formula:}
The simplices on the cosets $gH_i$, indexed by both the coset and its color, cover $Z$, and their nerve is $Y$. Conversely, cover $Y$ by the chambers
\[
Y_g=\set{gH_i:i\in I},\qquad g\in G,
\]
viewed as simplices and indexed by $g$. Their nerve is $Z$: the chambers $Y_{g_0},\dots,Y_{g_j}$ have a common vertex if and only if $g_0,\dots,g_j$ belong to a common coset of one of the $H_i$. In both covers, every member and every nonempty intersection is a simplex.

We first establish two general lemmas about nerves, and then apply them to these covers.

\begin{lemma} \label{lem:nerve-diameter}
Let $Z$ be a connected simplicial complex, and let $\mathcal Y$ be an indexed cover of its vertices and edges by nonempty connected subcomplexes of diameter at most $D_1$, where $D_1\geq 1$. If $\mathcal A$ is the nerve, then $\mathcal A$ is connected and
\[
\diam(\mathcal A)\leq\diam(Z)+1,\qquad
\diam(Z)\leq D_1(\diam(\mathcal A)+1).
\]
\end{lemma}

\begin{proof}
Let $Y,Y'\in\mathcal Y$, choose $u\in Y$ and $v\in Y'$, and take a shortest path from $u$ to $v$ in $Z$. Choose a member of $\mathcal Y$ containing each edge of the path. Consecutive members intersect, the first intersects $Y$, and the last intersects $Y'$. Thus $Y,Y'$ are joined in $\mathcal A$ by a path of length at most $\dist_Z(u,v)+1$. If the chosen path has length zero, $Y$ and $Y'$ already intersect. This proves connectivity and the first bound.

Conversely, let $(Y_0;\dots;Y_m)$ be a path in $\mathcal A$, and choose any vertices $u\in Y_0$ and $v\in Y_m$. Choose one vertex in every consecutive intersection $Y_{j-1}\cap Y_j$. Connecting these vertices, together with $u$ and $v$, by paths inside the respective $Y_j$ gives a path in $Z$ of length at most $(m+1)D_1$. Taking $m\leq\diam(\mathcal A)$ proves the second bound.
\end{proof}

\begin{lemma} \label{lem:nerve-fillings}
Let $Z$ be connected, and let $\mathcal Y$ be an indexed cover of its vertices and edges by nonempty connected subcomplexes. Suppose that every $Y\in\mathcal Y$ has diameter at most $D_1$, and that every nonempty intersection $Y\cap Y'$ is connected and has diameter at most $D_2$, where $D_1,D_2\geq 1$. Suppose further that every loop of length at most $2D_1+2D_2+1$ in any $Y\in\mathcal Y$ can be contracted inside $Y$ using at most $R$ triangle moves. If $\mathcal A$ is the nerve, then for every $t\geq 1$,
\[
\Area_Z(t)\leq 3R\bigl(\Area_{\mathcal A}(t)+t\bigr).
\]
\end{lemma}

\usetikzlibrary{fit,calc,shapes.geometric}
% requires: \usetikzlibrary{fit,calc,shapes.geometric}

\begin{proof}
For each $Y\in\mathcal Y$, choose a vertex $v_Y$ and, for every $u\in Y$, a path $P_{Y,u}$ from $v_Y$ to $u$ inside $Y$ of length at most $D_1$. For every edge $\set{Y,Y'}$ of $\mathcal A$, choose $u_{Y,Y'}=u_{Y',Y}\in Y\cap Y'$. We lift the oriented edge $(Y,Y')$ to the path
\[
L_{Y,Y'}=P_{Y,u_{Y,Y'}}\circ P_{Y',u_{Y,Y'}}^{-1}
\]
from $v_Y$ to $v_{Y'}$. In particular, $L_{Y',Y}=L_{Y,Y'}^{-1}$. For a path $P'=(Y_0;Y_1;\dots;Y_m)$ in $\mathcal A$, its \emph{canonical lift} is the concatenation
\[
L(P')=L_{Y_0,Y_1}\circ\cdots\circ L_{Y_{m-1},Y_m}.
\]
When $P'$ is a loop, so is $L(P')$.

We will prove that every loop $P$ of length $n$ in $Z$ can first be transformed, using at most $3Rn$ triangle moves, into a canonical lift of a loop $P'$ in $\mathcal A$ of length at most $n$, up to conjugation by a path. We will then show that each elementary move in a contraction of $P'$ can be performed on its canonical lift using at most $3R$ triangle moves. Importantly, we only count triangle moves: the lengths of intermediate paths need not remain bounded by $n$.

\smallskip
\noindent\textbf{Step 1: Canonicalizing loops.}
Let $P=(u_0;u_1;\dots;u_n=u_0)$ be a loop in $Z$. For each edge $(u_i,u_{i+1})$, choose $Y_i\in\mathcal Y$ containing it, with indices interpreted cyclically. Since $u_i\in Y_{i-1}\cap Y_i$, the cyclic sequence $(Y_0;Y_1;\dots;Y_{n-1};Y_0)$ is a closed walk in $\mathcal A$, possibly with consecutive repetitions. Removing these repetitions gives a loop $P'$ of length at most $n$.

We first replace each edge $(u_i,u_{i+1})$ of $P$ by
\[
P_{Y_i,u_i}^{-1}\circ P_{Y_i,u_{i+1}}.
\]
The original edge and its replacement form a loop inside $Y_i$ of length at most $2D_1+1$. By assumption, this replacement can be performed using at most $R$ triangle moves, so the total cost is at most $Rn$. For an illustration of this part see \pref{fig:canonicalization-1}.

After these replacements, the path travels through the base vertices of the sets $Y_i$. More precisely, define
\[
T_i=P_{Y_{i-1},u_i}\circ P_{Y_i,u_i}^{-1},
\]
which is a path from $v_{Y_{i-1}}$ to $v_{Y_i}$. Up to backtracking, the resulting loop is
\[
C\circ T_1\circ T_2\circ\cdots\circ T_{n-1}\circ T_0\circ C^{-1},
\qquad C=P_{Y_0,u_0}^{-1}.
\]
Thus, it remains to replace each transition $T_i$ by the canonical path $L_{Y_{i-1},Y_i}$.

If $Y_{i-1}=Y_i$, then $T_i$ cancels by backtracking. Otherwise, let $w=u_{Y_{i-1},Y_i}$. Since $u_i,w\in Y_{i-1}\cap Y_i$, we may choose a path $Q$ from $u_i$ to $w$ inside this intersection, of length at most $D_2$. Inserting the backtracking $Q\circ Q^{-1}$, we write
\[
T_i=P_{Y_{i-1},u_i}\circ Q\circ Q^{-1}\circ P_{Y_i,u_i}^{-1}.
\]
We now perform two contractions. First, inside $Y_{i-1}$, replace $P_{Y_{i-1},u_i}\circ Q$ by $P_{Y_{i-1},w}$. Second, inside $Y_i$, replace $Q^{-1}\circ P_{Y_i,u_i}^{-1}$ by $P_{Y_i,w}^{-1}$. Each replacement is justified by a loop of length at most $2D_1+D_2$ inside the corresponding covering set, and therefore costs at most $R$ triangle moves. The resulting path is precisely
\[
P_{Y_{i-1},w}\circ P_{Y_i,w}^{-1}=L_{Y_{i-1},Y_i}.
\]

There are at most $n$ transitions, so this second part costs at most $2Rn$ triangle moves. We have therefore transformed $P$ into
\[
C\circ L(P')\circ C^{-1}
\]
using at most $3Rn$ triangle moves. Conjugation by $C$ will not affect the number of triangle moves needed for the remaining contraction. For an illustration of this part see \pref{fig:canonicalization-2}.

\begin{figure}[ht]
\centering
\begin{tikzpicture}[scale=1.05, every node/.style={font=\small}]

% Covering sets
\draw[blue!70!black,dashed,thick] (2,2.2) ellipse (2.5 and 1.75);
\draw[blue!70!black,dashed,thick] (5.85,1.65) ellipse (2.55 and 2.5);
\node[blue!70!black] at (2,4.2) {$Y_0$};
\node[blue!70!black] at (7.3,4.1) {$Y_1$};

% Original loop vertices
\coordinate (a) at (0,3);
\coordinate (b) at (4,3);
\coordinate (c) at (7,0);
\coordinate (d) at (3.5,-2);
\coordinate (e) at (-0.9,-0.2);

% Original closed loop
\draw[thick] (a)--(b)--(c);
\draw[thick,dotted] (c)--(d)--(e);
\draw[thick] (e)--(a);

\fill (a) circle (1.5pt) node[above left] {$a$};
\fill (b) circle (1.5pt) node[above left] {$b$};
\fill (c) circle (1.5pt) node[right] {$c$};
\fill (d) circle (1.5pt);
\fill (e) circle (1.5pt);
\node at (1.5,-2.0) {$\cdots$};

% Replacement path from a to b through u_{Y_0}
\coordinate (a1) at (0.6,2.2);
\coordinate (a2) at (1.25,1.5);
\coordinate (u0) at (2,1.15);
\coordinate (b1) at (2.75,1.55);
\coordinate (b2) at (3.35,2.2);

\draw[blue!75!black,thick]
  (a)--(a1)--(a2)--(u0)--(b1)--(b2)--(b);

\foreach \p in {a1,a2,b1,b2}
  \fill (\p) circle (1.2pt);

\fill[blue!75!black] (u0) circle (1.7pt);
\node[below,blue!75!black] at (u0) {$u_{Y_0}$};

% Replacement path from b to c through u_{Y_1}
\coordinate (c1) at (4.75,3.45);
\coordinate (c2) at (5.6,3.5);
\coordinate (u1) at (6.4,3.05);
\coordinate (c3) at (7.2,2.35);
\coordinate (c4) at (7.7,1.15);

\draw[blue!75!black,thick]
  (b)--(c1)--(c2)--(u1)--(c3)--(c4)--(c);

\foreach \p in {c1,c2,c3,c4}
  \fill (\p) circle (1.2pt);

\fill[blue!75!black] (u1) circle (1.7pt);
\node[above right,blue!75!black] at (u1) {$u_{Y_1}$};

% Double-line arrows: original edge => replacement path
\draw[red,thick,double,->,double distance=1.2pt]
  (2,2.92)--(2,1.48);

\draw[red,thick,double,->,double distance=1.2pt]
  (5.5,1.5)--(6.72,2.72);

\end{tikzpicture}
\caption{Replacing the edges $(a,b)$ and $(b,c)$ by paths through the chosen base vertices of $Y_0$ and $Y_1$, respectively.}
\label{fig:canonicalization-1}
\end{figure}
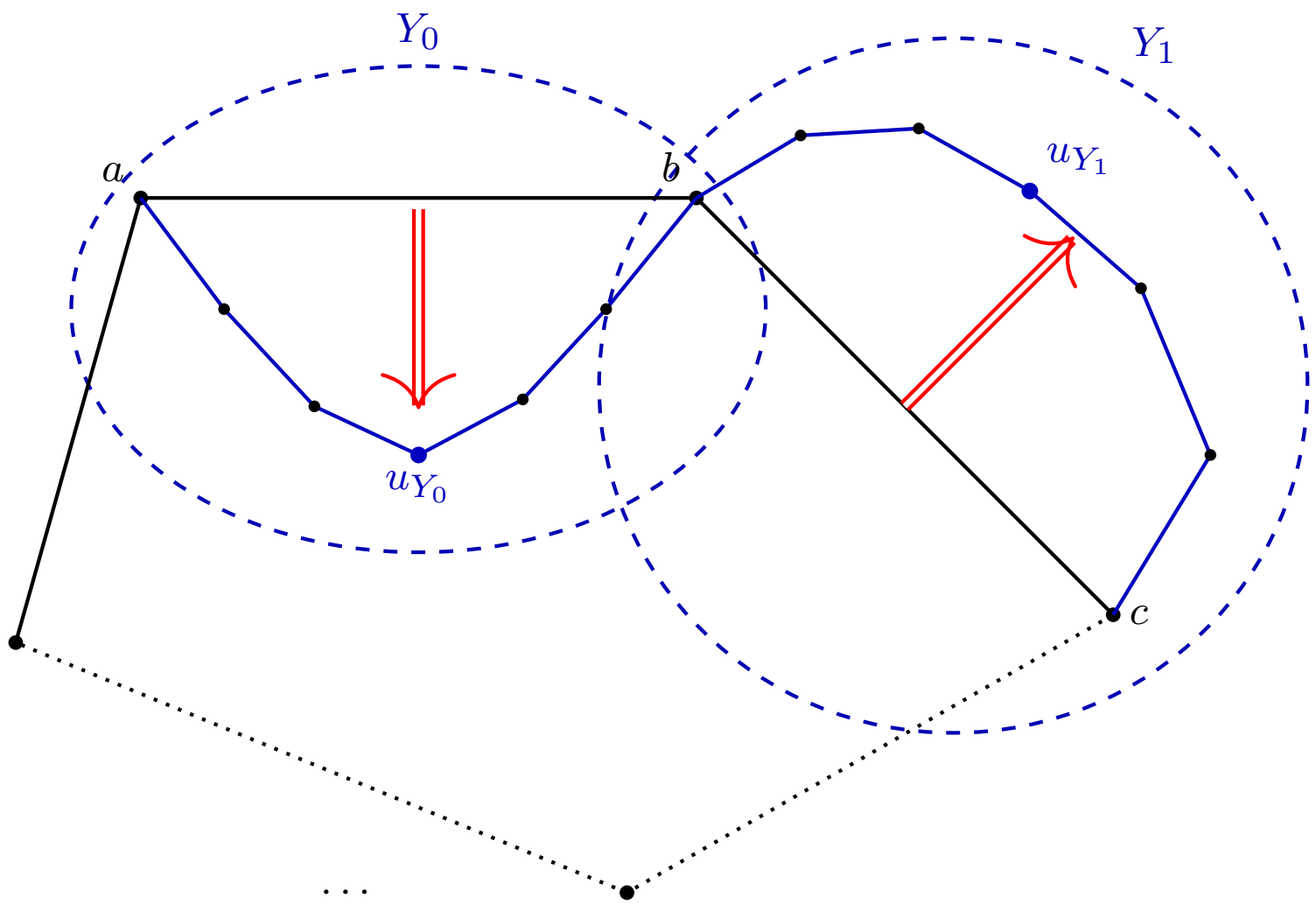

% requires: \usetikzlibrary{calc,shapes.geometric}
% requires: \usetikzlibrary{calc,shapes.geometric}
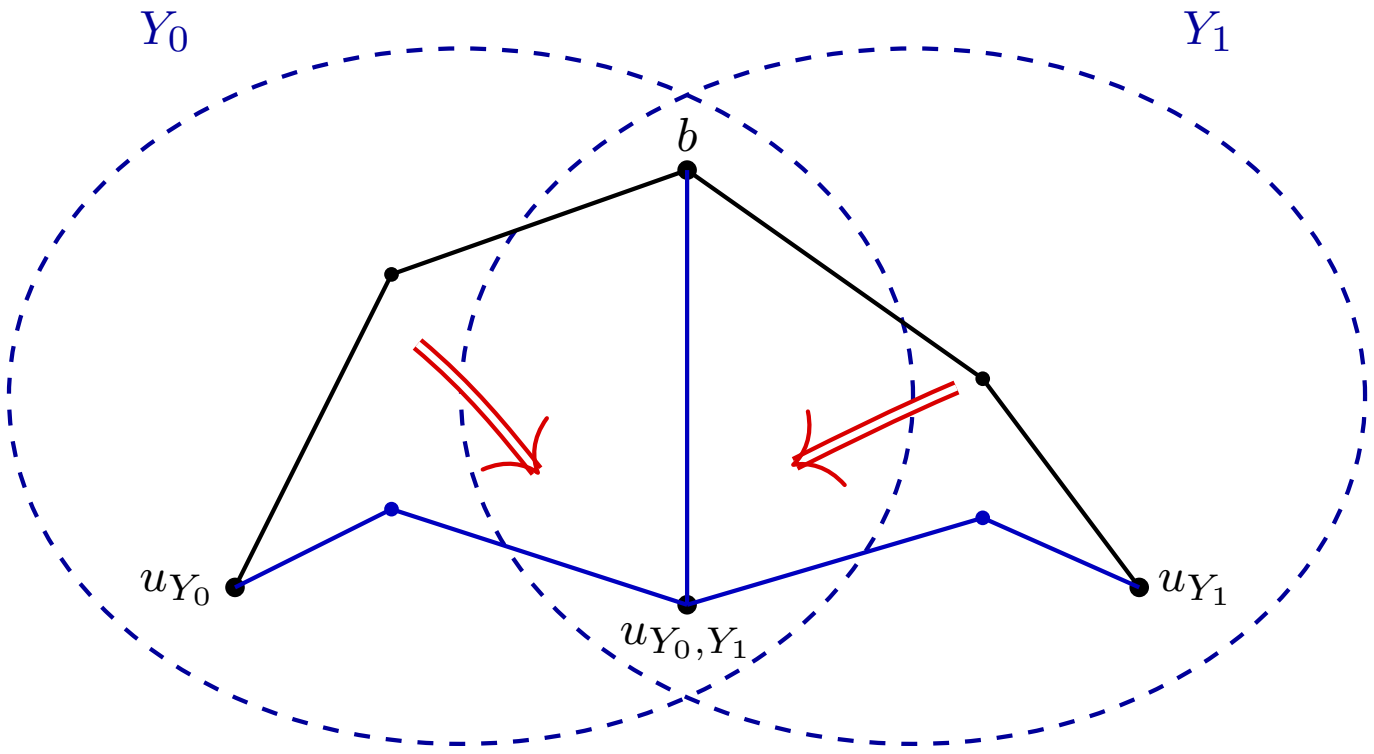
\begin{figure}[ht]
\centering
\begin{tikzpicture}[scale=1.15, every node/.style={font=\small}]

% Colors
\colorlet{oldpath}{black}
\colorlet{newpath}{blue!75!black}
\colorlet{cover}{blue!60!black}
\colorlet{contract}{red!85!black}

% Covering sets
\draw[cover,dashed,thick] (2.7,1.9) ellipse (2.6 and 2.0);
\draw[cover,dashed,thick] (5.3,1.9) ellipse (2.6 and 2.0);
\node[cover] at (1.0,4.0) {$Y_0$};
\node[cover] at (7.0,4.0) {$Y_1$};

% Main vertices
\coordinate (u0) at (1.4,0.8);
\coordinate (b)  at (4.0,3.2);
\coordinate (w)  at (4.0,0.7);
\coordinate (u1) at (6.6,0.8);

\fill (u0) circle (1.6pt);
\fill (b)  circle (1.6pt);
\fill (w)  circle (1.6pt);
\fill (u1) circle (1.6pt);

\node[left]  at (u0) {$u_{Y_0}$};
\node[above] at (b) {$b$};
\node[below] at (w) {$u_{Y_0,Y_1}$};
\node[right] at (u1) {$u_{Y_1}$};

% Old transition path (black)
\coordinate (l1) at (2.3,2.6);
\coordinate (r1) at (5.7,2.0);

\draw[oldpath,thick] (u0)--(l1)--(b)--(r1)--(u1);
\fill (l1) circle (1.2pt);
\fill (r1) circle (1.2pt);

% New canonical pieces (blue dashed)
\coordinate (lb) at (2.3,1.25);
\coordinate (rb) at (5.7,1.2);

\draw[newpath,thick] (u0)--(lb)--(w);
\draw[newpath,thick] (b)--(w);
\draw[newpath,thick] (w)--(rb)--(u1);

\fill[newpath] (lb) circle (1.2pt);
\fill[newpath] (rb) circle (1.2pt);

% Contraction arrows
\draw[contract,thick,double,->,double distance=1.2pt]
  (2.45,2.2) .. controls (2.75,1.95) and (2.95,1.7) .. (3.15,1.45);

\draw[contract,thick,double,->,double distance=1.2pt]
  (5.55,1.95) .. controls (5.2,1.8) and (4.9,1.65) .. (4.6,1.5);

\end{tikzpicture}
\caption{Second part of of the canonicalization. The original transition path from $u_{Y_0}$ to $u_{Y_1}$ passes through $b$. We replace its left part by a path from $u_{Y_0}$ to $b$ through the chosen intersection vertex $u_{Y_0,Y_1}$ inside $Y_0$, and similarly replace its right part inside $Y_1$. The two occurrences of the segment between $b$ and $u_{Y_0,Y_1}$ are then removed by backtracking moves.}
\label{fig:canonicalization-2}
\end{figure}

\smallskip
\noindent\textbf{Step 2: Contracting canonical lifts.}
Consider a contraction of $P'$ in $\mathcal A$ using $a$ triangle moves, together with free backtracking moves. We show that each of these moves can be performed on the corresponding canonical lifts with at most $3R$ triangle moves in $Z$. This applies equally to moves that remove a triangle or backtrack and to their inverses.

A backtracking move inserts or removes two consecutive edges $(Y,Y')$ and $(Y',Y)$. Their canonical lifts are $L_{Y,Y'}$ and $L_{Y',Y}=L_{Y,Y'}^{-1}$, so the corresponding modification in $Z$ consists only of backtracking moves and costs nothing.

Now consider a triangle move that replaces the path $(Y_0;Y_1;Y_2)$ by the edge $(Y_0;Y_2)$, where $\set{Y_0,Y_1,Y_2}$ is a triangle of $\mathcal A$. We need to replace $L_{Y_0,Y_1}\circ L_{Y_1,Y_2}$ by $L_{Y_0,Y_2}$, or equivalently contract
\[
L_{Y_0,Y_1}\circ L_{Y_1,Y_2}\circ L_{Y_2,Y_0}.
\]

For brevity, write $u_{ij}=u_{Y_i,Y_j}$ and $P_{i,u}=P_{Y_i,u}$. Since the three covering sets form a triangle in the nerve, there is a vertex $u\in Y_0\cap Y_1\cap Y_2$. For each unordered pair $\set{i,j}$, choose a path $Q_{ij}$ from $u_{ij}$ to $u$ inside $Y_i\cap Y_j$, of length at most $D_2$.

Expanding the definitions of the canonical lifts, the loop we wish to contract is
\[
P_{0,u_{01}}\circ
\bigl(P_{1,u_{01}}^{-1}\circ P_{1,u_{12}}\bigr)\circ
\bigl(P_{2,u_{12}}^{-1}\circ P_{2,u_{20}}\bigr)\circ
P_{0,u_{20}}^{-1}.
\]
First, inside $Y_1$, replace
\[
P_{1,u_{01}}^{-1}\circ P_{1,u_{12}}
\quad\text{by}\quad
Q_{01}\circ Q_{12}^{-1}.
\]
Both paths connect $u_{01}$ to $u_{12}$ inside $Y_1$, and together form a loop of length at most $2D_1+2D_2$. Thus, the replacement costs at most $R$ triangle moves.

Similarly, inside $Y_2$, replace
\[
P_{2,u_{12}}^{-1}\circ P_{2,u_{20}}
\quad\text{by}\quad
Q_{12}\circ Q_{20}^{-1},
\]
at a cost of at most $R$ triangle moves. After these replacements, our original loop has become
\[
P_{0,u_{01}}\circ Q_{01}\circ Q_{12}^{-1}\circ Q_{12}\circ Q_{20}^{-1}\circ P_{0,u_{20}}^{-1}.
\]
Canceling the backtracking $Q_{12}^{-1}\circ Q_{12}$ leaves
\[
P_{0,u_{01}}\circ Q_{01}\circ Q_{20}^{-1}\circ P_{0,u_{20}}^{-1}.
\]
This is a loop inside $Y_0$ of length at most $2D_1+2D_2$, so it contracts using at most $R$ additional triangle moves. See \pref{fig:example-move} for an illustration.

We have contracted the lifted boundary of the triangle using at most $3R$ triangle moves. Reversing these contractions also allows us to perform the reverse triangle move at the same cost.

Consequently, every elementary move in the contraction of $P'$ can be simulated on its canonical lift, and the entire contraction of $L(P')$ requires at most $3Ra$ triangle moves. This bound does not depend on the lengths of the intermediate paths in the contraction.

Combining the two steps, and taking $a\leq\Area_{\mathcal A}(n)$, gives
\[
\Area_Z(P)\leq 3Rn+3R\Area_{\mathcal A}(n).
\]
Taking the supremum over all loops $P$ of length at most $t$ proves
\[
\Area_Z(t)\leq 3R\bigl(\Area_{\mathcal A}(t)+t\bigr),
\]
as required.

% requires: \usetikzlibrary{calc,shapes.geometric}
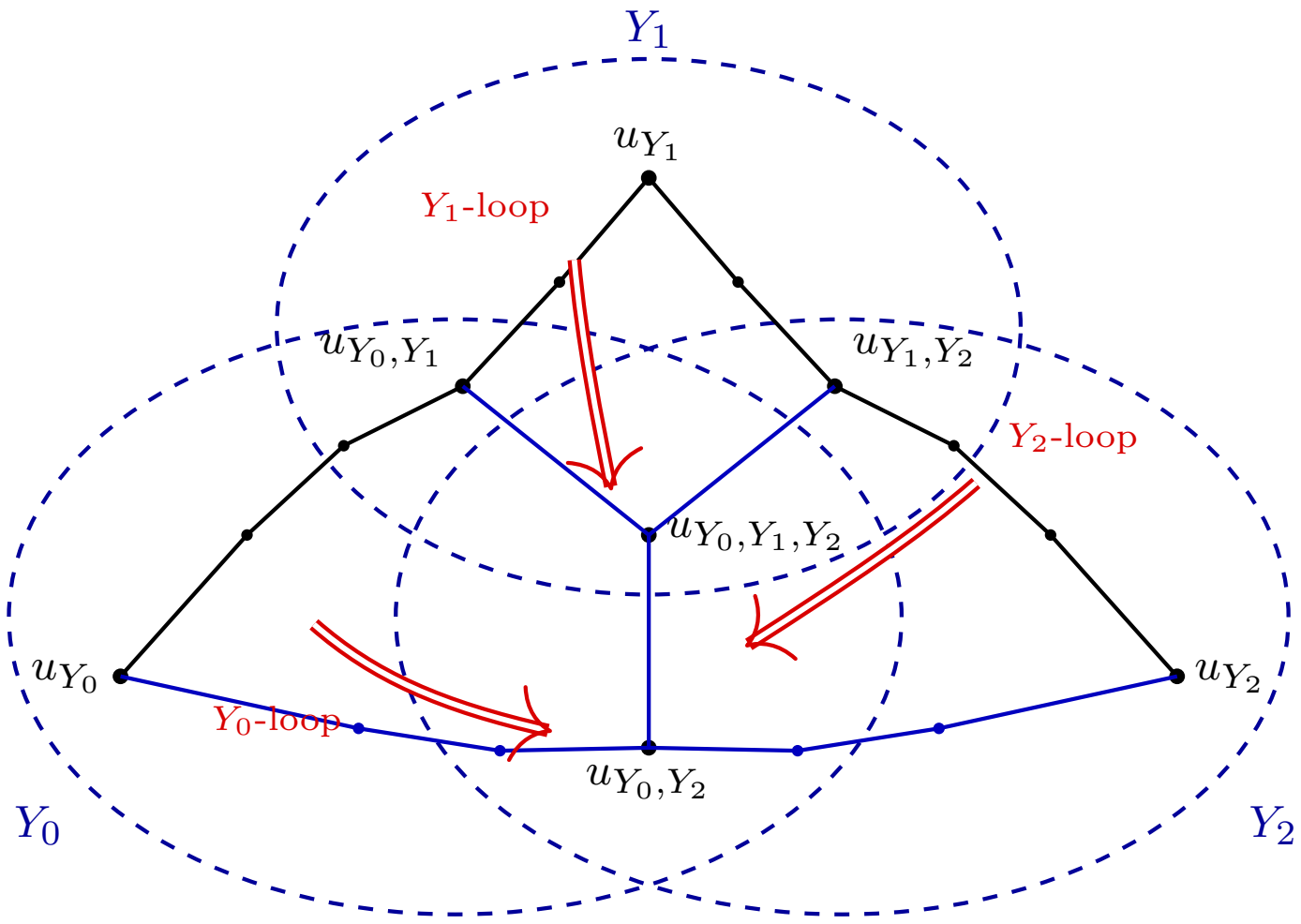
\begin{figure}[ht]
\centering
\begin{tikzpicture}[scale=1.05, every node/.style={font=\small}]

% colors
\colorlet{cover}{blue!60!black}
\colorlet{contract}{red!85!black}
\colorlet{oldpath}{black}
\colorlet{newpath}{blue!75!black}

% --- covering sets Y_0, Y_1, Y_2 ---
\draw[cover,dashed,thick] (3.3,1.6) ellipse (3.0 and 2.0);
\draw[cover,dashed,thick] (4.6,3.55) ellipse (2.5 and 1.8);
\draw[cover,dashed,thick] (5.9,1.6) ellipse (3.0 and 2.0);

\node[cover] at (0.5,0.2) {$Y_0$};
\node[cover] at (4.6,5.55) {$Y_1$};
\node[cover] at (8.8,0.2) {$Y_2$};

% --- main named vertices ---
\coordinate (u0)  at (1.05,1.2);
\coordinate (u01) at (3.35,3.15);
\coordinate (u1)  at (4.6,4.55);
\coordinate (u12) at (5.85,3.15);
\coordinate (u2)  at (8.15,1.2);
\coordinate (u02) at (4.6,0.72);
\coordinate (u012) at (4.6,2.15);

\fill (u0) circle (1.5pt);
\fill (u01) circle (1.5pt);
\fill (u1) circle (1.5pt);
\fill (u12) circle (1.5pt);
\fill (u2) circle (1.5pt);
\fill (u02) circle (1.5pt);
\fill (u012) circle (1.5pt);

\node[left] at (u0) {$u_{Y_0}$};
\node[above left] at (u01) {$u_{Y_0,Y_1}$};
\node[above] at (u1) {$u_{Y_1}$};
\node[above right] at (u12) {$u_{Y_1,Y_2}$};
\node[right] at (u2) {$u_{Y_2}$};
\node[below] at (u02) {$u_{Y_0,Y_2}$};
\node[right] at (u012) {$u_{Y_0,Y_1,Y_2}$};

% --- original canonical path u0 ... u01 ... u1 ... u12 ... u2 ---
\coordinate (a1) at (1.9,2.15);
\coordinate (a2) at (2.55,2.75);
\coordinate (b1) at (4.0,3.85);
\coordinate (b2) at (5.2,3.85);
\coordinate (c1) at (6.65,2.75);
\coordinate (c2) at (7.3,2.15);

\draw[oldpath,thick] (u0)--(a1)--(a2)--(u01)--(b1)--(u1)--(b2)--(u12)--(c1)--(c2)--(u2);

\foreach \p in {a1,a2,b1,b2,c1,c2}
  \fill (\p) circle (1.1pt);

% --- target canonical path u0 ... u02 ... u2 ---
\coordinate (d1) at (2.65,0.85);
\coordinate (d2) at (3.6,0.7);
\coordinate (d3) at (5.6,0.7);
\coordinate (d4) at (6.55,0.85);

\draw[newpath,thick] (u0)--(d1)--(d2)--(u02)--(d3)--(d4)--(u2);
\foreach \p in {d1,d2,d3,d4}
  \fill[newpath] (\p) circle (1.1pt);

% --- pairwise-intersection paths to the triple intersection ---
\draw[newpath,thick] (u01)--(u012);
\draw[newpath,thick] (u12)--(u012);
\draw[newpath,thick] (u02)--(u012);

% --- contraction arrows ---
\draw[contract,thick,double,->,double distance=1.1pt]
  (4.1,4.0) .. controls (4.15,3.45) and (4.25,2.95) .. (4.35,2.45);
\node[contract] at (3.5,4.35) {\scriptsize $Y_1$-loop};

\draw[contract,thick,double,->,double distance=1.1pt]
  (6.8,2.5) .. controls (6.4,2.15) and (5.95,1.85) .. (5.25,1.4);
\node[contract] at (7.45,2.8) {\scriptsize $Y_2$-loop};

\draw[contract,thick,double,->,double distance=1.1pt]
  (2.35,1.55) .. controls (2.75,1.2) and (3.2,1.0) .. (3.95,0.83);
\node[contract] at (2.1,0.9) {\scriptsize $Y_0$-loop};

\end{tikzpicture}
\caption{The triangle move on canonical lifts. The original canonical path from $u_{Y_0}$ to $u_{Y_2}$ passes through $u_{Y_0,Y_1}$, $u_{Y_1}$, and $u_{Y_1,Y_2}$. We connect each pairwise intersection point to the triple-intersection point $u_{Y_0,Y_1,Y_2}$. This yields three loops, contained respectively in $Y_1$, $Y_2$, and $Y_0$, which are contracted separately. The result is the canonical path from $u_{Y_0}$ to $u_{Y_2}$ through $u_{Y_0,Y_2}$.}
\label{fig:example-move}
\end{figure}

\end{proof}

We apply the two lemmas to the covers of the coset complex $Y$ and the auxiliary complex $Z$ described above.

\begin{corollary} \label{cor:coset-diameter-word-length}
Let $G$ be a finite group with subgroups $H_i$, $i\in I$, such that $Y=\CC(G;(H_i)_{i\in I})$ is connected. Let $L$ be the maximum, over $g\in G$, of the minimum number of letters from $\bigcup_i H_i$ whose product is $g$. Then
\[
L-1\leq\diam(Y)\leq L+1.
\]
In particular, $\diam(Y)\leq D$ implies that every $g\in G$ is a product of at most $D+1$ elements from the subgroups.
\end{corollary}

\begin{proof}
The Cayley graph underlying $Z$ has diameter $L$. The cover of $Z$ by its coset simplices has nerve $Y$, and each member has diameter at most $1$. Apply \pref{lem:nerve-diameter} with $D_1=1$.
\end{proof}

We next compare fillings of loops in $Z$ with the Dehn function associated with $Y$. We use the definition of the Dehn function given in \pref{subsec:dehn}, which counts applications of relators in the free group $\Free(S)$, where $S=\bigcup_i H_i\setminus\set{1}$.

There is a minor distinction between words in $\Free(S)$ and paths in $Z$. For $g\in S$, write $[g]$ for the corresponding generator of $\Free(S)$. The formal inverse $[g]^{-1}$ and the generator $[g^{-1}]$ indexed by the inverse group element are different symbols in $\Free(S)$. However, both represent the same element of $G$. In particular, the complete multiplication table of any subgroup containing $g$ includes the relation
\[
[g][g^{-1}]=1.
\]
Thus, a geometric backtracking move in $Z$ need not correspond to a free reduction in $\Free(S)$, but can be realized using one multiplication-table relation. We keep track of this distinction in the following lemma.

\begin{lemma} \label{lem:cayley-area-word-area}
Let $Y=\CC(G;(H_i)_{i\in I})$ be $1$-connected, and let $Z$ be the auxiliary complex defined above. With the complete multiplication-table presentation of $G$ and the Dehn function $\Dehn_Y$ defined in the preliminaries, for every $t\geq 1$,
\[
\Area_Z(t)\leq\Dehn_Y(t)\leq 3\Area_Z(t)+2t.
\]
\end{lemma}

\begin{proof}
An oriented edge $(x,xg)$ of $Z$ is labeled by an element $g\in S$. Thus a word in the alphabet $S\cup S^{-1}$ whose product is the identity defines a loop in $Z$, using the corresponding elements of $G$ as edge increments. A free reduction $[g][g]^{-1}$ becomes a backtracking move, while an application of a multiplication-table relation is realized by at most one triangle move in a coset simplex. In particular, a two-letter inverse relation is realized by a backtracking move and requires no triangle moves. Therefore, a reduction of a word gives a contraction of the corresponding loop in $Z$ with no more triangle moves than relation applications. Taking the positive edge-label word of an arbitrary loop in $Z$ proves $\Area_Z(t)\leq\Dehn_Y(t)$.

For the other inequality, consider a word $w\in\Free(S)$ of length $n\leq t$ representing the identity in $G$. First, replace each occurrence of a formal inverse $[g]^{-1}$ by the positive generator $[g^{-1}]$. This replacement uses one application of the relation $[g][g^{-1}]=1$. Thus, after at most $n$ relation applications, we obtain a positive word $w'$ of length $n$ whose letters are the edge labels of a loop $P$ in $Z$.

Let $a=\Area_Z(P)$. A contraction of $P$ with $a$ triangle moves gives a disk diagram consisting of $a$ triangles, possibly together with backtracking edges. We may remove redundant backtracking branches, so that the remaining diagram consists of the $a$ triangles and the edges needed to connect them to the original boundary. Each triangle contributes three edge occurrences, while the original boundary contributes $n$. Consequently, the diagram requires at most $(3a+n)/2$ identifications of oppositely oriented edge occurrences.

Each triangle has all its vertices in a common coset $xH_i$, so its oriented boundary labels form a three-letter multiplication-table relation in $H_i$. Thus the triangles contribute at most $a$ relation applications. Along an identification of oppositely oriented edges, the positive labels are $[g]$ and $[g^{-1}]$, which can be canceled using the relation $[g][g^{-1}]=1$. Hence these identifications contribute at most $(3a+n)/2$ additional relation applications.

Including the initial replacement of formal inverses, we obtain
\[
\Area_R(w)\leq n+a+\frac{3a+n}{2}\leq 3a+2n.
\]
Taking the maximum over all words of length at most $t$ proves
\[
\Dehn_Y(t)\leq 3\Area_Z(t)+2t.
\]
\end{proof}

\begin{corollary} \label{cor:coset-area-comparison}
Let $Y=\CC(G;(H_i)_{i\in I})$ be a $1$-connected coset complex with $|I|\geq 3$, and let $\Dehn_Y$ be as in \pref{def: dehn function of coset complex}. Then, for every $t\geq 1$,
\[
\Area_Y(t)\leq 15\bigl(\Dehn_Y(t)+t\bigr),
\]
and
\[
\Dehn_Y(t)\leq 45\Area_Y(t)+47t.
\]
\end{corollary}

\begin{proof}
In both covers of $Y$ and $Z$, all members and their nonempty intersections are simplices, so we may take $D_1=D_2=1$. Every loop of length at most $5$ in a simplex contracts with at most $5$ triangle moves. Applying \pref{lem:nerve-fillings} to the two covers gives
\[
\Area_Y(t)\leq 15\bigl(\Area_Z(t)+t\bigr),\qquad
\Area_Z(t)\leq 15\bigl(\Area_Y(t)+t\bigr).
\]
The first inequality of the corollary follows from $\Area_Z(t)\leq\Dehn_Y(t)$. For the second, we have
\[
\Dehn_Y(t)\leq 3\Area_Z(t)+2t\leq 45\Area_Y(t)+47t,
\]
as required.
\end{proof}

We conclude by translating the area bounds into bounds on cone radii, and conversely translating cone radii into bounds on the Dehn function. This forward direction of this theorem is the statement in \pref{thm: dehn implies cone}.

\begin{theorem} \label{thm:coset-cones-and-dehn}
Let $G$ be a finite group with subgroups $H_i$, $i\in I$, where $|I|\geq 3$, and assume that the coset complex $Y=\CC(G;(H_i)_{i\in I})$ is $1$-connected. Write $D=\diam(Y)$. Then $Y$ admits a $1$-cone $\Cone_Y$ satisfying
\[
\Rad_0(\Cone_Y)\leq D,\qquad
\Rad_1(\Cone_Y)\leq 15\bigl(\Dehn_Y(2D+1)+2D+1\bigr).
\]
Conversely, if $Y$ admits a $1$-cone $\Cone_Y$ with $\Rad_1(\Cone_Y)\leq A$, then for every $t\geq 1$,
\[
\Dehn_Y(t)\leq(45A+47)t.
\]
\end{theorem}

\begin{proof}
Choose an apex $v_0$ and shortest paths $P_u$ from $v_0$ to every vertex $u$, taking $P_{v_0}$ to be the trivial path. Then $\Rad_0(\Cone_Y)\leq D$. For each ordered edge $(u,v)$, the loop $P_u\circ(u;v)\circ P_v^{-1}$ has length at most $2D+1$. By \pref{cor:coset-area-comparison}, we can choose a contraction $T_{uv}$ with at most $15(\Dehn_Y(2D+1)+2D+1)$ triangle moves. These paths and contractions constitute a $1$-cone in the sense of \pref{def: non-abelian cone} and prove the first implication.

Conversely, let $P=(u_0;u_1;\dots;u_m=u_0)$ be any loop of length $m\leq t$. Insert the paths $P_{u_i}$ and their inverses between consecutive edges. Using only backtracking moves, we obtain a conjugate of the concatenation of the $m$ loops $P_{u_i}\circ(u_i;u_{i+1})\circ P_{u_{i+1}}^{-1}$. The cone contracts each of these using at most $A$ triangle moves, so $\Area_Y(P)\leq mA$ and hence $\Area_Y(t)\leq At$. Applying \pref{cor:coset-area-comparison} gives
\[
\Dehn_Y(t)\leq45\Area_Y(t)+47t\leq(45A+47)t,
\]
as required.
\end{proof}

\end{document}